\documentclass[a4paper,11pt,bibtotoc,tocindent,idxtotoc,makeidx]{article}

\usepackage[english]{babel}
\usepackage{amsmath}
\usepackage{amsthm}
\usepackage{amssymb}
\usepackage{amsfonts}
\usepackage{amsxtra}
\usepackage{color}
\usepackage[breaklinks]{hyperref}
\usepackage{enumitem}
\usepackage{tocloft}
\usepackage{makeidx}
\usepackage{textcomp}
\usepackage{rotating}
\usepackage{multirow}
\usepackage{booktabs,longtable}
\usepackage[all]{xy}
\usepackage{float}
\usepackage{rotfloat}

\setenumerate{itemsep=0pt}
\setitemize{itemsep=0pt}
\setdescription{itemsep=0pt}
\setlist{noitemsep}

\numberwithin{equation}{section}
\newcommand{\td}{\,\mathrm{d}}

\newcommand{\Stab}{\textup{Stab}}
\newcommand{\const}{\textup{const}}
\newcommand{\linspan}{\textup{span}}

\newcommand{\Ad}{\textup{Ad}}
\newcommand{\ad}{\textup{ad}}
\renewcommand{\Re}{\textup{Re}}
\renewcommand{\Im}{\textup{Im}}
\renewcommand{\det}{\textup{det}}
\newcommand{\id}{\textup{id}}
\newcommand{\tr}{\textup{tr}}

\newcommand{\Tr}{\textup{Tr}}
\newcommand{\Sym}{\textup{Sym}}

\newcommand{\Herm}{\textup{Herm}}
\newcommand{\Str}{\textup{Str}}
\newcommand{\str}{\mathfrak{str}}

\newcommand{\aut}{\mathfrak{aut}}

\newcommand{\Co}{\textup{Co}}
\newcommand{\co}{\mathfrak{co}}

\newcommand{\gl}{\mathfrak{gl}}
\newcommand{\SL}{\textup{SL}}
\renewcommand{\sl}{\mathfrak{sl}}

\newcommand{\Sp}{\textup{Sp}}
\renewcommand{\sp}{\mathfrak{sp}}
\newcommand{\Mp}{\textup{Mp}}

\newcommand{\SO}{\textup{SO}}
\newcommand{\so}{\mathfrak{so}}
\newcommand{\SU}{\textup{SU}}
\newcommand{\su}{\mathfrak{su}}

\newcommand{\RR}{\mathbb{R}}
\newcommand{\KK}{\mathbb{K}}
\newcommand{\CC}{\mathbb{C}}
\newcommand{\ZZ}{\mathbb{Z}}
\newcommand{\NN}{\mathbb{N}}
\newcommand{\QQ}{\mathbb{Q}}
\renewcommand{\SS}{\mathbb{S}}
\newcommand{\HH}{\mathbb{H}}
\newcommand{\FF}{\mathbb{F}}
\newcommand{\OO}{\mathbb{O}}

\newcommand{\XX}{\mathbb{X}}
\newcommand{\BB}{\mathbb{B}}
\newcommand{\PP}{\mathbb{P}}

\renewcommand{\1}{\mathbf{1}}

\newcommand{\calF}{\mathcal{F}}
\newcommand{\calO}{\mathcal{O}}

\newcommand{\calW}{\mathcal{W}}
\newcommand{\calB}{\mathcal{B}}
\newcommand{\calH}{\mathcal{H}}
\newcommand{\calP}{\mathcal{P}}
\newcommand{\calU}{\mathcal{U}}

\newcommand{\calD}{\mathcal{D}}
\newcommand{\calT}{\mathcal{T}}
\newcommand{\calV}{\mathcal{V}}

\newcommand{\calR}{\mathcal{R}}
\newcommand{\calE}{\mathcal{E}}

\newcommand{\frakg}{\mathfrak{g}}

\newcommand{\frakk}{\mathfrak{k}}
\newcommand{\frakp}{\mathfrak{p}}
\newcommand{\fraks}{\mathfrak{s}}
\newcommand{\frakn}{\mathfrak{n}}

\newcommand{\frakm}{\mathfrak{m}}
\newcommand{\frakl}{\mathfrak{l}}
\newcommand{\frakt}{\mathfrak{t}}

\newcommand{\frake}{\mathfrak{e}}
\newcommand{\frakf}{\mathfrak{f}}

\newcommand{\diag}{\textup{diag}}
\newcommand{\calBe}{\mathcal{B}_{\bf e}}
\renewcommand{\min}{\textup{min}}
\newcommand{\Spin}{\textup{Spin}}

\newcommand{\minholrep}{minimal representation}
\newcommand{\minholreps}{minimal representations}

\theoremstyle{plain}
\newtheorem{theorem}{Theorem}[section]
\newtheorem{thmalph}{Theorem}

\newtheorem{proposition}[theorem]{Proposition}
\newtheorem{lemma}[theorem]{Lemma}
\newtheorem{corollary}[theorem]{Corollary}

\newtheorem{fact}[theorem]{Fact}

\theoremstyle{definition}
\newtheorem{definition}[theorem]{Definition}
\newtheorem{example}[theorem]{Example}

\newtheorem{remark}[theorem]{Remark}

\theoremstyle{remark}
\renewenvironment{pf}{\begin{proof}[\sc Proof]}{\end{proof}}

\numberwithin{equation}{section}
\renewcommand{\theequation}{\arabic{section}.\arabic{equation}} 

\begin{document}

\title{Fock model and Segal--Bargmann transform for minimal representations of Hermitian Lie groups}
\author{Joachim Hilgert, Toshiyuki Kobayashi\footnote{Partially supported by Grant-in-Aid for Scientific Research (B) (22340026), Japan Society for the Promotion of Science, and the Alexander Humboldt Foundation.}, Jan M\"ollers, Bent \O rsted}
\maketitle
\begin{abstract}

For any Hermitian Lie group $G$ of tube type
 we construct a Fock model of its \minholrep. 
The Fock space is defined on the minimal nilpotent $K_\CC$-orbit $\XX$ in $\frakp_\CC$ and the $L^2$-inner product involves a K-Bessel function as density. Here $K\subseteq G$ is a maximal compact subgroup
 and $\frakg_{\CC}=\frakk_{\CC}+\frakp_{\CC}$
 is a complexified Cartan decomposition. 
In this realization the space
 of $\frakk$-finite vectors consists
 of holomorphic polynomials on $\XX$. 
The reproducing kernel of the Fock space is calculated explicitly in terms of an I-Bessel function. 
We further find 
 an explicit formula
 of a generalized Segal--Bargmann transform
 which intertwines the Schr\"{o}dinger and Fock model. 
Its kernel involves the same I-Bessel function. 
Using the Segal--Bargmann transform
 we also determine the integral kernel
 of the unitary inversion operator
 in the Schr{\"o}dinger model
 which is given by a J-Bessel function.\\



\end{abstract}

\newpage

\tableofcontents

\newpage

\addcontentsline{toc}{section}{Introduction}
\section*{Introduction}

The classical Segal--Bargmann transform is the integral operator
\begin{align*}
 \BB u(z) &= e^{-\frac{1}{2}z^2}
\int_{\RR^n}{e^{2z\cdot x}e^{-x^2}u(x)\td x}.
\end{align*}
It induces a unitary isomorphism
$$     
     \BB:L^2(\RR^n) \overset {\sim}\to\calF(\CC^n),
$$
where $\calF(\CC^n)$ denotes
 the classical Fock space on $\CC^n$ consisting of entire functions,
 square integrable with respect to
 the Gaussian measure $e^{-|z|^2}\td z$.  
This transform has many remarkable properties, 
e.g., 
 it intertwines the harmonic oscillator with the Euler operator, 
and it has been widely used 
 in physics problems
 such as field theory. 
For a detailed introduction to the classical Segal--Bargmann transform and the classical Fock space we refer the reader to the book of G. B. Folland \cite{Fol89}.

In representation theory
 the unitary operator $\BB$ intertwines 
 two prominent models of the same unitary representation
 of the metaplectic group  $\Mp(n,\RR)$,
 a double cover of the symplectic group, 
 namely,
 the Schr\"odinger and the Fock model
 of the Weil representation,
 which is also referred to
 as the harmonic--Segal--Shale--Weil--oscillator--metaplectic
  representation. 

We highlight the fact 
 that the Weil representation consists
 of two minimal representations
 (see Definition \ref{def:minrep})
 of the simple Lie group
 $\Mp(n,\RR)$.
The aim of this article is 
 to construct complete analogues
 of the above-mentioned theory 
 in the generality that the Weil representation $\varpi$ is replaced by
 a minimal representation
 of an arbitrary Hermitian Lie group $G$
 of tube type.  
This includes the construction of
\begin{itemize}
\item the `Schr\"{o}dinger model' of minimal representations,
\item the `Fock model' of minimal representations, and
\item the `Segal--Bargmann transform' intertwining them.
\end{itemize}
The \lq{Schr\"odinger model}\rq\
 of the corresponding unitary representations
 was constructed earlier 
 in this setting
 (see Vergne--Rossi \cite{RV76}), 
and has been extensively studied in a more general setting
(see e.g.\ \cite{HKM11,KM07a,KMa11,KO03c}).

In order to construct the latter two objects,
 we recall the Kirillov--Kostant--Duflo orbit philosophy,
which suggests to understand minimal representations in relation with
 {\it{real}} minimal nilpotent coadjoint orbits
 $\mathbb{O}_{\min}^G$
 in a functorial manner. 
In fact,
 our key idea
 that underlies the construction
 of the `Fock model' and the \lq{Segal--Bargmann transform}\rq\
 is to define a {\it{geometric quantization
 of the Kostant--Sekiguchi correspondence}}
 \cite{Seki87}
 of minimal nilpotent orbits
$\mathbb{O}_{\min}^G
 \leftrightarrow \mathbb{O}_{\min}^{K_{\mathbb{C}}}$,
which is summarized in the diagram below.  

$$
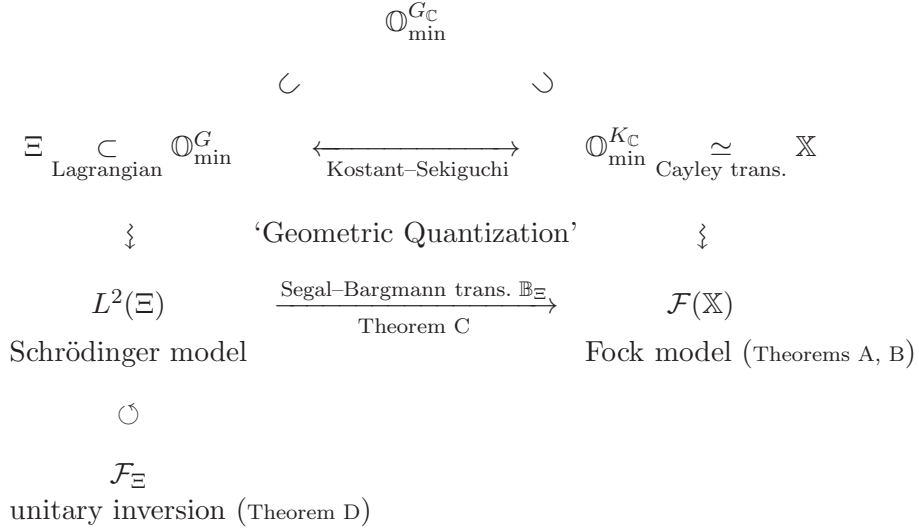
\begin{matrix}
& 
\mathbb{O}_{\min}^{G_{\mathbb{C}}}
\\[2ex]
& \rotatebox[origin=c]{45}{$\subset$}
    \hfill
  \rotatebox[origin=c]{135}{$\subset$}
\\[2ex]
\Xi \underset{\text{Lagrangian}}{\subset} \mathbb{O}_{\min}^{G}
& \rlap{$\leftarrow$}\xrightarrow[\text{Kostant--Sekiguchi}]{}
& \mathbb{O}_{\min}^{K_{\mathbb{C}}} 
\underset{\text{Cayley trans.}}{\simeq} \mathbb{X}
\\[4ex]
\rotatebox[origin=c]{-90}{$\rightsquigarrow$}
&\hidewidth\text{`Geometric Quantization'}\hidewidth
& \rotatebox[origin=c]{-90}{$\rightsquigarrow$}
\\[2ex]
L^2(\Xi)
&\xrightarrow[\text{Theorem C}]{\text{Segal--Bargmann trans.\ $\mathbb{B}_{\Xi}$}}
&\mathcal{F}(\mathbb{X})
\\
\text{Schr\"{o}dinger model}
&& \text{Fock model (}{\scriptstyle\text{Theorems A, B}}\text{)}\hidewidth
\\[2ex]
\circlearrowleft
\\[2ex]
\mathcal{F}_{\Xi}
\\
\text{unitary inversion (}{\scriptstyle\text{Theorem D}}{)}\hidewidth
\end{matrix}
$$

The techniques
 of our proofs are 
 twofold:
discrete branching laws
 of minimal representations
 with respect to a distinguished subalgebra
${\mathfrak{s}} \simeq {\mathfrak{sl}}(2,\RR)$
 (see \eqref{eqn:sl2})
 and the theory
 of Jordan algebras.
We shall see 
 for Hermitian groups of tube type
 how the Jordan algebra structure allows generalizations of many aspects of the classical case.\\

Let us explain our results 
more precisely.  
Suppose that $V$ is 
 a simple Euclidean Jordan algebra.  
We denote by 
 $\Co(V)$ and $\operatorname{Str}(V)$
 the conformal group
 and the structure group
 of $V$, 
 respectively. 
We set $G:=\Co(V)_0$
 and $L:=\operatorname{Str}(V)_0$
the identity component groups.  
We denote by $\vartheta$ the Cartan involution of $G$ given by conjugation with the conformal inversion $j(x)=-x^{-1}$ and let $K=G^\vartheta$ be the corresponding maximal compact subgroup of $G$. 
Then $G$ is the group of biholomorphic transformations
 on an irreducible Hermitian symmetric space $G/K$
 of tube type.
Conversely,
 any simple Hermitian Lie group
 of tube type with trivial center
 arises in this fashion
 (see Table \ref{tb:Groups}
 in the Appendix).

The Lie algebra $\frakg$ of $G$,
also known as the Kantor--Koecher--Tits algebra,
has a Gelfand--Naimark decomposition $\frakg=\frakn+\frakl+\overline{\frakn}$, 
where the abelian Lie algebra
 $\frakn\simeq V$ acts on $V$
 by constant vector fields, 
$\frakl:=\str(V)\subseteq\gl(V)$ is the structure algebra acting
 by linear vector fields,
 and $\overline{\frakn}=\vartheta\frakn$ acts
 by quadratic vector fields. 

Let $G_{\mathbb{C}}$ be the complexification of $G$, 
 and $K_{\mathbb{C}}
$ that of $K$.  
There is a unique minimal nilpotent coadjoint orbit
 of $G_\CC$, 
 which we denote by $\mathbb{O}_\min^{G_\CC}\subseteq\frakg_\CC^*$.  
We write 
$
   {\mathfrak {g}}_{\mathbb{C}}
  ={\mathfrak {k}}_{\mathbb{C}}+{\mathfrak {p}}_{\mathbb{C}}
$
 for the complexified Cartan decomposition,
 and 
$
   {\mathfrak {g}}_{\mathbb{C}}^{\ast}
  ={\mathfrak {k}}_{\mathbb{C}}^{\ast}
  +{\mathfrak {p}}_{\mathbb{C}}^{\ast}
$
 for the dual.  
Then $\mathbb{O}_\min^{G_\CC} \cap {\mathfrak {p}}_\CC^*$
 splits into two equi-dimensional $K_\CC$-orbits.  
We choose one of them,
 and denote it by $\mathbb{O}_\min^{K_\CC}\subseteq\frakp_\CC^*$.  
The intersection 
$\mathbb{O}_{\min}^{G_{\mathbb{C}}} \cap \mathfrak{g}^*$
also splits into two equi-dimensional $G$-orbits.
We write $\mathbb{O}_{\min}^G$ for the component corresponding to
$\mathbb{O}_{\min}^{K_{\mathbb{C}}}$
via the Kostant--Sekiguchi correspondence.

In what follows 
 we denote by $\widetilde {I}_{\alpha}(z)$, 
 $\widetilde {J}_{\alpha}(z)$, and $\widetilde {K}_{\alpha}(z)$,
 by the renormalized Bessel functions
 (see Appendix \ref{app:BesselFcts}).  

\subsubsection*{The Schr{\"o}dinger model}
Let us recall some known results on the $L^2$-model
 (the \lq{Schr{\"o}dinger model}\rq)
 for \minholreps\
 (see Subsection \ref{sec:SchrödingerForMinimalHolomorphic}
 for details). 
For the statement of the results we assume that $\dim\,V>1$. 
The case $V=\RR$ is discussed 
 at the end of this introduction. 

We identify $\frakg^{\ast}$ with $\frakg$
 by the Killing form. 
Then the intersection $\Xi:=\mathbb{O}_{\operatorname{min}}^{G} \cap
{\mathfrak{n}}$
 is a Lagrangian submanifold
 of $\mathbb{O}_{\operatorname{min}}^{G}$
 endowed with the Kostant--Kirillov--Souriau
 symplectic form (cf.~\cite[Theorem 2.9]{HKM11}).  

Let $L^2(\Xi, \td \mu)$
 be the Hilbert space
 consisting of square integrable functions
 on $\Xi$
 with respect to a unique
 (up to scalar multiples)
 $L$-equivariant Radon measure
 $\td \mu$ on $\Xi$.  
Then the natural action of
 $L \ltimes \exp({\mathfrak{n}})$
 extends to an irreducible unitary representation, 
 to be denoted by $\pi$, 
of a finite cover $G^\vee$ of $G$.  
The resulting representation 
 is a minimal representation
 unless $\frakg_{\CC}$ is of type $A$.  
The corresponding differential representation of $\frakg$
 is given by differential operators up to order $2$. 
Its underlying $(\frakg,\frakk)$-module
 is a lowest weight module
 of scalar type
 whose parameter $\lambda$ is the smallest non-zero discrete point of the Wallach set
 (see \eqref{eqn:Wal1}).  
The corresponding one-dimensional 
 minimal $\frakk$-type 
 is given by $\CC\psi_0$ in this $L^2$-model,
 where $\psi_0$ is defined by
\[
     \psi_0(x)=e^{-\tr(x)}.  
\]
Here, 
 $\tr(-)$ denotes the trace function of the real Jordan algebra.  
We further remark that on $\Xi$ we have $\tr(x)=|x|$, where $|x|=(x|x)^{\frac{1}{2}}$ is the norm on $V$ induced by the trace form $(-|-)$. The trace form is extended $\CC$-bilinearly to $V_\CC$.

\subsection*{The Fock space}
Let us explain the construction
 of a new model
 for the same representation
 on a space of holomorphic $L^2$-functions which resembles the classical Fock space. The geometry for this space is given by the complexification $\XX$ in $V_\CC$ of the orbit $\Xi$, 
 on which the complexified structure group $L_\CC$
 acts transitively.  
Up to scalar multiples there is a unique $L_\CC$-equivariant measure $\td\nu$ on $\XX$. We further define a density
\begin{align*}
 \omega(z) &= \widetilde{K}_{\lambda-1}(|z|), & z\in\XX,
\end{align*}
in terms of the renormalized K-Bessel function.  
Here,
 $|z|=(z|\overline{z})^{\frac{1}{2}}$
 and $\lambda$ is given by \eqref{eqn:Wal1}
 as explained above. 
Denoting by $\calO(\XX)$ the space of holomorphic functions on the complex manifold $\XX$, we then define a \lq{Fock space}\rq\ by 
\begin{align}
\label{eqn:Fintro}
 \calF(\XX) &= \left\{F\in\calO(\XX):\int_\XX{|F(z)|^2\omega(z)\td\nu(z)}<\infty\right\}
\end{align}
endowed with a pre-Hilbert structure 
 by 
\begin{align}
\label{eqn:Finner}
 \langle F,G\rangle &:= \int_\XX{F(z)\overline{G(z)}\omega(z)\td\nu(z)}, & F,G\in\calF(\XX).
\end{align}

We then have:
\begin{thmalph}
[Theorems \ref{thm:Frepro}, \ref{thm:NaturalFockSpace},
 and \ref{thm:FockAsRestriction}]
\label{thm:a}
\begin{enumerate}
\item[\textup{(1)}] The Fock space $\calF(\XX)$ is a Hilbert space.
\item[\textup{(2)}] 
The reproducing kernel of $\calF(\XX)$ is given by
\begin{align*}
 \KK(z,w) &= \Gamma(\lambda)\widetilde{I}_{\lambda-1}(\sqrt{(z|\overline{w})}).
\end{align*}
\item[\textup{(3)}] 
Every function in $\calF(\XX)$ can be extended to an entire holomorphic function on the ambient space $V_\CC$.  
Further,
the space $\calP(\XX)$ of restrictions of holomorphic polynomials on $V_\CC$ to $\XX$ is dense in $\calF(\XX)$. 
\end{enumerate}
\end{thmalph}

Note that the renormalized I-Bessel function
 $\widetilde{I}_\alpha(t)$
 is an even function
 and hence $\widetilde{I}_\alpha(\sqrt{t})$
 is an entire function on $\CC$.\\

Since the Cayley transform (see \eqref{eq:IsoKCLC})
 induces an isomorphism of complex groups
 $c: K_{\CC} \overset \sim \to L_{\CC}$
 and a biholomorphic map 
 $\mathbb{O}_\min^{K_\CC} \overset \sim \to \XX$,
 the right-hand side of \eqref{eqn:Fintro}
 gives an {\it{intrinsic}} definition of the Fock space
 built on the minimal $K_{\CC}$-nilpotent
 orbit $\mathbb{O}_\min^{K_\CC}$,
 namely, 
the space of holomorphic,
 square integrable functions
 on $\mathbb{O}_\min^{K_\CC}$
 against the measure $\omega \td \nu$.  

Theorem \ref{thm:a} (1) and (3) assert
 that the intrinsic definition \eqref{eqn:Fintro}
 of the Fock space $\calF(\XX)$
 coincides with the extrinsic definition
 built on the embedding $\XX \subset V_{\CC}$.  
This feature is noteworthy
 even in the classical Fock model
 of the Weil representation,
 where $V_{\CC}=\operatorname{Sym}(n,\CC)$
 and $\XX$ is a submanifold
 consisting of rank one matrices
 in $\operatorname{Sym}(n,\CC)$;
 Theorem \ref{thm:a} (3) says
 that any holomorphic,
 square integrable function
 on the $n$-dimensional complex submanifold $\XX$
 extends holomorphically
 to the $\frac 1 2 n(n+1)$-dimensional space
 $\operatorname{Sym}(n,\CC)$.  

\subsection*{Fock model of minimal representations}

A remarkable property
 of our Fock space
 ${\mathcal{F}}({\mathbb{X}})$
 is that the conformal group $G$ 
 or its covering group acts on ${\mathcal{F}}({\mathbb{X}})$
 as an irreducible unitary representation.  
To construct the action $\rho$,
we begin with a \lq{holomorphic continuation
 $\pi_{\CC}$}\rq\ of the Schr{\"o}dinger model
 $(\pi,L^2(\Xi, \td \mu))$.  
The differential representation $\td \pi_{\CC}$
 of the Lie algebra $\frakg$
 is well-defined as a representation 
 on the space $\calP(\XX)$
 of regular functions,
 but unfortunately,
 $\calP(\XX)$ does not contain non-zero
 $\frakk$-finite vectors.  
Our idea is 
 to define the infinitesimal representation 
by 
\[
     \td \rho := \td \pi_{\CC} \circ c,
\]
 where $c \in \operatorname{Int}(\frakg_{\CC})$
 is a Cayley-type transform
 (see \eqref{eq:DefCayleyTransform}).
The resulting action $\td \rho(Y)$
 ($Y \in \frakg$) on $\calP(\XX)$ is still 
 a differential operator 
 up to order two.  
We then obtain:

\begin{thmalph}
[Fock model]
\label{thm:b}
\begin{enumerate}
\item[{\upshape{(1)}}]
The representation $(\td\rho,\calP(\XX))$
 is an irreducible ${\mathfrak{g}}$-module
 such that $\rho(Y)$ is skew-Hermitian 
 with respect to the $L^2$-inner product
 \eqref{eqn:Finner}
 for any $Y \in {\mathfrak {g}}$.  
\item[{\upshape{(2)}}]
The Lie algebra ${\mathfrak {k}}$ acts
 locally finitely, 
 and we have the following $\frakk$-type decomposition
\begin{align*}
 \calP(\XX) &= \bigoplus_{m=0}^\infty{\calP^m(\XX)},
\end{align*}
where $\calP^m(\XX)$ denotes the subspace of homogeneous polynomials of degree $m$. 
\item[{\upshape{(3)}}]
The $(\frakg,\frakk)$-module integrates to an irreducible unitary representation $\rho$ of the finite cover $G^\vee$ of $G$ on $\calF(\XX)$.
\end{enumerate}
\end{thmalph}

We give a direct proof for the irreducibility and unitarizability
 of $(\td \rho,\calP(\XX))$
 in Propositions \ref{prop:irredPX} and \ref{prop:PXinfuni}
 by using an explicit formula 
 of $\td\rho$,
 for which the crucial part is given by means of 
 the {\it{Bessel operators}}
 in the Jordan algebra
 as introduced by H. Dib \cite{Dib90}
 (see also \cite{HKM11}).   
The minimal cover $G^\vee$
 to which the representation integrates was 
 determined in \cite[Theorem 2.30]{HKM11}.\\

\subsection*{Segal--Bargmann transform}

The two irreducible unitary representations
 of the group $G^{\vee}$, 
 $(\pi,L^2(\Xi,\td \mu))$
 (the Schr{\"o}dinger model)
 and $(\rho,\calF(\XX))$ (the Fock model)
 are isomorphic to each other.  
We prove this by constructing 
 an explicit intertwining operator
 as follows (see Theorem \ref{thm:SBunitary}):

\begin{thmalph}[Segal--Bargmann transform]
\label{thm:c}
For $f\in L^2(\Xi,\td\mu)$ the integral
\begin{align*}
 (\BB_\Xi f) (z) &:= \Gamma(\lambda)e^{-\frac{1}{2}\tr(z)}\int_\Xi{\widetilde{I}_{\lambda-1}(2\sqrt{(z|x)})e^{-\tr(x)}f(x)\td\mu(x)}, & z\in V_\CC,
\end{align*}
converges uniformly on compact subsets in $V_\CC$
 and defines a holomorphic function
 $\BB_\Xi f\in\calF(\XX)$. This gives a unitary isomorphism
\begin{align*}
 \BB_\Xi:L^2(\Xi,\td\mu)\overset \sim \to\calF(\XX)
\end{align*}
intertwining the representations $\pi$ and $\rho$.
\end{thmalph}

 Here are some remarks on related works:

Once 
the Schr{\"o}dinger model $L^2(\Xi,\td \mu)$
 and the Fock model $\calF(\XX)$ are properly constructed,
 there is an alternative method
 to obtain the Segal--Bargmann transform $\BB_\Xi$.
That is, $\mathbb{B}_\Xi$ appears as the unitary part in the
 polar decomposition 
$\calR_\Xi^*=\BB_\Xi\circ\sqrt{\calR_\Xi\calR_\Xi^*}$,
where $\calR_\Xi^*$ is the adjoint of the unbounded operator
 $\calR_\Xi:\calF(\XX)\to L^2(\Xi,\td\mu)$,
$F(x) \mapsto e^{-\frac{1}{2}\tr(x)}F(x)$.
This method has been used before to construct analogues of the
 classical Segal--Bargmann transform or the Laplace transform
(see e.g.\ \cite{ADO07,HZ09,OO96}),
and we revisit this point in Subsection \ref{subsec:5.2}
 where we also find the heat kernel.

In works of Brylinski--Kostant \cite{BK94}
 and Achab--Faraut \cite{AF11} models
 of minimal representations on spaces
 of holomorphic functions were constructed
 by using the work of Kronheimer and Vergne.   
Aside from the fact that their cases
(non-Hermitian Lie groups) and our cases (Hermitian of tube type) are
 disjoint, 
there are one basic common point (0)
 and two major differences (1),
 (2) between their construction and our construction of the Fock model
(up to a Cayley transform):
\begin{itemize}
\item[(0)]
($K$-structure)\enspace
In both models, the space of regular functions on the minimal
nilpotent $K_{\mathbb{C}}$-orbit $\mathbb{O}_{\min}^{K_{\mathbb{C}}}$
is the space of $\mathfrak{k}$-finite vectors
 of the minimal representation;
\item[(1)]
(Lie algebra action)\enspace
The action of $\mathfrak{p}_{\mathbb{C}}$ on their Fock-type model is
given via pseudodifferential operators,
but our action is given by differential operators up to order two;
\item[(2)]
(measure on $\mathbb{O}_{\min}^{K_{\mathbb{C}}}$)\enspace
Their measure on $\mathbb{O}_{\min}^{K_{\mathbb{C}}}$ defining the
$L^2$-structure is not positive,
whereas our measure is given by a $K$-Bessel function and positive.  
\end{itemize}
In fact these two major points have enabled us to construct
explicitly an intertwiner between the $L^2$-model and the holomorphic
model,
whereas an analogous explicit operator is not known for their model.

\subsection*{Unitary inversion operator $\calF_{\Xi}$}
In the Schr{\"o}dinger model $L^2(\Xi,\td\mu)$,
 there is a distinguished unitary operator 
 $\calF_{\Xi}$,
 that is,
 the {\it{unitary inversion operator}}
 (see Section \ref{sec:UnitaryInversionOperator} for the precise definition). 
{}From the representation-theoretic viewpoint,
 the operator $\calF_{\Xi}$ generates
 the action $\pi$
 on  $L^2(\Xi,\td\mu)$
 of the whole group $G^\vee$ 
 together with the relatively simple action
 of a maximal parabolic subgroup.  

As a corollary of Theorem \ref{thm:c}, 
 we find an explicit integral kernel
 for the unitary inversion operator
 $\calF_{\Xi}$ in Theorem \ref{thm:UnitaryInversionKernel}
 as follows.  
\begin{thmalph}
[unitary inversion operator]
\label{thm:D}
The operator $\calF_\Xi$ is given by
\begin{align*}
 \calF_\Xi u(x) &= 2^{-r\lambda}\Gamma(\lambda)\int_\Xi{\widetilde{J}_{\lambda-1}(2\sqrt{(x|y)})u(y)\td\mu(y)}.  
\end{align*}
\end{thmalph}
Theorem \ref{thm:D}
 was established earlier 
 in the case of $\frakg=\so(2,n)$
 by two different methods
 by Kobayashi--Mano
 (\cite[Theorem 6.1.1]{KM07a}
 and \cite[Theorem 5.1.1]{KMa11}).  
Our approach 
 using the Segal--Bargmann transform
 gives yet another new approach
 even in the case $\frakg={\frak{so}}(2,n)$
 (see Section \ref{sec:exso2n}).

\subsection*{A theory of spherical harmonics}
A crucial role in the study of the structure of our representations is played
 by a subalgebra
 $\fraks \simeq {\mathfrak {sl}}(2,\RR)$
 spanned by a specific ${\mathfrak{sl}}_2$-triple $(E,F,H)$ in $\frakg$.  
A distinguished property
 is that the representation $\td \rho$ of $\frakg$
 is discretely decomposable
 in the sense of \cite{K98}
 when we restrict it 
 to the subalgebra $\fraks$.

In the Schr{\"o}dinger realization, 
 $\calBe:=-\sqrt{-1}\td\pi(F)$ is an elliptic second order differential operator on $\Xi$
 which extends to a self-adjoint operator
 on $L^2(\Xi,\td\mu)$, 
 and further to a holomorphic differential operator on $\XX$. We define an analogue of the space of spherical harmonics of degree $m$ by
\begin{align*}
 \calH^m(\XX) &:= \{p|_\XX:p\in\calP^m(V_\CC),\calBe p=0\}.
\end{align*}
Then $\calH^m(\XX)$ is naturally acted upon by the centralizer $Z_{G^\vee}(\fraks)$ of $\fraks$ in $G^\vee$ which turns out to be compact. This action is irreducible and we give explicit formulas for its highest weight vector and its spherical vector (see Propositions \ref{prop:HmSphericalVector} and \ref{prop:HmHighestWeightVector}). We further define $\calH^m(\Xi)$ to be the restriction to $\Xi\subseteq\XX$ of all elements of $\calH^m(\XX)$. The orbit $\Xi$ has a polar decomposition
\begin{align*}
 \Xi &= \RR_+ \times \SS, & \SS &= \{x\in\Xi:|x|=1\}.
\end{align*}
Since polynomials in $\calH^m(\Xi)$ are homogeneous,
 they are already uniquely determined
 by their values on $\SS$
 and we let $\calH^m(\SS)$ be
 the space of all restrictions
 to $\SS$ of polynomials
 in $\calH^m(\Xi)$. 
Then the embeddings $\SS\subseteq\Xi\subseteq\XX$
 induce isomorphisms 
$
     \calH^m(\SS)\simeq\calH^m(\Xi)\simeq\calH^m(\XX)
$
 of $Z_{G^\vee}(\fraks)$-representations. 

Let $\frakk^\frakl$ be the Lie algebra of $K^L:=L\cap K$ which equals the identity component of $Z_{G^\vee}(\fraks)$. 
Then $(\fraks,\frakk^\frakl)$ forms a reductive dual pair (see Proposition \ref{prop:DualPair}). Using the harmonics $\calH^m(\XX)$ resp. $\calH^m(\SS)$ we obtain an explicit branching law
 of our minimal representation with respect to the dual pair $(\fraks,\frakk^\frakl)$, generalizing an earlier result in \cite[Section 1.3]{KM07a}:

\begin{thmalph}
\label{thm:e}
\begin{enumerate}
\item[\textup{(1)}] \textup{(Theorem \ref{thm:DualPairDecompositionFock})}
The Fock space decomposes as a multiplicity-free direct
sum of irreducible $(\fraks,Z_{G^\vee}(\fraks))$-bimodules
\begin{align*}
 \calP(\XX) &\simeq \bigoplus_{m=0}^\infty{\calW_{r\lambda+2m}\boxtimes\calH^m(\XX)},
\end{align*}
where $\fraks$ acts on the first tensor factor and $ Z_{G^\vee}(\fraks)$ on the second. Here $\calW_{r\lambda+2m}=\linspan\{\tr^k(z):k\in\NN\}$ and $\fraks\simeq\sl(2,\RR)$ acts irreducibly on $\calW_{r\lambda+2m}$ with lowest weight $r\lambda+2m$, the integer $r$ being the rank of $G/K$.
\item[\textup{(2)}] \textup{($L^2$-dual pair correspondence,
 see Theorem \ref{thm:sl2decoS})}
The Schr\"odinger model $L^2(\Xi,\td\mu)$ decomposes discretely into a multiplicity-free sum of irreducible unitary representations of $\widetilde{\SL(2,\RR)}\times Z_{G^\vee}(\fraks)$
\begin{align*}
 L^2(\Xi,\td\mu) &\simeq {\sum_{a=0}^\infty}\raisebox{0.15cm}{$^\oplus$}{\,L^2(\RR_+,x^{r\lambda-1}\td x)\boxtimes\calH^m(\SS)},
\end{align*}
with respect to the natural homomorphism $\widetilde{\SL(2,\RR)}\times Z_{G^\vee}(\fraks)\to G^\vee$. Here $\widetilde{\SL(2,\RR)}$ acts irreducibly on the first factor of each summand $L^2(\RR_+,x^{r\lambda-1}\td r)\boxtimes\calH^m(\SS)$ with lowest weight $r\lambda+2m$.
\end{enumerate}
\end{thmalph}

We may think of Theorem E as giving a variant of theorems that functions spaces are isomorphic to the tensor product of invariants of a group $H$ and the set of all $H$-harmonic functions.

Alternatively,
 the above irreducible decomposition may be obtained
 by using the see-saw dual pair
 in the case $\frakg={\mathfrak{su}}(k,k)$
and ${\mathfrak{so}}^{\ast}(4k)$.  
Our proof for Theorem \ref{thm:e}
 is uniformly for all Hermitian Lie algebras
 of tube type including the dual pair
 $(\sl(2,\RR),\frakf_4)$ 
 in the exceptional Lie algebra
$\frakg=\frake_{7(-25)}$.  
In this case $\calH^m(\SS)$
 is the space of spherical harmonics defined on the octonionic projective space $\PP^2(\OO)\simeq F_4/\Spin(9)$.

\subsection*
{Folding maps
 and relations 
 with the classical Fock space}
 
 The classical Schr\"odinger model for the Weil representation
of the metaplectic group $\Mp(n,\RR)$ is realized in $L^2(\RR^n)$, 
whereas our Schr\"odinger model 
 in the case of the Jordan algebra 
 $V=\Sym(n,\RR)$
 is realized 
 in a somewhat different space, 
 that is, 
 the Hilbert space 
 $L^2(\Xi)$ where $\Xi$ is the set of all positive definite real symmetric matrices of rank one.
Since the folding map 
\[
     p:\RR^n\setminus\{0\}\to\Xi,\,x\mapsto x\,{}^t\!x
\]
is double covering,
 it induces an isomorphism $p^{\ast}:
     L^2(\Xi,\td\mu)
     \overset \sim \to L^2_{\textup{even}}(\RR^n)
$, 
the even part of 
$
      L^2(\RR^n), 
$
see Subsection \ref{subsec:foldingmaps}.

Further we discuss in Subsection \ref{subsec:relclassic}
 that there is a close relation
 of our Fock space $\calF(\XX)$
 and the Segal--Bargmann transform $\BB_\Xi$
 with the classical objects
 in the case 
 where $V= \operatorname{Herm}(n,\FF)$, 
 $\FF = \RR$, $\CC$ or $\HH$,
 via the complexified folding map, 
 and we give in Theorem \ref{thm:foldB}
 an alternative proof 
 to a part of Theorems \ref{thm:a} and \ref{thm:c}
 in these special cases.  

However, 
 our main approach to Theorems \ref{thm:a} -- \ref{thm:e}
 is uniform
 for all Hermitian Lie algebras $\frakg$
 of tube type,
 and 
especially for the case
 of the indefinite orthogonal group
 $\frakg=\so(2,n)$
 we obtain a new Fock model
 for the \minholrep\ as well as a Segal--Bargmann transform between the Schr\"odinger model (constructed by T. Kobayashi and B. {\O}rsted in \cite{KO03c}) and the new Fock model.
We examine this case 
 in more detail in Section \ref{sec:exso2n}.  
\\

\subsection*{One-dimensional Jordan algebra
 $V=\RR$}
For $V=\RR$ we have $\frakg=\sl(2,\RR)$ for which the Wallach set is given by $\calW=\{0\}\cup(0,\infty)$. 
In this case,
 our results still hold,
 with continuous parameter $\lambda\in(0,\infty)$ (see Theorems \ref{thm:Frepro}, \ref{thm:UnitaryRepOnFock}\,(2), \ref{thm:NaturalFockSpace}, \ref{thm:SBunitary} and \ref{thm:UnitaryInversionKernel}).
 For the Schr{\"o}dinger model,
 we may use 
 an $L^2$-model of the lowest weight representations
 $\pi_\lambda$ of $\widetilde{\SL(2,\RR)}$
 on $L^2(\RR_+,x^{\lambda-1}\td x)$
 with lowest weight $\lambda$
 studied by B.~Kostant
 \cite{Kos00} 
 and Ranga Rao \cite{RR77}
 (see Subsection \ref{sec:L2branching}). 
For these representations we also obtain a new Fock space realization on $\XX=\CC\setminus\{0\}$ and a Segal--Bargmann transform intertwining Schr\"odinger and Fock model as in the cases explained above. The Fock space $\calF(\XX)$ consists of holomorphic functions on $\CC\setminus\{0\}$ which have finite norm coming from the inner product
\begin{align*}
 \langle F,G\rangle &= \int_\CC{F(z)\overline{G(z)}\widetilde{K}_{\lambda-1}(|z|)|z|^{2(\lambda-1)}\td z}.  
\end{align*}
Theorem \ref{thm:NaturalFockSpace}
 applied to this special case
 means 
 that the origin $0$ is a removable singularity
 of any $F(z) \in  
 \calF(\XX) = \calO(\CC \setminus \{0\}) \cap L^2(\CC \setminus \{0\}, 
 |z|^{2 (\lambda-1)}\td z)$.   
Further the reproducing kernel of the Fock space $\calF(\XX)$ is given by
\begin{align*}
 \KK(z,w) &= \Gamma(\lambda)\widetilde{I}_{\lambda-1}(\sqrt{z\overline{w}}), & z,w\in\CC\setminus\{0\}, 
\end{align*}
and the Segal--Bargmann transform takes the form
\begin{align*}
 \BB_\Xi f(z) &= \Gamma(\lambda)e^{-\frac{1}{2}z}\int_0^\infty{\widetilde{I}_{\lambda-1}(2\sqrt{xz})e^{-x}f(x)x^{\lambda-1}\td x}, & z\in\CC.
\end{align*}
Theorem D in the case where $\frakg=\sl(2,\RR)$
 shows that the unitary inversion operator $\calF_\Xi$ is simply a Hankel transform (see also \cite[Theorem 5.8]{Kos00})
\begin{align*}
 \calF_\Xi u(x) &= 2^{-\lambda}\Gamma(\lambda)\int_0^\infty{\widetilde{J}_{\lambda-1}(2\sqrt{xy})u(y)y^{\lambda-1}\td y}, & x\in\RR_+.  
\end{align*}

{\bf{Notation:}}
\enspace
 $\NN=\{0,1,2,\ldots\}$, $\mathbb{R}_+=\{x\in\mathbb{R}:x>0\}$.
\section{The Schr\"odinger model for minimal representations}
\label{sec:schmodel}

In this section we set up the notation
 and recall briefly the construction of $L^2$-models
 for minimal representations
 of covering groups $G^{\vee}$
 of conformal groups associated to simple real Jordan algebras
{}from \cite{HKM11}. 
Using elliptic,
 self-adjoint differential operators
 $\calBe$ on $(L^2(\Xi), \td \mu)$
 (Lemma \ref{lem:BesselProdRule}), 
 we develop a theory of spherical harmonics
 in this context, 
 and thus find the branching law
 for all \minholreps\ when restricting
 to a distinguished subalgebra $\fraks \simeq \sl(2,\RR)$.
For the basic notion of the Jordan algebra,
 we refer the reader to an excellent book 
 of Faraut--Kor{\'a}nyi
 \cite{FK94}.  

\subsection{Minimal representations}
\label{subsec:minrep}
For a complex simple Lie algebra $\frakg_{\CC}$
 not of type $A_n$,
 A.~Joseph
 \cite{J76}
 introduced a unique completely prime indeal 
 ${\mathcal{J}}$ in $U(\frakg_{\CC})$
 such that the associated variety ${\mathcal{V}}({\mathcal{J}})$
 is equal to the closure of the (complex) minimal 
 nilpotent coadjoint orbit 
 ${\mathbb{O}}_{\operatorname{min}}^{G_{\CC}}$
 in $\frakg_{\CC}^{\ast}$.  
This ideal is primitive 
and called a Joseph ideal.  
\begin{definition}
\label{def:minrep}
For an irreducible unitary representation
 $\pi$ of a real simple Lie group $G$,
 we say $\pi$ is a 
 {\it{minimal representation}}
 if the annihilator
 of the differential representation
 $d \pi$ 
 is the Joseph ideal.  
\end{definition}
In this paper 
 we apply this terminology
 only to Hermitian Lie groups.  
Here,
 a simple non-compact Lie group $G$
 or its Lie algebra $\frakg$
 is said to be {\it{Hermitian}}
 if the associated Riemannian symmetric space $G/K$
 is a Hermitian symmetric space,
 or equivalently,
 if the center ${\frak{z}}(\frakk)$
 of the Lie algebra $\frakk$ is one-dimensional.  
We say a Hermitian Lie group is of tube type
 if $G/K$ has a realization as a tube domain.  

We write $\frakg_{\CC} =\frakk_{\CC}+\frakp_{\CC}$
 for the complexified Cartan decomposition.  
If $G$ is a Hermitian Lie group,
 then the $K$-module $\frakp_{\CC}$ decomposes
 into two irreducible $K$-modules
 ${\mathfrak{p}}_{\CC} = {\mathfrak {p}}_+ + {\mathfrak {p}}_-$, 
and ${\mathfrak {p}}_-$ can be identified 
 with the holomorphic tangent space
 at the base point in $G/K$.  
An irreducible $(\frakg,\frakk)$-module $X$
 is said to be a lowest weight $(\frakg,\frakk)$-module
 if 
\[
     X^{\mathfrak {p}_+}:=\{v \in X: Yv=0\quad
\text { for any }Y \in {\mathfrak {p}}_+\}
\]
 is non-zero.  
We say $X$ is of {\it{scalar type}}
 if $\dim X^{\mathfrak{p}_+}=1$. 
Likewise,
 it is a highest weight $(\frakg,\frakk)$-module
 if $X^{\mathfrak{p}_-}$ is non-zero.  
For a Hermitian group $G$,  
 it is known 
 that minimal representations
 are either highest weight modules
 or lowest weight modules.  

\subsection{The Schr\"odinger model
 for \minholreps}
\label{sec:SchrödingerForMinimalHolomorphic}

First, we recall the construction of the Schr\"odinger model
 ($L^2$-model) for \minholreps.

\subsubsection*{Jordan algebras and related groups}

Let $V$ be a simple Euclidean Jordan algebra, 
 and $\Co(V)_0$ the conformal group of $V$. 
We set $G:=\Co(V)_0$, 
the identity component group.  
We denote by $\vartheta$ the Cartan involution of $G$
 given by $\vartheta (g)=j \circ g \circ j$,
 where $j$ is the conformal inversion $j(x)=-x^{-1}$.  
Then $K:=G^\vartheta$ is a maximal compact subgroup of $G$, 
 and $G$ is the group of biholomorphic transformations
 on an irreducible Hermitian symmetric space $G/K$
 of tube type.  
In particular,
 $G$ is the adjoint group.  
Conversely,
 any simple Hermitian Lie group
 of tube type
 with trivial center arises
 in this fashion.

The Lie algebra $\frakg$ of $G$
 has a Gelfand--Naimark decomposition $\frakg=\frakn+\frakl+\overline{\frakn}$, 
where $\frakn\simeq V$ (abelian Lie algebra), 
$\frakl=\str(V)\subseteq\gl(V)$ is the structure algebra
 and $\overline{\frakn}=\vartheta\frakn$.  
As the differential action
 of the conformal transformations
 of $G$ on $V$, 
the Lie algebra $\frakg$ acts on $V$ as follows:
for $X=(u,T,v) \in V \oplus \frakl \oplus V \simeq 
\frakn \oplus \frakl \oplus \overline{\frakn}=\frak g$,
 the vector field on $V$ is given by 
\begin{align*}
 X(z) = u+Tz-P(z)v,\qquad & z\in V,
\intertext{where}
      P(x)=2L(x)^2-L(x^2), 
      \qquad
      &x\in V,
\end{align*}
 denotes the quadratic representation of $V$ and $L(x)$ the left multiplication by $x$. Thus, $\frakn=\{(u,0,0):u\in V\}$ acts via constant vector fields, $\frakl$ via linear vector fields and $\overline{\frakn}=\{(0,0,v):v\in V\}$ by quadratic vector fields. 
The Cartan involution $\vartheta$ leaves
 $L$ invariant,
 and $K^L:= L \cap K \equiv L^{\vartheta}$
 is a maximal compact subgroup of $L$.  
Correspondingly,
 we have a Cartan decomposition 
 $\frakl = \frakk^{\frakl} + \frakp^{\frakl}$
 of the structure Lie algebra $\frakl=\frak{str}(V)$, 
 where 
\begin{alignat*}{2}
 \frakk^\frakl &= \aut(V) &&= \{D\in\gl(V):D(x\cdot y)=Dx\cdot y+x\cdot Dy\,\forall\,x,y\in V\},\\
 \frakp^\frakl &= L(V) &&= \{L(x):x\in V\}.
\end{alignat*}
The Cartan involution $\vartheta$ acts 
 on $\frakg={\frakn}+{\mathfrak {l}}+\overline{\frakn}$
 by the following formula: 
\begin{align*}
 \vartheta(u,D+L(a),v) &= (-v,D-L(a),-u), 
\end{align*}
and hence
 $\frakk=\mathfrak {g}^{\vartheta}$ is given by 
\begin{align*}
 \frakk &= \{(u,D,-u):u\in V,D\in\frakk^\frakl\}.
\end{align*}
We set 
\begin{align}
 E &:= ({\bf e},0,0), & H &:= (0,2\,\id,0), & F &:= (0,0,{\bf e}).\label{eq:sl2triple}
\end{align}
Then $\{H,E,F\}$ forms
 an $\sl_2$-triple in $\frakg$.  
We define a subalgebra $\fraks$ of $\frakg$ 
 by 
\begin{equation}
\label{eqn:sl2}
  \fraks:= \RR H+ \RR E + \RR F,
\end{equation}
denote by $S\subseteq G$ the corresponding subgroup of $G$, $S\simeq\PP\SL(2,\RR)$. 
Let $\widetilde G$
 be the universal covering group 
 of $G$,
 and denote by $\widetilde{S}\simeq\widetilde{\SL(2,\RR)}$
 the corresponding subgroup in $\widetilde{G}$.  
We define 
 $c \in \operatorname{Int}(\frakg_\CC)$
 by the formula
\begin{align}
 c &:= \exp(-\tfrac{1}{2}\sqrt{-1}\ad(E))\exp(-\sqrt{-1}\ad(F)).\label{eq:DefCayleyTransform}
\end{align}
It is then routine to check the following formulas:
\begin{align}
 c(a,0,0) &= \textstyle(\frac{a}{4},\sqrt{-1}L(a),a),\label{eq:CayleyTransform1}\\
 c(0,L(a)+D,0) &= \textstyle(\sqrt{-1}\frac{a}{4},D,-\sqrt{-1}a),\label{eq:CayleyTransform2}\\
 c(0,0,a) &= \textstyle(\frac{a}{4},-\sqrt{-1}L(a),a).\label{eq:CayleyTransform3}
\end{align}
Therefore,
 the complexifications $\frakk_\CC$ and $\frakl_\CC$ are isomorphic
 to each other by the transform 
$c$:
\begin{lemma}
[Cayley transform]
\label{lem:Cayley}
The transform $c$ induces an isomorphism
 of complex Lie algebras:
\begin{align}
 c: \frakk_\CC \overset \sim \to \frakl_\CC,\,
 (u,D,-u)\mapsto D+2\sqrt{-1}L(u).\label{eq:IsoKCLC}
\end{align}
\end{lemma}
It is convenient to give
 the corresponding matrix computation 
 in $SL(2,\CC)$.  
For this,
 we take the standard $\sl_2$-triple $\{h,e,f\}$
 in $\sl(2,\RR)$ as follows: 
\begin{align*}
 e &:= \left(\begin{array}{cc}0&1\\0&0\end{array}\right), & h &:= \left(\begin{array}{cc}1&0\\0&-1\end{array}\right), & f &:= \left(\begin{array}{cc}0&0\\1&0\end{array}\right).
\end{align*}
Then the element
\begin{align*}
 C &= \exp(-\tfrac{1}{2}\sqrt{-1}e)\exp(-\sqrt{-1}f) = \left(\begin{array}{cc}\tfrac{1}{2} & -\tfrac{1}{2}\sqrt{-1}\\-\sqrt{-1} & 1\end{array}\right)\in\SL(2,\CC)
\end{align*}
defines a Cayley transform
 on $\sl(2,\CC)$ by $\Ad(C)$,
 that is,
\[
\operatorname{Ad}(C)
(\sqrt{-1} \begin{pmatrix} 0 & -1 \\ 1 & 0\end{pmatrix})=
\begin{pmatrix} 1 & 0 \\ 0 & -1\end{pmatrix}
(=h).  
\] 
By the inverse Cayley transform,
 we get another $\sl_2$-triple 
 $\{\widetilde {h}, \widetilde {e}, \widetilde {f}\}$ by
\begin{align*}
 \widetilde{e} &:= \Ad(C)^{-1}e = (e+f-\sqrt{-1}h),\\
 \widetilde{h} &:= \Ad(C)^{-1}h = -\sqrt{-1}(e-f),\\
 \widetilde{f} &:= \Ad(C)^{-1}f = \frac{1}{4}(e+f+\sqrt{-1}h).
\end{align*}
Likewise,
 by the inverse Cayley transform,
 we get another $\sl_2$-triple 
 $\{\widetilde {H}, \widetilde {E}, \widetilde {F}\}$ 
 in $\frakg$ by:
\begin{alignat*}{2}
 \widetilde{E} &:= c^{-1}(E) = (E+F-\sqrt{-1}H) &&= ({\bf e},-2\sqrt{-1}\id,{\bf e}),\\
 \widetilde{H} &:= c^{-1}(H) = -\sqrt{-1}(E-F) &&= -\sqrt{-1}({\bf e},0,-{\bf e}),\\
 \widetilde{F} &:= c^{-1}(F) = \frac{1}{4}(E+F+\sqrt{-1}H) &&= \frac{1}{4}({\bf e},+2\sqrt{-1}\id,{\bf e}).
\end{alignat*}

\subsection*{The Lagrangian submanifold ${\Xi}$}

The structure group $L=\Str(V)_0\subseteq G$
 acts on the Jordan algebra $V$ by linear transformations,
 and $V$ decomposes into finitely many $L$-orbits. The open orbit $\Omega=L\cdot{\bf e}$ is a symmetric cone. 
There is a unique minimal non-zero $L$-orbit in the closure of $\Omega$, 
 which we denote by $\Xi$. 
Note that for $V=\RR$, i.e. $\frakg=\sl(2,\RR)$, we have $\Omega=\RR_+$ and hence also $\Xi=\RR_+$. The orbit $\Xi$ is the orbit of any primitive idempotent $c_1\in V$. From now on we fix such an idempotent $c_1$.

For $x\in V$,
 we denote by $V(x,\nu)$ the eigenspace of $L(x)$
 with eigenvalue $\nu$.  
Then the tangent space of $\Xi$
 is given as follows:
\begin{lemma}
\label{lem:TxO}
For any $x \in \Xi$,
\[
  T_x \Xi = V(x,|x|) \oplus V(x,\frac 1 2 |x|).  
\]
\end{lemma}
\begin{proof}
We begin with the case $|x|=1$.  
Then $x$ is a primitive idempotent,
 and therefore
\[
  V=V(x,0) \oplus V(x,\frac 1 2) \oplus V(x,1).  
\]
Since $L(V)\subseteq\frakl$, 
the tangent space $T_x\Xi=\frakl\cdot x$ contains
 at least $V(x,1)$ and $V(x,\frac{1}{2})$. 
But $\dim\,V(x,1)=1$, $\dim\,V(x,\frac{1}{2})=(r-1)d$, 
and by \cite[Lemma 1.6]{HKM11} we also have $\dim\Xi=1+(r-1)d$. 
Hence $T_x\Xi=V(x,1) \oplus V(x,\frac 1 2)$.  
For general $x \in \Xi$,
 we note that $\frac x{|x|}$ is a primitive idempotent.  
Further,
 since $\Xi$ is a cone,
 the tangent space $T_x\Xi$ is identified with 
 $T_{\frac {x}{|x|}}\Xi$,
 which equals 
$
     V(\frac{x}{|x|},1) \oplus V(\frac{x}{|x|},\frac 1 2)
  =
     V(x,|x|) \oplus V(x,\frac 1 2|x|).  
$
Thus the Lemma is proved.  
\end{proof}

Let 
$\calP(\Xi)$ be the space of restrictions of polynomials on $V$ to the orbit $\Xi$. The space $\calP(\Xi)$
 has a natural grading 
\begin{equation}
\label{eqn:Pm}
     \calP(\Xi)=\bigoplus_{m=0}^\infty{\calP^m(\Xi)}
\end{equation}
 where $\calP^m(\Xi)$ is the space
 of restrictions of homogeneous polynomials
 of degree $m$. 

The orbit $\Xi$ is conical,
 and we have a polar decomposition 
\begin{equation}
\label{eqn:polarO}
 \RR_+ \times \SS \overset \sim \to \Xi,
 (t,x) \mapsto tx,
\end{equation}
where we set 
\begin{equation}
\label{eqn:S}
 \SS := \{ x \in \Xi : |x|=1\}
      = K^L \cdot c_1.  
\end{equation}
Let $\td k$ be the Haar measure on $K^L$
 such that $\int_{K^L} dk =1$. 

We define a Radon measure $\td \mu_{\lambda}$
 on $\Xi$ 
 by using the polar coordinates
 \eqref{eqn:polarO}:
\begin{align}
  \int_\Xi{f(x)\td\mu_{\lambda}(x)} 
 &= \frac{2^{r\lambda}}{\Gamma(r\lambda)}\int_{K^L}{\int_0^\infty{f(ktc_1)t^{r\lambda-1}\td t}\td k}
 \quad
 \text{for }
 f \in C_c(\Xi).
\label{eq:IntFormulaO}
\end{align}
Let $d$ be the multiplicity of the short roots for $G/K$ and put
\begin{alignat}{2}
&\lambda=\frac{d}{2}
&&\text{for }r>1,
\label{eqn:Wal1}
\\
&\lambda\in(0,\infty)
\qquad
&&\text{for }r=1.
\label{eqn:Wal2}
\end{alignat}
See Table \ref{tb:Groups} in the Appendix for the explicit value of $\lambda$.  
Then $\td\mu_\lambda$ is (up to scalar multiples)
 the unique Radon measure on $\Xi$ transforming
 by 
\[
     \td\mu_{\lambda}(gx)
    =\det_V(g)^{\frac{r\lambda}{n}}\td\mu_{\lambda}(x)
\]
 for $g \in L$.  
We have normalized the measure $\td\mu_\lambda$ such that
 the $L^2$-norm of the function $\psi_0(x)=e^{-\tr(x)}$
 (see \eqref{eqn:psi0} below)
 is equal to $1$. We will often write simply $\td \mu$
 when $\dim V >1$, as in this case $\lambda$ is determined uniquely by $V$.

\subsubsection*{Construction of the representation}

The Bessel operator $\calB \equiv \calB_\lambda$ is a vector-valued
 second order differential operator on $V$
 defined by
\begin{align}
\label{eqn:Bessel}
 \calB_\lambda &= P\left(\frac{\partial}{\partial x}\right)x+\lambda\frac{\partial}{\partial x},
\end{align}
where $\frac{\partial}{\partial x}$ denotes the gradient with respect to the trace form $(-|-)$. In an orthonormal basis $(e_\alpha)_\alpha$ of $V$ with respect to $(-|-)$ and coordinates $x=\sum_\alpha{x_\alpha e_\alpha}$ this means
\begin{align*}
 \calB_\lambda u &= \sum_{\alpha,\beta}{\frac{\partial^2 u}{\partial x_\alpha\partial x_\beta}P(e_\alpha,e_\beta)x}+\sum_\alpha{\frac{\partial u}{\partial x_\alpha}e_\alpha}.
\end{align*}
The Bessel operator $\calB_\lambda$ is tangential to the orbit $\Xi$ and defines a differential operator acting on $C^\infty(\Xi)$.

Applying the construction in \cite{HKM11}
 to our setting 
 (that is, $V$ is Euclidean), 
 we obtain an irreducible unitary representation $\pi$
 of the universal covering
 $\widetilde G$
 of $G$
 on the Hilbert space 
 $L^2(\Xi,\td\mu)$, 
where $\td\mu:=\td \mu_{\lambda}$ is defined
 by \eqref{eq:IntFormulaO}.  
It is a minimal representation 
of $\widetilde G$ 
 if $\frakg_{\CC}$ is not of type $A$
(see \cite[Theorem B]{HKM11}).  

Note that for $r>1$ the representation $\pi$ descends to a finite cover $G^\vee$ of $G$,
 whereas for $r=1$ it descends to a finite cover
 of $G=\PP\SL(2,\RR)$ 
 if and only if $\lambda\in\QQ$.

The corresponding Lie algebra action $\td\pi$ is given by
\begin{align}
 \td\pi(u,0,0) &= \sqrt{-1}(u|x),
\notag
\\
 \td\pi(0,T,0) &= D_{T^*x}+\frac{r\lambda}{2n}\Tr(T^*),
\label{eqn:dpiL}
\\
 \td\pi(0,0,v) &= \sqrt{-1}(v|\calB).
\notag
\end{align}
Here $n=\dim\,V$, 
\[
      (x|y)=\frac{r}{n}\Tr(L(xy))
\]
 denotes
 the trace form on the Jordan algebra $V$, 
 and $T^*$ is the adjoint of $T$
 with respect to the trace form $(-|-)$. Further, the notation $D_u$, $u\in V$, is used for the directional derivative:
\begin{align}
 D_uf(x) &= \left.\frac{\td}{\td t}\right|_{t=0}f(x+tu).  
\label{eq:DefDirectionalDerivative}
\end{align}

We define a function $\psi_0$
 on $\Xi$
 by 
\begin{equation}
\label{eqn:psi0}
\psi_0(x):=e^{-\tr(x)}=e^{-|x|}.  
\end{equation}
Then the function $\psi_0$ 
 transforms
 by a character under the action of $K$ and constitutes the minimal $\frakk$-type $W_0=\CC\psi_0$. 
Our normalization of the measure \eqref{eq:IntFormulaO}
 shows that 
\begin{equation}
\label{eqn:normpsi}
\int_{\Xi}|\psi_0(x)|^2d\mu_{\lambda}(x)
=
\frac{2^{r \lambda}}{\Gamma(r\lambda)}
\int_{0}^{\infty} e^{-2r}t^{r\lambda-1} dt 
=1.  
\end{equation}

Let $L^2(\Xi, \td \mu)_{\frakk}$
 denote
 by the space of $\frakk$-finite vectors
 of the representation $(\pi,L^2(\Xi,\td\mu))$, 
 and $\calP(\Xi)$ denote
 by the space of restrictions of polynomials on $V$ to $\Xi$. 
Then we have the following description 
 (see e.g. \cite[Proposition XIII.3.2]{FK94}):
\begin{equation}
\label{eqn:L2K}
  L^2(\Xi, \td \mu)_{\frakk}
 =
 \{f(x) \psi_0(x)
  :
   f \in \calP(\Xi)\}.   
\end{equation}

Let $\alpha_0$ denote the highest weight
 of the one-dimensional representation $W_0$. 
Then the complete decomposition of $\pi$ into $\frakk$-types is given by
\begin{align*}
 L^2(\Xi,\td\mu) &\simeq {\sum_{m=0}^\infty}\raisebox{0.15cm}{\!$^\oplus$}{W_m},
\end{align*}
where $W_m$ is of highest weight $\alpha_0+m\gamma$
 for a certain root $\gamma$.

In particular,
 the underlying $(\frakg,\frakk)$-module
 on $L^2(\Xi, \td \mu)_{\frakk}$
 is an irreducible,
 unitarizable, 
lowest weight module
 of scalar type.  
We recall that irreducible unitary lowest weight representations
 of simple Hermitian groups
 were classified
 by Enright, Howe, and Wallach
 and also by Jakobsen,
 independently.  
Among them,
 those with scalar minimal $\frakk$-type
 ({\it{scalar type}}) of the universal covering group $\widetilde{G}$
 are parameterized by the so-called Berezin--Wallach set $\calW$. 
The parameter  of our representation 
amounts to $\lambda$ defined in 
 \eqref{eqn:Wal1} or \eqref{eqn:Wal2}
 in the normalization 
 of \cite{FK94}, 
 where $\calW$ is given by:
\begin{align*}
 \calW = \{0,\textstyle\frac{d}{2},\ldots,(r-1)\textstyle\frac{d}{2}\}\cup((r-1)\textstyle\frac{d}{2},\infty).
\end{align*}
Here $r$ is the rank 
 of the Hermitian symmetric space $G/K$.  
Thus,
 our representation corresponds
 to the smallest non-zero element of $\calW$ if $r >1$.  

For $\frakg_{\CC}$ is of type $A_n$,
 the Joseph ideal is not defined.  
For $\frakg ={\mathfrak{su}}(n,n)$,
 the representation $\pi$ still has the property
 that the associated variety
 ${\mathcal{V}}(\operatorname{Ann}d \pi)$
 of the annihilator 
 $\operatorname{Ann}(\td \pi)$
in $U(\frakg_{\CC})$
 is the closure
 of the complex minimal nilpotent coadjoint orbit
 $\overline{{\mathbb{O}}_{\operatorname{min}}^{G_{\CC}}}$.  
By an abuse of notation,
 we shall say $\pi$ is a minimal representation
 also for $\frakg={\mathfrak{su}}(n,n)$.  

\subsection{The Schr\"odinger model for complex groups}
\label{sec:MinRepCplx}

Let $V_\CC$ be the complexified Jordan algebra of $V$,
 and $\Co(V_\CC)$ the structure group of $V_\CC$.  
Then $G_\CC:=\Co(V_\CC)_0$ is a natural complexification of $G=\Co(V)_0$.  
In \cite{HKM11} 
 the authors also construct
 an $L^2$-model of an irreducible unitary representation
 of the complex group $G_\CC=\Co(V_\CC)_0$
 that attains the minimum
 of the Gelfand--Kirillov dimensions
 among all infinite-dimensional irreducible unitary representations
 of $G_{\CC}$. 
It should be noted
 that this representation is not a lowest weight representation
 in contrast to the \minholrep\ 
 of $G=\Co(V)_0$
 in Subsection \ref{subsec:minrep}.  
We review the construction briefly.  

As in the Euclidean case, 
 the space $V_\CC$ decomposes into finitely many $L_\CC$-orbits, 
 where the corresponding structure group $L_\CC=\Str(V_\CC)_0$ is a natural complexification of $L$. 
The orbit $L_\CC\cdot{\bf e}$ of the identity element
 is an open cone
 and we denote by $\XX$ the unique minimal non-zero $L_\CC$-orbit
 in its boundary. 

Comparing this complex setting
 with the setting of Subsection \ref{sec:SchrödingerForMinimalHolomorphic}, 
 we see that 
\begin{itemize}
\item[$\bullet$]
     $\XX=L_\CC\cdot c_1$, 
\item[$\bullet$]
$\Xi = L \cdot c$ is a totally real submanifold in $\XX$, 
\item[$\bullet$]
$\XX$ is closed under the complex conjugation
 with respect to $V \subset V_\CC$. 
\end{itemize}
Let $\det_W(g)$ be the real determinant of the $\RR$-linear action of $g\in L_\CC$ on $W=V_\CC$, viewed as a real vector space. 
Again,
 there is a unique $L_\CC$-equivariant measure $\td\nu$
 on $\XX$ 
 up to scaling,
 subject to the equivariant condition: 
 $\td\nu(gz)=\det_W(g)^{\frac{r\lambda}{n}}\td\nu(z)$.  
In terms of the polar decomposition $\XX=K^{L_\CC}\RR_+c$, 
 where $K^{L_\CC}$ is a maximal compact subgroup 
 in $L_\CC$, we normalize $\td\nu$ by
\begin{align}
 \int_\XX{F(z)\td\nu(z)} &= \frac{1}{c_{r,\lambda}}\int_{K^{L_\CC}}{\int_0^\infty{F(utc_1)t^{2r\lambda-1}\td t}\td u},\label{eq:IntFormulaX}
\end{align}
where $\td u$ denotes the normalized Haar measure on $K^{L_\CC}$ and
\begin{align*}
 c_{r,\lambda} &= 2^{2r\lambda-2}\Gamma\left(r\lambda\right)\Gamma\left((r-1)\lambda+1\right).
\end{align*}
(This normalization is chosen such that in the Fock model the constant polynomial $\1$ is a unit vector.) 
Then we can construct
 an irreducible unitary representation, 
 to be denoted by $\tau$, 
 of $G_\CC$ on the Hilbert space 
 $L^2(\XX,\td\nu)$ (\cite{HKM11}). 
The Lie algebra action $\td\tau$ is given by
\begin{align*}
 \td\tau(u,0,0) &= \sqrt{-1}(u|z)_W,\\
 \td\tau(0,T,0) &= D_{T^*z}+\frac{r\lambda}{2n}\Tr_W(T^*),\\
 \td\tau(0,0,v) &= \sqrt{-1}(v|\calB^W)_W.
\end{align*}
Here $(z|w)_W=2\left((\Re(z)|\Re(w))+(\Im(z)|\Im(w))\right)$ defines an inner product on the real Jordan algebra $W=V_\CC$, by $\Tr_W$ we mean the real trace of an operator on the real vector space $W$, and $\calB^W=\calB^W_\lambda$ denotes the real Bessel operator of $W$. The Bessel operator can be defined using a real basis $(e_\alpha)_\alpha$ of $W$ and its dual basis $(\overline{e}_\alpha)_\alpha$ with respect to the non-degenerate bilinear form 
 ({\it{trace form}})
 $(z,w)\mapsto2\left((\Re(z)|\Re(w))-(\Im(z)|\Im(w))\right)$:
\begin{align*}
 \calB_\lambda^W &= \sum_{\alpha,\beta}{\frac{\partial^2}{\partial x_\alpha\partial x_\beta}P(\overline{e}_\alpha,\overline{e}_\beta)x}+\sum_\alpha{\frac{\partial}{\partial x_\alpha}\overline{e}_\alpha}.
\end{align*}
Note that $\td\tau$ does not act via holomorphic differential operators, but via real differential operators up to second order on $\XX$.

We shall find an action
 on the Fock space
 by holomorphic differential operators later.  
For this, 
 we observe that $\td\pi$ acts on $C^\infty(V)$ by polynomial differential operators. Using the Wirtinger derivative
\begin{align*}
 \frac{\partial}{\partial z} &= \frac{1}{2}\left(\frac{\partial}{\partial x}-\sqrt{-1}\frac{\partial}{\partial y}\right)
\end{align*}
the action $\td\pi$ extends uniquely to a $\CC$-linear action $\td\pi_\CC$ of $\frakg_\CC$ on $C^\infty(V_\CC)$ by holomorphic differential operators. In particular,
\begin{align*}
 \td\pi_\CC(a,0,0) &= \sqrt{-1}(a|z),\\
 \td\pi_\CC(0,0,a) &= \sqrt{-1}(a|\calB),
\end{align*}
where $(-|-)$ denotes the extension of the trace form on $V$ to a $\CC$-bilinear form on $V_\CC$ and the Bessel operator $\calB$ extends
 to a holomorphic differential operator on $V_\CC$ by the formula
\begin{align}
 \calB &= P\left(\frac{\partial}{\partial z}\right)z+\lambda\frac{\partial}{\partial z}.\label{eq:DefHolBessel}
\end{align}
To show that $\td\pi_\CC$ indeed defines an action on $C^\infty(\XX)$, we use the following result which expresses $\td\pi_\CC$ also as a kind of Wirtinger derivative.

\begin{proposition}\label{prop:PiCC}
For $X\in\frakg_\CC$ we have
\begin{align}
 \td\pi_\CC(X) &= \frac{1}{2}\left(\td\tau(X)-\sqrt{-1}\td\tau(\sqrt{-1}X)\right).\label{eq:ComplexifiedRep}
\end{align}
In particular, for every $X=(u,T,v)\in\frakg_\CC$ and all $F,G\in C^\infty(\XX)$ we have
\begin{align*}
 \int_\XX{\td\pi_\CC(u,T,v)F(z)\cdot G(z)\td\nu(z)} &= \int_\XX{F(z)\cdot\td\pi_\CC(u,-T,v)G(z)\td\nu(z)}.
\end{align*}
\end{proposition}

\begin{proof}
Let $(e_j)$ be any orthonormal basis of $V$ with respect to the trace form $(-|-)$. Write $x=\sum_j{x_je_j}$. Then $\frac{\partial}{\partial x}=\sum_j{\frac{\partial}{\partial x_j}e_j}$.\\
We now view $W=V_\CC$ as a real Jordan algebra. Then $f_j:=\frac{1}{\sqrt{2}}e_j$ and $g_j:=\frac{1}{\sqrt{2}}\sqrt{-1}e_j$ constitute an $\RR$-basis of $V_\CC$ with dual basis with respect to the trace form given by $(\overline{f}_j:=f_j)_j\cup(\overline{g}_j:=-g_j)_j$. We write $z=\sum_j{z_je_j}=\sum_j{(a_jf_j+b_jg_j)}$ with $a_j,b_j\in\RR$ and $z_j=\sqrt{2}(a_j+\sqrt{-1}b_j)$. Hence, $\frac{\partial}{\partial a_j}=\frac{1}{\sqrt{2}}\frac{\partial}{\partial x_j}$ and $\frac{\partial}{\partial b_j}=\frac{1}{\sqrt{2}}\frac{\partial}{\partial y_j}$ with $z_j=x_j+\sqrt{-1}y_j$, $x_j,y_j\in\RR$. Then the gradient in $W$, viewed as real Jordan algebra, is given by
\begin{align*}
 \sum_j{\frac{\partial}{\partial a_j}\overline{f}_j}+\sum_j{\frac{\partial}{\partial b_j}\overline{g}_j} &= \frac{1}{2}\sum_j{\left(\frac{\partial}{\partial x_j}-\sqrt{-1}\frac{\partial}{\partial y_j}\right)e_j}
\end{align*}
which is the same as the Wirtinger derivative in the complex Jordan algebra $V_\CC$. Now \eqref{eq:ComplexifiedRep} is easily verified by explicit computations.\\
The integral formula follows from \eqref{eq:ComplexifiedRep} using the fact that $\td\tau(X)$ is given by skew-adjoint real differential operators operators on $L^2(\XX,\td\nu)$ with real coefficients if $X=(0,T,0)\in\frakg$ and purely imaginary coefficients if $X=(u,0,v)\in\frakg$.
\end{proof}

Since $\td\tau$ restricts to an action on $C^\infty(\XX)$ by differential operators, the same is true for $\td\pi_\CC$ by the previous proposition. 
Therefore $\td\pi_\CC$ is a representation of $\frakg_\CC$ on $C^\infty(\XX)$ by differential operators of order at most $2$.

\subsection{The Bessel operator and a related second order ODE}

Let $\calB \equiv \calB_{\lambda}$
 be the Bessel operator
 defined in \eqref{eqn:Bessel}.  
We first recall the product rule of the Bessel operator which will be used repeatedly.

\begin{lemma}[{\cite[Lemma 1.7.1]{Moe10}}]\label{lem:BesselProdRule}
\begin{align*}
 \calB(f(z)g(z)) &= \calB f(z)g(z)+2P\left(\frac{\partial f}{\partial z}(z),\frac{\partial g}{\partial z}(z)\right)z+f(z)\calB g(z).
\end{align*}
\end{lemma}

We further introduce the identity component
 of the Bessel operator by 
\begin{align}
 \calBe &:= -\sqrt{-1} \td \pi(F) = ({\bf e}|\calB),\label{eq:IdentityBessel}
\end{align}
where we recall $F=(0,0,{\bf e})$.  
Here we give some basic properties
 of the operator $\calBe$.  
An analogue of the heat kernel
 corresponding to $\calBe$
 will be discussed in Section \ref{sec:heat}.  

\begin{lemma}
\label{lem:BesseleElliptic}
$\calBe$ is an elliptic differential operator 
 of second order on $\Xi$.
Further,
 it defines a self-adjoint operator
 on $L^2(\Xi, \td \mu)$.  
\end{lemma}

\begin{proof}
We already know 
 that $\calBe$ is a differential operator
 of second order
 along $\Xi$.  
Let us compute the principal symbol.  
For $x\in\Xi$ 
 we identify $T_x^*\Xi$ with $T_x\Xi$
 via the trace form $(-|-)$ on $V$. 
Then the principal symbol of $\calBe$ at $x\in\Xi$ in direction $\xi\in T_x^*\Xi$ is given by
\begin{align*}
 (P(\xi)x|{\bf e}) &= (x|P(\xi){\bf e}) = (x|\xi^2) = (L(x)\xi|\xi).
\end{align*}
By Lemma \ref{lem:TxO},
 $L(x)$ has eigenvalues $|x|$ and $\frac{1}{2}|x|$ on $T_x^*\Xi$.  
Therefore,
\begin{align*}
 |(L(x)\xi|\xi)| \geq \frac{|x|}{2}|\xi|^2
\end{align*}
which implies that $\calBe$ is an elliptic operator.

The last statement follows from the fact that
 $(\pi,L^2(\Xi, \td \mu))$
 is a unitary representation.  
\end{proof}

\begin{remark}
In \cite{Men11} G. Meng studies some physics-related aspects
 of the lowest weight representations
 of Hermitian Lie groups of tube type. 
In particular,
 he constructs a Riemannian metric
 on the cone $\Xi$. 
Its corresponding Laplace operator $\Delta$
 has principal symbol $\frac{1}{\tr(x)}(\xi|L(x)\xi)$
 and hence the principal symbols
 of $\tr(x)\Delta$ and $\calBe$ agree.
\end{remark}

Using the product rule one can calculate the action of $\calBe$ on products of homogeneous polynomials and powers of the trace. For this we let $\calP^m(V)$ be the space of homogeneous polynomials on $V$ of degree $m$.

\begin{lemma}\label{lem:BesselOnProductWithTrace}
For every $p\in\calP^m(V)$ and every $k\in\NN$ we have
\begin{align*}
 \calBe(\tr^k(x)p) &= k(r\lambda+2m+k-1)\tr^{k-1}(x)p+\tr^k(x)\calBe p.
\end{align*}
\end{lemma}

\begin{proof}
First note that the commutator relation $[\calBe,\tr(x)]=2\calE+r\lambda$ holds where $\calE=(x|\frac{\partial}{\partial x})$ is the Euler operator. In fact, using Lemma \ref{lem:BesselProdRule} we obtain
\begin{align*}
 \calBe(\tr(x)f(x)) &= \calBe\tr(x)\cdot f(x)+2\left(\left.P\left(\frac{\partial}{\partial x}\tr(x),\frac{\partial f}{\partial x}(x)\right)x\right|{\bf e}\right)+\tr(x)\cdot\calBe f(x)\\
 &= r\lambda f(x)+2\calE\!f(x)+\tr(x)\cdot\calBe f(x).
\end{align*}
To show the claim we now proceed by induction on $k$. For $k=0$ the claim is trivial. For $k>0$ we find on $\calP^m(V)$:
\begin{align*}
 [\calBe,\tr^k(x)] &= [\calBe,\tr^{k-1}(x)]\tr(x)+\tr^{k-1}(x)[\calBe,\tr(x)]\\
 &= (k-1)(r\lambda+2(m+1)+k-2)\tr^{k-1}(x)+(r\lambda+2m)\tr^{k-1}(x)\\
 &= k(r\lambda+2m+k-1)\tr^{k-1}(x)
\end{align*}
since $\calE$ acts on $\calP^m(V)$ by the scalar $m$.
\end{proof}

Further we need the following simple result on the quadratic representation on the orbit $\XX$:

\begin{lemma}\label{lem:QuadrProj}
For $z\in\XX$ we have
\begin{align*}
 P(z) &= (z|-)z.
\end{align*}
\end{lemma}

\begin{proof}
Write $z=gc_1$ with $g\in L_\CC$. We recall that $P(c_1)$ is the orthogonal projection onto $\CC c_1$ and hence given by $P(c_1)=(c_1|-)c_1$. Thus, we obtain for $w\in V_\CC$:
\begin{align*}
 P(z)w &= P(gc_1)w = gP(c_1)g^*w = (c_1|g^*w)gc_1 = (z|w)z.\qedhere
\end{align*}
\end{proof}

For the statement of the next formula denote by $\calB_z$ the Bessel operator $\calB$ acting in the variable $z\in V_\CC$.

\begin{lemma}
For $z\in V_\CC$ and $w\in\XX$ we have
\begin{align}
 \calB_z(z|w)^k &= k(\lambda+k-1)(z|w)^{k-1}w.\label{eq:BesselOnPowers}
\end{align}
\end{lemma}

\begin{proof}
Since $\frac{\partial}{\partial z}(z|w)^k=k(z|w)^{k-1}w$, we obtain
\begin{align*}
 \calB_z(z|w)^k &= k(k-1)(z|w)^{k-2}P(w)z+k\lambda(z|w)^{k-1}w.
\end{align*}
Now, $w\in\XX$ and hence $P(w)z=(z|w)w$ by Lemma \ref{lem:QuadrProj}. This yields
\begin{align*}
 \calB_z(z|w)^k &= k(\lambda+k-1)(z|w)^{k-1}w.\qedhere
\end{align*}
\end{proof}

Since the map 
\[
      \Xi\times\Xi\to\RR,\,(x,y)\mapsto(x|y), 
\]
has nowhere vanishing derivatives,
  the pullback $u((x|y))\in\calD'(\Xi\times\Xi)$ is well-defined
 for any distribution $u\in\calD'(\RR)$
 on $\RR$. 
We write $u(x|y)$ for short.

\begin{proposition}
\label{prop:BesselHypergeomEq}
Let $u\in\calD'(\RR)$. Then
\begin{align*}
 \calB_xu(x|y) &= yu(x|y)
\end{align*}
if and only if $u$ solves the differential equation
\begin{align}
 tu''+\lambda u'-u &= 0.\label{eq:HyperGeomEq}
\end{align}
Moreover, for $\lambda\notin(-\NN)$
 the renormalized Bessel functions 
(see Appendix \ref{app:BesselFcts})
\begin{align*}
 \widetilde{I}_{\lambda-1}(2\sqrt{t}) 
 &= \frac{1}{\Gamma(\lambda)}{_0F_1}(\lambda;t) 
& \mbox{and}
\\
 \widetilde{K}_{\lambda-1}(2\sqrt{t})
\end{align*}
form a fundamental system of solutions 
to \eqref{eq:HyperGeomEq} on $\RR_+$. In particular
\begin{itemize}
\item $\widetilde{I}_{\lambda-1}(2\sqrt{t})$ is the unique
 (up to scalar multiples) entire solution
 of \eqref{eq:HyperGeomEq} and
\item $\widetilde{K}_{\lambda-1}(2\sqrt{t})$ is the unique
 (up to scalar multiples) solution
 to \eqref{eq:HyperGeomEq}
 which decays exponentially as $t$ tends to $+\infty$.
\end{itemize}
\end{proposition}

\begin{proof}
Since $\calB_x=P(\frac{\partial}{\partial x})x+\lambda\frac{\partial}{\partial x}$ and $\frac{\partial}{\partial x}(x|y)=y$, we obtain
\begin{align*}
 \calB_xu(x|y)-yu(x|y) &= P(y)xu''(x|y)+\lambda yu'(x|y)-yu(x|y).
\end{align*}
Now, for $y\in\Xi$ we have $P(y)x=(x|y)y$ by Lemma \ref{lem:QuadrProj}. 
Therefore we obtain
\begin{align*}
 \calB_xu(x|y)-yu(x|y) &= \left((x|y)u''(x|y)+\lambda u'(x|y)-u(x|y)\right)y
\end{align*}
which gives the claim.
\end{proof}

\subsection{A theory of spherical harmonics}
\label{sec:SphericalHarmonics}

The geometry $\Xi$
 for the Schr{\"o}dinger model
 has the polar coordinates
\[
  \Xi \simeq \RR_+ \times \SS,
\]
where $\SS:=K^L\cdot c_1$
 is a compact manifold.  
As a homogeneous space,
 $\SS \simeq K^L/K^L_{c_1}$
 where $K^L_{c_1}=\Stab_{K^L}(c_1)$. 
Since $K^L/K_{c_1}^L$ is a compact symmetric space,
 the irreducible decomposition
 of $L^2(\SS)$
 is multiplicity-free 
 and we can tell explicitly
 the highest weights
 of these irreducible representations
 of $K^L$
 by the Cartan--Helgason theorem for $K^L$ semisimple and by usual Fourier expansion for $K^L\simeq\SO(2)$ (no other cases occur).
 This subsection takes
 a Jordan theoretic approach 
 to define the space of spherical harmonics
\[
  \calH^m(\SS) \simeq \calH^m(\Xi)
  \qquad
  (m \in \NN)
\]
by using the elliptic differential operator $\calBe$, 
 introduced in \eqref{eq:IdentityBessel}, 
 and give a concrete decomposition
 of the left-regular representation of $K^L$ on $L^2(\SS)$. 

In the case 
$
     V = \operatorname{Herm}(k,\FF),
$
 $\FF = \RR$, $\CC$, or $\HH$,
 we discuss in Subsection \ref{subsec:foldingmaps}
 a relation of our spherical harmonics
 with the classical spherical harmonics
 on the sphere.  

Throughout this subsection
 we assume $r > 1$
 because for $r=1$ the orbit $\SS$ consists of a single point.  

We complete the idempotent $c_1$ to a Jordan frame $c_1,c_2,\ldots,c_r$ and denote the corresponding Peirces spaces by $V_{ij}$, $1\leq i,j\leq r$. Choose any $x_0\in V_{12}$ with $\|x_0\|^2=2$ and put $X_0:=[L(c_1),L(x_0)]$ and $\frakt:=\RR X_0\subseteq\frakk^\frakl$. We further define $\gamma\in\frakt_\CC^*$ by $\gamma(X_0)=\frac{1}{4}\sqrt{-1}$.

\begin{proposition}\label{prop:HighestWeightsSphericalReps}
$\frakt$ is a maximal abelian subspace
 in the orthogonal complement of $\frakk^\frakl_{c_1}$ in $\frakk^\frakl$. 
Moreover, the restricted root system $\Phi(\frakk^\frakl_\CC,\frakt_\CC)$
 is given by
\begin{align*}
 \Phi(\frakk^\frakl_\CC,\frakt_\CC) &= \begin{cases}\emptyset & \mbox{ for $d=1$, $r=2$,}\\\{\pm\gamma\} & \mbox{ for $d=1$, $r>2$,}\\\{\pm2\gamma\} & \mbox{ for $d>1$, $r=2$,}\\\{\pm\gamma,\pm2\gamma\} & \mbox{ for $d>1$.}\end{cases}
\end{align*}
\end{proposition}

\begin{proof}
We note that
\begin{align*}
 \frakk^\frakl &= \frakk^\frakl_0\oplus\bigoplus_{i<j}{\frakk^\frakl_{ij}},\\
\intertext{where}
 (\frakk^\frakl)_0 
 &= \{D\in\frakk^\frakl:Dc_i=0\,\quad\text{for all }\,i=1,\ldots,r\},\\
 \frakk^\frakl_{ij} &= \{[L(c_i),L(x)]=-[L(c_j),L(x)]:x\in V_{ij}\}.
\end{align*}
In this notation
\begin{align*}
 \frakk^\frakl_{c_1} &= \frakk^\frakl_0\oplus\bigoplus_{1<i<j}{\frakk^\frakl_{ij}}, & \frakm &:= (\frakk^\frakl_{c_1})^\perp = \bigoplus_{i=2}^r{\frakk^\frakl_{1i}}.
\end{align*}
\begin{enumerate}
\item[\textup{(1)}] We first show
 that $\frakt$ is a maximal abelian subspace
 in $\frakm$. 
Suppose $y\in\bigoplus_{i=2}^r{V_{1i}}$ with $[X_0,[L(c_1),L(y)]]=0$. Write $y=y'+y''$ with $y'\in V_{12}$ and $y''\in\bigoplus_{i=3}^r{V_{1i}}$. Then it is easy to see that $X_0y'\in V_{11}\oplus V_{22}$ and $X_0y''\in\bigoplus_{i=3}^r{V_{2i}}$. 
Therefore $[L(c_1),L(X_0y)]=0$ and hence
\begin{align*}
 0 &= [X_0,[L(c_1),L(y)]]\\
 &= [L(X_0c_1),L(y)]+[L(c_1),L(X_0y)]\\
 &= -[L(\textstyle\frac{x}{4}),L(y)].
\end{align*}
We obtain $[L(x),L(y)]=0$. Applying this operator to $c_2$ gives
\begin{align*}
 0 &= [L(x),L(y)]c_2 = -\frac{xy''}{2}.
\end{align*}
By \cite[Lemma IV.2.2]{FK94} we obtain $\|xy''\|^2=\frac{1}{8}\|x\|^2\|y''\|^2$ and hence $y''=0$. 
Therefore $y=y'\in V_{12}$. Applying the operator $[L(x),L(y)]=0$ to $x$ this time gives
\begin{align*}
 0 &= [L(x),L(y)]x = \frac{(x|y)}{2}x-y
\end{align*}
and hence $y=\frac{(x|y)}{2}x$ and therefore $[L(c_1),L(y)]\in\frakt$. 
Thus $\frakt$ is maximal in $\frakm$.
\item[\textup{(2)}] We now calculate the root system $\Phi(\frakk^\frakl_\CC,\frakt_\CC)$. For this we observe the following mapping properties of $\ad(X_0)$:
\begin{itemize}
\item $D\in\frakk^\frakl_0$:
\begin{align*}
 \ad(X_0)D &= -[L(c_1),L(Dx_0)]\in\frakk^\frakl_{12}.
\end{align*}
\item $[L(c_1),L(y)]\in\frakk^\frakl_{12}$:
\begin{align*}
 \ad(X_0)[L(c_1),L(y)] &= -\frac{1}{4}[L(x_0),L(y)]\in\frakk^\frakl_0.
\end{align*}
\item $[L(c_1),L(y)]\in\frakk^\frakl_{1i}$, $i\geq3$:
\begin{align*}
 \ad(X_0)[L(c_1),L(y)] &= \frac{1}{2}[L(c_2),L(x_0y)]\in\frakk^\frakl_{2i}.
\end{align*}
\item $[L(c_2),L(y)]\in\frakk^\frakl_{2i}$, $i\geq3$:
\begin{align*}
 \ad(X_0)[L(c_2),L(y)] &= \frac{1}{2}[L(c_1),L(x_0y)]\in\frakk^\frakl_{1i}.
\end{align*}
\item $D\in\frakk^\frakl_{ij}$, $2<i<j$:
\begin{align*}
 \ad(X_0)D &= 0.
\end{align*}
\end{itemize}
From this one easily obtains the following root spaces:
\begin{align*}
 (\frakk^\frakl_\CC)_{\pm\gamma} &= \left\{2[L(c_1),L(x_0y)]\pm\sqrt{-1}[L(c_2),L(y)]:y\in\bigoplus_{i=3}^r{(V_{2i})_\CC}\right\},\\
 (\frakk^\frakl_\CC)_{\pm2\gamma} &= \left\{2[L(c_1),L(y)]\pm\sqrt{-1}[L(x_0),L(y)]:y\in(V_{12})_\CC,y\perp x_0\right\},\\
 (\frakk^\frakl_\CC)_0 &= \left\{D\in(\frakk^\frakl_0)_\CC:Dx_0=0\right\}\oplus\CC X_0\oplus\bigoplus_{2<i<j}^r{(\frakk^\frakl_{ij})_\CC},
\end{align*}
which shows the claim.\qedhere
\end{enumerate}
\end{proof}

\begin{remark}
The case $r=2$ and $d=1$ occurs exactly for the Jordan algebra $V=\Sym(2,\RR)\simeq\RR^{1,2}$ with conformal Lie algebra $\frakg=\sp(2,\RR)\simeq\so(2,3)$. In this case $\frakk^\frakl\simeq\so(2)$. In all other cases $\frakk^\frakl$ is semisimple (see e.g. Table \ref{tb:Groups}).
\end{remark}

Recall the identity component $\calBe=({\bf e}|\calB)$
 of the Bessel operator (see \eqref{eq:IdentityBessel}).  

Suppose $p \in \calP^m(\Xi)$,
 a homogeneous polynomial on $\Xi$
 (see \eqref{eqn:Pm}).  
We say $p$ is a {\it{spherical harmonic}}
 on $\Xi$ of degree $m$
 if $\calBe p=0$.  
We define the space
 of \textit{harmonic polynomials} on $\Xi$ of degree $m$
by 
\begin{align}
 \calH^m(\Xi) &:= \{p\in\calP^m(\Xi):\calBe p=0\}, 
\label{eq:DefHarmonicsOnO}
\end{align}
and set 
 ${\mathcal{H}}^m(\SS)
:=\{p|_{\SS}: p \in {\mathcal{H}}^m({\Xi})\}$. 
Since homogeneous polynomials are already determined by their values on the sphere $\SS=\{x\in\Xi:|x|=1\}$,
we have an obvious isomorphism
${\mathcal{H}}^m(\Xi)
 \overset{\sim}\to {\mathcal{H}}^m(\SS)$
 by restriction.  

Recall the subalgebra $\fraks\simeq\sl(2,\RR)$ of $\frakg$ introduced in \eqref{eqn:sl2}. We will later see (Proposition \ref{prop:DualPair}) that $Z_\frakg(\fraks)=\frakk^\frakl$ and hence $Z_G(\fraks)$ is a possibly disconnected subgroup of $G$ with Lie algebra $\frakk^\frakl$.

\begin{proposition}\label{prop:SphericalHarmonicsReps}
\begin{enumerate}
\item[\textup{(1)}] Let $r>2$ or $d>1$. Then each $\calH^m(\SS)$ is an irreducible $K^L_{c_1}$-spherical representation of $K^L$ of highest weight $2m\gamma$.
\item[\textup{(2)}] Let $r=2$ and $d=1$, i.e. $V\simeq\RR^{1,2}\simeq\Sym(2,\RR)$. Then $\calH^m(\SS)$ is an irreducible representation of $Z_G(\fraks)\simeq O(2)$. It decomposes into two irreducible non-isomorphic representations of $K^L\simeq\SO(2)$ for $m>0$ and is the trivial representation of $K^L$ for $m=0$.
\end{enumerate}
\end{proposition}

\begin{proof}
\begin{enumerate}
\item[\textup{(1)}] Let $r>2$ or $d>1$. Then $\frakk^\frakl$ is semisimple. The Killing form on $\frakk^\frakl$ is
 a scalar multiple of the trace form
\begin{align*}
 \frakk^\frakl\times\frakk^\frakl\to\RR,\,(D,D')\mapsto\Tr_V(DD').
\end{align*}
Examining closely the mapping properties of elements in $\frakm=(\frakk^\frakl_{c_1})^\perp\subseteq\frakk^\frakl$ one can show
 that on $\frakm$ the trace form is again a scalar multiple
 of the $\frakk^\frakl_{c_1}$-invariant inner product
\begin{align*}
 \kappa([L(c_1),L(x)],[L(c_1),L(y)]) &:= (x|y), & \mbox{for $x,y\in V(c,\frac{1}{2})$.}
\end{align*}
We denote by $\kappa$ also its extension
 to the scalar multiple of the Killing form
 on $\frakk^\frakl$. 
Let $\Omega\in\calU(\frakk^{\frakl})$ be the Casimir operator corresponding to this form and denote its action on a function $p\in C^\infty(\SS)$ by $\Omega\cdot p$. We show that $\Omega$ acts on the spherical $\frakk^\frakl$-representation with highest weight $k\gamma$ by the scalar $-\frac{1}{32}k(rd+k-2)$ and that $\Omega$ acts on $\calH^m(\SS)$ by the scalar $-\frac{1}{8}m(\frac{rd}{2}+m-1)$. Note that the scalars $-\frac{1}{32}k(rd+k-2)$, $k\in\NN$, are all distinct. Since $\calH^m(\SS)\subseteq L^2(\SS)$ is a $K^L$-invariant subspace, it then has to be an irreducible $K^L$-representation with highest weight $2m\gamma$.
\begin{enumerate}
\item[\textup{(1)}] The Casimir operator $\Omega$ acts
 on the irreducible $\frakk^\frakl$-representation
 with highest weight $\alpha$
 as a scalar $\widetilde{\kappa}(\alpha,\alpha+2\rho)$,
 where $2\rho$ is the sum of all positive roots. 
Here $\widetilde{\kappa}$ denotes the bilinear form
 on $\frakm_\CC^*$ corresponding to $\kappa$
 under the identification $\frakm^*\simeq\frakm$ via $\kappa$. 
We find that
\begin{align*}
 \rho &= \frac{1}{2}(m_{2\gamma}\cdot2\gamma+m_\gamma\cdot\gamma)\\
 &= \frac{1}{2}((d-1)\cdot2\gamma+(r-2)d\cdot\gamma)\\
 &= \left(\frac{rd}{2}-1\right)\gamma.
\end{align*}
For $\alpha=k\gamma$ we then obtain
\begin{align*}
 \widetilde{\kappa}(\alpha,\alpha+2\rho) &= k(rd+k-2)\widetilde{\kappa}(\gamma,\gamma) = -\frac{1}{64}k(rd+k-2)\kappa(X_0,X_0)\\
 &= -\frac{1}{32}k(rd+k-2).
\end{align*}
\item[\textup{(2)}] Since $\Omega$ is $K^L$-invariant,
 it suffices to show
 that 
\[
     (\Omega\cdot p)(c_1)=-\frac{1}{8}m(\frac{rd}{2}+m-1)p(c_1)
\]
 for every $p\in\calH^m(\SS)$. For each $j=2,\ldots,r$ we choose an orthonormal basis $(e_{jk})_{k=1,\ldots,d}$ of $V_{1j}$. Then the elements $[L(c_1),L(e_{jk})]$, $j=2,\ldots,r$, $k=1,\ldots,d$, form an orthonormal basis of $\frakm$ with respect to the inner product $\kappa$. Since for every element $X\in\frakk^\frakl_{c_1}$ we have $(X\cdot p)(c_1)=\left.\frac{\td}{\td t}\right|_{t=0}p(e^{tX}c_1)=0$, we obtain
\begin{align*}
 (\Omega\cdot p)(c_1) &= \sum_{j=2}^r{\sum_{k=1}^d{\left.D_{[L(c_1),L(e_{jk})]x}^2p(x)\right|_{x=c_1}}}
\end{align*}
(see \eqref{eq:DefDirectionalDerivative} for the definition of $D_u$, $u\in V$). Fix $j$ and $k$ and put $A=[L(c_1),L(e_{jk})]$. Then $Ac_1=-\frac{1}{4}e_{jk}$ and $Ae_{jk}=\frac{1}{4}(c_1-c_j)$ and we find
\begin{align*}
 D_{Ax}^2p(c_1) &= \left(Ac_1\left|\frac{\partial}{\partial x}\left(Ax\left|\frac{\partial p}{\partial x}\right.\right)_{x=c_1}\right.\right)\\
 &= -\frac{1}{4}\left(e_{jk}\left|\frac{\partial}{\partial x}\left(Ax\left|\frac{\partial p}{\partial x}\right.\right)_{x=c_1}\right.\right)\\
 &= -\frac{1}{4}\frac{\partial}{\partial x_{jk}}\left(Ax\left|\frac{\partial p}{\partial x}\right.\right)_{x=c_1}\\
 &= -\frac{1}{4}\left(Ae_{jk}\left|\frac{\partial p}{\partial x}(c_1)\right.\right)-\frac{1}{4}\left(Ac_1\left|\frac{\partial}{\partial x_{jk}}\frac{\partial p}{\partial x}_{x=c_1}\right.\right)\\
 &= -\frac{1}{16}\left(c_1-c_j\left|\frac{\partial p}{\partial x}(c_1)\right.\right)+\frac{1}{16}\left(e_{jk}\left|\frac{\partial}{\partial x_{jk}}\frac{\partial p}{\partial x}_{x=c_1}\right.\right)\\
 &= -\frac{1}{16}\left(\frac{\partial p}{\partial x_1}(c_1)-\frac{\partial f}{\partial x_j}(c_1)\right)+\frac{1}{16}\frac{\partial^2p}{\partial x_{jk}^2}(c_1),
\end{align*}
where $x_1$ and $x_j$ denote the coordinates of $c_1$ and $c_j$ and $x_{jk}$ the coordinates of $e_{jk}$ in an arbitrary extension to an orthonormal basis. Hence,
\begin{align*}
 (\Omega\cdot p)(c_1) &= \frac{1}{16}\sum_{j=2}^r{\sum_{k=1}^d{\frac{\partial^2p}{\partial x_{jk}^2}}}+\frac{d}{16}\sum_{j=1}^r{\frac{\partial p}{\partial x_j}(c_1)}-\frac{rd}{16}\frac{\partial p}{\partial x_1}(c_1).
\end{align*}
We now use that $p$ is harmonic and calculate the action of the Bessel operator $\calBe$ at the point $x=c_1$. Let $(e_\alpha)_\alpha\subseteq V$ be an orthonormal basis of $V$ extending $(x_{jk})_{j,k}\cup(c_j)_j$ such that every element is contained in some $V_{jk}$. For this note that
\begin{align*}
 (P(e_\alpha,e_\beta)c_1|{\bf e}) &= \begin{cases}1 & \mbox{ if $e_\alpha=e_\beta\in V_{11}$,}\\\frac{1}{2} & \mbox{ if $e_\alpha=e_\beta\in V_{1j}$, $j=2,\ldots,r$,}\\0 & \mbox{ else,}\end{cases}\\
\intertext{and}
 (e_\alpha|{\bf e}) &= \begin{cases}1 & \mbox{ if $e_\alpha\in V_{jj}$ for $j=1,\ldots,r$,}\\0 & \mbox{ else.}\end{cases}
\end{align*}
We find
\begin{align*}
 0 = \calBe p(c_1) &= \sum_{\alpha,\beta}{\frac{\partial^2p}{\partial x_\alpha\partial x_\beta}(c_1)(P(e_\alpha,e_\beta)c_1|{\bf e})}+\lambda\sum_\alpha{\frac{\partial p}{\partial x_\alpha}(c_1)e_\alpha}\\
 &= \frac{\partial^2p}{\partial x_1^2}(c_1)+\frac{1}{2}\sum_{j=2}^r{\sum_{k=1}^d{\frac{\partial^2p}{\partial x_{jk}^2}(c_1)}}+\frac{d}{2}\sum_{j=1}^r{\frac{\partial p}{\partial x_j}(c_1)}.
\end{align*}
Since $\frac{\partial p}{\partial x_1}(c_1)=mp(c_1)$ and $\frac{\partial^2p}{\partial x_1^2}(c_1)=m(m-1)p(c_1)$ this gives
\begin{align*}
 \frac{1}{2}\sum_{j=2}^r{\sum_{k=1}^d{\frac{\partial^2p}{\partial x_{jk}^2}(c_1)}}+\frac{d}{2}\sum_{j=1}^r{\frac{\partial p}{\partial x_j}(c_1)} &= -m(m-1)p(c_1)
\end{align*}
which shows the claim.
\end{enumerate}
\item[\textup{(2)}] For $V=\Sym(2,\RR)$ we have $G=\Sp(2,\RR)/\{\pm\1\}$, $K^L=\SO(2)/\{\pm\1\}$ and $Z_G(\fraks)=O(2)/\{\pm\1\}$. It is easy to check that via the folding map $p:\RR^2\setminus\{0\}\to\Xi$ the space $\calH^m(\Xi)$ becomes $\calH^{2m}(\RR^2)$, the classical spherical harmonics of degree $2m$ on $\RR^2$. This is certainly an irreducible $O(2)$-representation which decomposes into two non-isomorphic irreducible $\SO(2)$-representations.\qedhere
\end{enumerate}
\end{proof}

\begin{theorem}\label{thm:DecompOfS}
The left-regular representation
 of $K^L$ on $L^2(\SS)$ decomposes 
 into a multiplicity-free direct sum 
 of irreducible representations of $K^L$. More precisely,
\begin{align*}
 L^2(\SS) &= {\sum_{m=0}^\infty}\raisebox{0.15cm}{\!$^\oplus$}{\calH^m(\SS)}.
\end{align*}
with $\calH^m(\SS)$ irreducible for $r>2$ or $d>1$ or $m=0$ and $\calH^m(\SS)$ decomposing into two non-isomorphic irreducible components for $r=2$, $d=1$ and $m>0$.
\end{theorem}

\begin{proof}
\begin{enumerate}
\item[\textup{(1)}] Let us first assume that $r>2$ or $d>1$ so that $\frakk^\frakl$ is semisimple. Since $\SS \simeq K^L/K_{c_1}^L$
 is a semisimple symmetric space,
 the space $L^2(\SS)$ is
 the multiplicity-free direct sum
 of all $K^L_{c_1}$-spherical $K^L$-representations
 by a theorem of E. Cartan. 
By Proposition \ref{prop:HighestWeightsSphericalReps} the only possible highest weights that can appear in $L^2(\SS)$ are given by
\begin{align*}
 \begin{cases}\{m\gamma:m\in\NN\} & \mbox{ for $d=1$,}\\\{2m\gamma:m\in\NN\} & \mbox{ for $d>1$.}\end{cases}
\end{align*}
If $d>1$ then by Proposition \ref{prop:SphericalHarmonicsReps} these representations in fact appear in $L^2(\SS)$ and hence their direct sum has to be dense in $L^2(\SS)$. The same argument applies to the case of $d=1$ (i.e. $V\simeq\Sym(n,\RR)$, $n\geq2$) where it only remains to show that the weights $m\gamma$ for $m\in\NN$ odd do not appear as highest weights of $K^L_{c_1}$-spherical representations. This is done in the next lemma.
\item[\textup{(2)}] In the case $V=\Sym(2,\RR)$ we have $\SS\simeq\PP^1(\RR)=S^1/\{\pm1\}$ and the decomposition is simply the expansion into Fourier coefficients since $\calH^m(\SS)\simeq\calH^{2m}(S^1)=\CC e^{2m\sqrt{-1}\theta}\oplus\CC e^{-2m\sqrt{-1}\theta}$ for $m>0$ and $\calH^0(\SS)\simeq\CC$.\qedhere
\end{enumerate}
\end{proof}

\begin{lemma}
Let $V=\Sym(k,\RR)$, $k\geq3$. Then the weights $m\gamma$ for $m\in\NN$ odd are not highest weights of $K^L_{c_1}$-spherical irreducible representations of $K^L$.
\end{lemma}

\begin{proof}
In this case $K^L=\PP\SO(k)$ acting by conjugation on $V=\Sym(k,\RR)$. Since $c_1=\diag(1,0\ldots,0)$, its stabilizer in $K^L$ is given by $K^L_{c_1}=S(O(1)\times O(k-1))$. It is known that the irreducible $\SO(k)$-representation of highest weight $m\gamma$ is the representation on the space $\calH^m(\RR^k)$ of spherical harmonics of degree $m$ in $\RR^k$. Obviously, the group $K^L_{c_1}$ fixes a non-zero vector in $\calH^m(\RR^k)$ if and only if $m$ is even which shows the claim.
\end{proof}

Now we also determine the spherical vectors and the highest weight vectors (with respect to a maximal torus in $\frakk^\frakl_\CC$ containing $\frakt_\CC$) in each $\calH^m(\SS)$.

\begin{proposition}\label{prop:HmSphericalVector}
Assume $r>2$ or $d>1$. For every $m\in\NN$ the space $\calH^m(\SS)^{K^L_{c_1}}$ of $K^L_{c_1}$-invariant harmonics of degree $m$ is one-dimensional and spanned by the function $\varphi_m\in\calH^m(\SS)$ given by
\begin{align*}
 \varphi_m(x) &= {_2F_1}(-m,m+r\lambda-1;\lambda;(x|c_1)), & x\in\SS,
\end{align*}
where $r$ is the rank of $V$ (or the rank of the symmetric space $G/K$) and $\lambda=\frac{d}{2}$ is the smallest non-zero discrete Wallach point (see \eqref{eqn:Wal1}).
\end{proposition}

\begin{remark}
The function ${_2F_1}(-m,m+r\lambda-1;\lambda;z)$ is
 a polynomial of $z$ of degree $m$, 
 which can be expressed
 in terms of the Jacobi polynomials $P_n^{(\alpha,\beta)}(z)$:
\begin{align*}
 {_2F_1}(-m,m+r\lambda-1;\lambda;z) &= \frac{m!\,\Gamma(\lambda)}{\Gamma(m+\lambda)}P_m^{(\lambda-1,(r-1)\lambda-1)}(1-2z).
\end{align*}
In the case $V=\Sym(k,\RR)$, 
 we have $\lambda=\frac{d}{2}=\frac{1}{2}$
 and $r=k$
 and the spherical vector reduces
\begin{align*}
 {_2F_1}(-m,m+r\lambda-1;\lambda;z^2) &= (-1)^m\frac{(2m)!\,(r\lambda-\frac{1}{2})_m}{(\frac{1}{2})_m(2r\lambda-2)_m(m+r\lambda-\frac{1}{2})_m}C_{2m}^{r\lambda-1}(z),
\end{align*}
by the change
 of the coordinates $z \mapsto z^2$
 (see \eqref{eq:2F1asCnlambda} in the Appendix).  
Here $C_n^\lambda(z)$ denotes the Gegenbauer polynomial. 
This corresponds to the well-known fact
 that 
 if $f\in C^{\infty}(S^{n-1})$
 is invariant by $O(n-1)$
 acting on the last $n-1$ coordinates and
 satisfies
 $\Delta_{S^{n-1}} f = -k (k+n-2)f$,
 then $f$ is a scalar multiple
 of the Gegenbauer polynomial $C_k^{\frac{n-2}{2}}(x_1)$,
 $(x_1,\ldots,x_n)\in S^{n-1}$.
\end{remark}

\begin{proof}[{Proof of Proposition \ref{prop:HmSphericalVector}}]
Since $\varphi_m$ is clearly $K^L_{c_1}$-invariant, it only remains to show that the $m$-homogeneous extension
\begin{align*}
 \overline{\varphi}_m(x) &= \tr(x)^m\varphi_m(\textstyle\frac{x}{|x|}), & x\in V,
\end{align*}
to a polynomial $\overline{\varphi}_m\in\calP^m(V)$ is harmonic. Note that $\overline{\varphi}_m(x)=\tr(x)^mu(\frac{(x|c_1)}{(x|{\bf e})})$ with $u(z)={_2F_1}(-m,m+r\lambda-1;\lambda;z)$. Then
\begin{align*}
 \calBe\overline{\varphi}_m(x) ={}& \calBe\tr(x)^m\cdot u\left(\frac{(x|c_1)}{(x|{\bf e})}\right)+2\left(\left.P\left(\frac{\partial\,\tr^m}{\partial x}(x),\frac{\partial}{\partial x}u\left(\frac{(x|c_1)}{(x|{\bf e})}\right)\right)x\right|{\bf e}\right)\\
 & +\tr(x)^m\cdot\calBe\left[u\left(\frac{(x|c_1)}{(x|{\bf e})}\right)\right].
\end{align*}
We have
\begin{align*}
 \frac{\partial\,\tr^m}{\partial x}(x) &= m\tr(x)^{m-1}{\bf e}, & \calBe\tr(x)^m &= m(r\lambda+m-1)\tr(x)^{m-1}.
\end{align*}
Therefore
\begin{align*}
 \left(\left.P\left(\frac{\partial\,\tr^m}{\partial x}(x),\frac{\partial}{\partial x}u\left(\frac{(x|c_1)}{(x|{\bf e})}\right)\right)x\right|{\bf e}\right) &= m\tr(x)^{m-1}\left(x\left|\frac{\partial}{\partial x}\right.\right)u\left(\frac{(x|c_1)}{(x|{\bf e})}\right)\\
 &= m\tr(x)^{m-1}\calE\left[u\left(\frac{(x|c_1)}{(x|{\bf e})}\right)\right],
\end{align*}
where $\calE=(x|\frac{\partial}{\partial x})$ is the Euler operator. Since the Euler operator is a radial operator and $u\left(\frac{(x|c_1)}{(x|{\bf e})}\right)$ is invariant under dilations, this term vanishes and we find
\begin{align*}
 \calBe\overline{\varphi}_m(x) &= m(r\lambda+m-1)\tr(x)^{m-1}u\left(\frac{(x|c_1)}{(x|{\bf e})}\right)+\tr(x)^m\calBe\left[u\left(\frac{(x|c_1)}{(x|{\bf e})}\right)\right].
\end{align*}
To calculate the action of $\calBe$ on $u\left(\frac{(x|c_1)}{(x|{\bf e})}\right)$, let $(e_\alpha)$ be an orthonormal basis of $V$ and denote by $x_\alpha$ the coordinates of $x\in V$ with respect to this basis. Then
\begin{align*}
 \frac{\partial}{\partial x_\alpha}u\left(\frac{(x|c_1)}{(x|{\bf e})}\right) ={}& \left(\frac{(e_\alpha|c_1)}{(x|{\bf e})}-\frac{(e_\alpha|{\bf e})(x|c_1)}{(x|{\bf e})^2}\right)u'\left(\frac{(x|c_1)}{(x|{\bf e})}\right),\\
 \frac{\partial^2}{\partial x_\alpha\partial x_\beta}u\left(\frac{(x|c_1)}{(x|{\bf e})}\right) ={}& \left(\frac{(e_\alpha|c_1)}{(x|{\bf e})}-\frac{(e_\alpha|{\bf e})(x|c_1)}{(x|{\bf e})^2}\right)\left(\frac{(e_\beta|c_1)}{(x|{\bf e})}-\frac{(e_\beta|{\bf e})(x|c_1)}{(x|{\bf e})^2}\right)u''\left(\frac{(x|c_1)}{(x|{\bf e})}\right)\\
 & +\left(2\frac{(e_\alpha|{\bf e})(e_\beta|{\bf e})(x|c_1)}{(x|{\bf e})^3}-\frac{(e_\alpha|c_1)(e_\beta|{\bf e})}{(x|{\bf e})^2}-\frac{(e_\beta|c_1)(e_\alpha|{\bf e})}{(x|{\bf e})^2}\right)u'\left(\frac{(x|c_1)}{(x|{\bf e})}\right).
\end{align*}
Hence,
\begin{align*}
 & \calBe\left[u\left(\frac{(x|c_1)}{(x|{\bf e})}\right)\right]\\
 ={}& \left(\frac{(P(c_1,c_1)x|{\bf e})}{(x|{\bf e})^2}-2\frac{(x|c_1)(P(c_1,{\bf e})x|{\bf e})}{(x|{\bf e})^3}+\frac{(x|c_1)^2(P({\bf e},{\bf e})x|{\bf e})}{(x|{\bf e})^4}\right)u''\left(\frac{(x|c_1)}{(x|{\bf e})}\right)\\
 & +\left(2\frac{(x|c_1)(P({\bf e},{\bf e})x|{\bf e})}{(x|{\bf e})^3}-2\frac{(P(c_1,{\bf e})x|{\bf e})}{(x|{\bf e})^2}\right)u'\left(\frac{(x|c_1)}{(x|{\bf e})}\right)\\
 & +\lambda\left(\frac{(c_1|{\bf e})}{(x|{\bf e})}-\frac{({\bf e}|{\bf e})(x|c_1)}{(x|{\bf e})^2}\right)u'\left(\frac{(x|c_1)}{(x|{\bf e})}\right)\\
 ={}& \frac{1}{(x|{\bf e})}\left[\frac{(x|c_1)}{(x|{\bf e})}\left(1-\frac{(x|c_1)}{(x|{\bf e})}\right)u''\left(\frac{(x|c_1)}{(x|{\bf e})}\right)+\lambda\left(1-r\frac{(x|c_1)}{(x|{\bf e})}\right)u'\left(\frac{(x|c_1)}{(x|{\bf e})}\right)\right]\\
 ={}& -\frac{1}{(x|{\bf e})}m(r\lambda+m-1)u\left(\frac{(x|c_1)}{(x|{\bf e})}\right)
\end{align*}
by the differential equation for the hypergeometric function. Hence, $\calBe\overline{\varphi}_m=0$.
\end{proof}

Let $\widetilde{\frakt}_\CC\subseteq\frakk^\frakl_\CC$ be a maximal torus containing $\frakt_\CC$. Abusing notation, we denote by $\gamma$ also the element in $\widetilde{\frakt}_\CC^*$ which vanishes on $\frakk^\frakl_{c_1}$ and equals $\gamma$ on $(\frakk^\frakl_{c_1})^\perp$. We choose an ordering on $\Phi(\frakk^\frakl_\CC,\widetilde{\frakt}_\CC)$ such that $\gamma$ is positive.

\begin{proposition}\label{prop:HmHighestWeightVector}
Assume $r>2$ or $d>1$. For every $m\in\NN$ the function $\phi_m$ on $\SS$ defined by
\begin{align*}
 \phi_m(x) &:= (x|c_1+\sqrt{-1}x_0-c_2)^m, & x\in\SS,
\end{align*}
is a highest weight vector with respect to $\widetilde{\frakt}_\CC$ in $\calH^m(\SS)$.
\end{proposition}

\begin{proof}
Let us write $a=c_1+\sqrt{-1}x_0-c_2$ for short. We divide the proof into two steps.
\begin{enumerate}
\item[\textup{(1)}] \textbf{Claim:} $\phi_m\in\calH^m(\SS)$. Clearly $\phi_m$ defines a homogeneous polynomial of degree $m$ by the same formula $\phi_m(x)=(x|a)^m$. Since $\frac{\partial}{\partial x}(x|a)=a$ we obtain
\begin{align*}
 \calBe\phi_m(x) &= m(m-1)(x|a)^{m-2}(P(a)x|{\bf e})+m\lambda(x|a)^{m-1}(a|{\bf e})\\
 &= m(m-1)(x|a)^{m-2}(x|a^2)+m\lambda(x|a)^{m-1}\tr(a) = 0
\end{align*}
since $a^2=0$ and $\tr(a)=0$. Hence, $\phi_m\in\calH^m(\XX)$.
\item[\textup{(2)}] \textbf{Claim:} $\phi_m$ is a weight vector of weight $2m\gamma$. Note that
\begin{align*}
 X_0a &= X_0c_1+\sqrt{-1}X_0x_0-X_0c_2 = -\frac{x_0}{4}+\sqrt{-1}\frac{c_1-c_2}{2}-\frac{x_0}{4} = \frac{1}{2}\sqrt{-1}a.
\end{align*}
Hence,
\begin{align*}
 X_0\cdot\phi_m(x) &= -D_{X_0x}\phi_m(x) = -\left(X_0x\left|\frac{\partial\phi_m}{\partial x}(x)\right.\right)\\
 &= -m(x|a)(X_0x|a) = m(x|a)(x|X_0a)\\
 &= \frac{m}{2}\sqrt{-1}(x|a)^m = 2m\gamma(X_0)\phi_m(x).
\end{align*}
\end{enumerate}
Hence, $\phi_m\in\calH^m(\SS)_{2m\gamma}$. Since $\calH^m(\SS)$ is a highest weight representation with highest weight $2m\gamma$,
 the claim follows.
\end{proof}

\subsection{Branching laws with respect to ${\mathfrak{sl}}(2,\RR)$}
\label{sec:L2branching}

We recall from \eqref{eqn:sl2}
 that there is a distinguished subalgebra
 $\fraks\simeq {\mathfrak {sl}}(2,\RR)$
 of $\frakg = \co(V)$.  
Let $S^{\vee}= SL(2,\RR)^{\vee}$
denote the connected subgroup 
 of $G^{\vee}$
 with Lie algebra $\fraks \simeq {\mathfrak{sl}}(2,\RR)$.  
We shall prove
 that our \minholrep\ splits discretely
 into a direct sum 
 of irreducible unitary representations
 when we restrict it
 to the subalgebra $\fraks$.   
For each irreducible representation $\pi_t$
 of $SL(2,\RR)^{\vee}$, 
 the multiplicity space 
 $\operatorname{Hom}_{SL(2,\RR)^{\vee}}(\pi_t,L^2(\Xi))$
 is described
 by the space of generalized spherical harmonics, 
 on which the compact subgroup $K^L$ of $G$ acts naturally. We will even see that the possibly disconnected compact subgroup $Z_{G^\vee}(\fraks)\subseteq G^\vee$ with Lie algebra $\frakk^\frakl$ acts on the space of generalized spherical harmonics.
Thus we shall find the branching law
 of the \minholrep\
 with respect to ${SL(2,\RR)}^{\vee} \times Z_{G^\vee}(\fraks)$
 in Theorem \ref{thm:sl2decoS}
 for the Schr{\"o}dinger 
 model and in Theorem \ref{thm:DualPairDecompositionFock}
 for the Fock model.

\nocite{K98}
\nocite{KM05}

\subsubsection*{Discretely decomposable restriction}

First of all,
we give a quick review of a general theory
 of discretely decomposable  representations.

We begin with a general setting.
Let $\frakg'=\frakk'+\frakp'$ be a Cartan decomposition
 of a semisimple Lie algebra over $\RR$.  
The following notion singles out an algebraic property of unitary
representations that split into irreducible representations
without continuous spectra.
\begin{definition}
\label{def:DDR}
Let $(\varpi,X)$ be a $(\mathfrak{g}',\frakk')$-module.  
\begin{enumerate}
\item[\textup{(1)}]
We say $\varpi$ is $\frakk'$-admissible if\/
$\dim\operatorname{Hom}_{\frakk'}(\tau,\varpi) < \infty$
 for any irreducible,
 finite dimensional representation
 $\tau$ of $\frakk'$.
\item[\textup{(2)}]
{\upshape (\cite[Definition 1.1]{K98})}
We say $\varpi$ is a discretely decomposable
 if there exist an increasing sequence of\/
$(\mathfrak{g}',\frakk')$-modules $\{X_j\}_{j \in {\mathbb{N}}}$
 of finite length
 such that
\begin{equation*}
X = \bigcup_{j=0}^\infty X_j,
\end{equation*}
\item[\textup{(3)}]
We say $\varpi$ is infinitesimally unitary
 with respect to a Hermitian inner product $(\ ,\ )$
 on $X$ if 
\[
(\varpi(Y)u,v) = -(u,\varpi(Y)v)
\quad\text{for any $Y \in \mathfrak{g}'$ and any $u,v\in X$}.
\]
\end{enumerate}
\end{definition}
We collect some basic results on discretely decomposable
 $(\mathfrak{g}',\frakk')$-modules: 
\begin{fact}
[see {\cite{K98}}]\label{fact:deco}
Let $(\varpi,X)$ be a $(\mathfrak{g}',\frakk')$-module.
\begin{enumerate}
\item[\textup{(1)}]
If $\varpi$ is $\frakk'$-admissible,
then $\varpi$ is discretely decomposable
 as a $(\mathfrak{g}',\frakk')$-modules.
\item[\textup{(2)}]
Suppose $\varpi$ is discretely decomposable as a
$(\mathfrak{g}',\frakk')$-module. 
If $\varpi$ is infinitesimally unitary,
 then $\varpi$ is isomorphic to
an algebraic direct sum
 of irreducible $(\mathfrak{g}',\frakk')$-modules.
\end{enumerate}
\end{fact}

The point here 
 is that $(\varpi, X)$ is not necessarily
 of finite length
 as a $(\frakg',\frakk')$-module.  
We shall apply this concept to the specific situation where
$\mathfrak{g}'=\mathfrak{sl}(2,{\mathbb{R}})$.  

\subsubsection*{The dual pair}

Recall that $\fraks\simeq\sl(2,\RR)$ is the subalgebra
 of $\frakg=\co(V)$
 spanned by $E$, $F$ and $H$ and that the compact group $K^L$ has Lie algebra $\frakk^\frakl=\aut(V)$.

\begin{proposition}\label{prop:DualPair}
$(\fraks,\aut(V))$ is a reductive dual pair in $\co(V)$.
\end{proposition}

\begin{proof}
Let $Z_{\co(V)}(\fraks)$ denote
 the centralizer of $\fraks$ in $\co(V)$.  
We first show that $Z_{\co(V)}(\fraks)=\aut(V)$. An element $(u,T,v)\in\co(V)$ is in the centralizer of $\fraks$ in $\co(V)$ if and only if the following three equations hold:
\begin{align}
 0 &= [(u,T,v),E] = (T{\bf e},-2{\bf e}\Box v,0)\label{eq:uTvCentEqE}\\
 0 &= [(u,T,v),H] = (-u,0,v)\label{eq:uTvCentEqH}\\
 0 &= [(u,T,v),F] = (0,2u\Box{\bf e},-T^\#{\bf e}).\label{eq:uTvCentEqF}
\end{align}
Equation \eqref{eq:uTvCentEqH}
 implies immediately
 that $u=0=v$
 and \eqref{eq:uTvCentEqE}
 yields $T{\bf e}=0$
 which is equivalent to $T\in\aut(V)$. 
In this case,
 all equations are satisfied
 and we obtain $Z_{\co(V)}(\fraks)=\aut(V)$.

Conversely let us prove
 that $Z_{\co(V)}(\aut(V))=\fraks$. We have $(u,T,v)\in Z_{\co(V)}(\aut(V))$ if and only if
\begin{align*}
 0 &= [(u,T,v),(0,S,0)] = (-Su,[T,S],S^\# v) & \forall\,S\in\aut(V).
\end{align*}
First, by the next lemma we obtain that $u,v\in\RR{\bf e}$. It remains to show that $T\in\RR\,\id$. Write $T=L(x)+D$. Then $[T,S]=0$ for all $S\in\aut(V)$ implies $Sx=0$ and $[S,D]=0$ for all $S\in\aut(V)$. Again by the next lemma, we obtain $x\in\RR{\bf e}$, so it remains to show that $D=0$. Let $Z(\frakl)$ denote the center of $\frakl$. We know that $Z(\frakl)=\RR\,\id_V$ and $\frakl=[\frakl,\frakl]\oplus Z(\frakl)$ with $[\frakl,\frakl]$ semisimple. Since $\aut(V)$ is generated by the derivations $[L(x),L(y)]$, $x,y\in V$, we have $\aut(V)\subseteq[\frakl,\frakl]$. 
Therefore $\aut(V)$ is a symmetric subalgebra of the semisimple Lie algebra $[\frakl,\frakl]$ and hence semisimple itself. Thus, $[S,D]=0$ for all $S\in\aut(V)$ implies $D=0$ and the proof is complete.
\end{proof}

\begin{lemma}
Let $x\in V$. If $Dx=0$ for all $D\in\aut(V)$ then $x\in\RR{\bf e}$.
\end{lemma}

\begin{proof}
Write $x=\sum_{i\leq j}{x_{ij}}$, $x_{ij}\in V_{ij}$. For convenience we put $x_{ji}:=x_{ij}$ for $i<j$. For $D=[L(c_i),L(y)]$, $y\in V_{ij}$, $i\neq j$, we obtain
\begin{align*}
 0 &= Dx = c_i(xy)-y(c_ix)\\
 &= c_i\left(\sum_{k\neq j}{x_{ik}y}+\sum_{k\neq i}{x_{kj}y}+x_{ij}y\right)-y\left(x_{ii}+\frac{1}{2}\sum_{k\neq i}{x_{ik}}\right)\\
 &= \frac{1}{2}x_{ii}y+\frac{1}{2}\sum_{k\neq i}{x_{kj}y}+\frac{(x_{ij}|y)}{2}c_i-x_{ii}y-\frac{1}{2}\sum_{k\neq i}{x_{ik}y}.
\end{align*}
The $c_i$-component of this expression is $\frac{1}{2}(x_{ij}|y)c_i$ which has to vanish for every $y\in V_{ij}$. Since $\tau$ is non-degenerate, this means that $x_{ij}=0$ for all $i\neq j$. Hence, $x=\sum_{i=1}^{r_0}{t_ic_i}$. From the above calculation we obtain
\begin{align*}
 0 &= \frac{t_i}{4}y+\frac{t_j}{4}y-\frac{t_i}{2}y = \frac{t_j-t_i}{4}y
\end{align*}
which implies $t_i=t_j$. Since this has to hold for all $i,j=1,\ldots,r_0$, we obtain $x\in\RR{\bf e}$.
\end{proof}

\subsubsection*{The $\sl(2,\RR)$-representations}

The integral formula \eqref{eq:IntFormulaO}
 with respect to the polar decomposition 
$
     \Xi\simeq \RR_+ \times \SS
$
 leads us to the isomorphism
 of the Hilbert space:
\begin{align*}
 L^2(\Xi,\td\mu) \simeq L^2(\RR_+,t^{r\lambda-1}\td t)
\widehat \otimes L^2(\SS).
\end{align*}
By using Theorem \ref{thm:DecompOfS} 
 we obtain:
\begin{align}
\label{eqn:pideco}
  L^2(\Xi,\td\mu)
 \simeq 
{\sum_{m=0}^\infty}\raisebox{0.15cm}{\!$^\oplus$}{L^2(\RR_+,t^{r\lambda-1}\td t)\otimes\calH^m(\SS)}.
\end{align}
We show in Theorem \ref{thm:sl2decoS}
that this is the decomposition of the representation $\pi$
 into irreducible $S^{\vee} \times Z_{G^\vee}(\fraks)$-representations. 
The proof essentially coincides with Kobayashi--Mano
 \cite[Section 1.3]{KM07a} 
 in the case $V=\Sym(k,\RR)$ and $V=\RR^{1,k-1}$
 or in the deformed setting \cite[Theorem 3.28]{BKO09},
 but we use the Jordan algebra
 to carry out necessary computations.

First of all,
 we review a realization of 
 a series $\widetilde{\pi}_s$ ($-1<s<\infty$)
 of lowest weight representations of the universal covering group of $\SL(2,\RR)$ on $L^2(\RR_+)$ 
(see B.~Kostant \cite{Kos00}
or an earlier work by Ranga Rao \cite{RR77}). 
Following \cite{Kos00}
 (but the vectors $e$ and $f$ are replaced
 by $-e$ and $-f$
 which amounts to a Lie algebra isomorphism),
 we define the differential action $\td\widetilde{\pi}_s$
 of the Lie algebra $\sl(2,\RR)$ on $L^2(\RR_+)$
 by the following skew-adjoint operators, $t$ denoting the variable in $\RR_+$:
\begin{align*}
 \td\widetilde{\pi}_s(e) &= \sqrt{-1}t,\\
 \td\widetilde{\pi}_s(h) &= 2t\frac{\td}{\td t}+1,\\
 \td\widetilde{\pi}_s(f) &= \sqrt{-1}\left(t\frac{\td^2}{\td t^2}+\frac{\td}{\td t}-\frac{s^2}{4t}\right).
\end{align*}
Further,
 the underlying $(\frakg,\frakk)$-module is spanned
 over $\CC$
 by the functions
\begin{align*}
 \widetilde{\phi}_k^s(t) &= t^{\frac{s}{2}}e^{-t}L_k^s(2t), & k\in\NN,
\end{align*}
where $L_n^\alpha(z)$ denote the Laguerre polynomials.

Let $\mu \in \RR$
 (we shall take a specific $\mu$ later), 
 and transfer these representations to $L^2(\RR_+,t^\mu\td t)$, 
through the unitary isomorphism
\begin{align*}
 \calU:L^2(\RR_+)\to L^2(\RR_+,t^\mu\td t),\,\calU f(t)=t^{-\frac{\mu}{2}}f(t).
\end{align*}
Define the representation $\pi_s$ on $L^2(\RR_+,t^\mu\td t)$ by 
\[
     \pi_s:=\calU\circ\widetilde{\pi}_{s-1}\circ\calU^{-1}. 
\]
(Note that the parameterization of $\pi_s$ follows \cite{FK94} and the parametrization of $\widetilde{\pi}_s$
follows \cite{Kos00}. The parameterizations are such that $s$ is the lowest weight of $\pi_s\simeq\widetilde{\pi}_{s-1}$.) Then its differential representation is given by
\begin{align*}
 \td\pi_s(e) &= \sqrt{-1}t,\\
 \td\pi_s(h) &= 2t\frac{\td}{\td t}+(\mu+1),\\
 \td\pi_s(f) &= \sqrt{-1}\left(t\frac{\td^2}{\td t^2}+(\mu+1)\frac{\td}{\td t}-\frac{(s-1)^2-\mu^2}{4t}\right).
\end{align*}
The underlying $({\mathfrak {sl}}(2,\RR), {\mathfrak {so}}(2))$-module
 $\calV_s$ is spanned by the functions
\begin{align*}
 \phi_k^s(t) &:= \calU\widetilde{\phi}_k^{s-1}(t) = t^{\frac{s-\mu-1}{2}}e^{-t}L_k^s(2t), & k\in\NN.
\end{align*}
We obtain representations $(\pi_s,L^2(\RR_+,t^\mu\td t))$ of $\widetilde{\SL(2,\RR)}$ for $s\in(0,\infty)$ with underlying Lie algebra modules $(\td\pi_s,\calV_s)$.

Now put $\mu:=r\lambda-1$ and $s:=r\lambda+2m$, $m\in\NN$. We collect some additional information on the resulting representations. The action of the inverse Cayley transformed $\sl_2$-triple $(\widetilde{e},\widetilde{f},\widetilde{h})$ (see Subsection \ref{sec:SchrödingerForMinimalHolomorphic}) on the basis $\phi_k^s$ is given by
\begin{align*}
 \td\pi_s(\widetilde{e})\phi_k^s &= 2\sqrt{-1}\phi_{k+1}^s,\\
 \td\pi_s(\widetilde{h})\phi_k^s &= (r\lambda+2m+2k)\phi_k^s,\\
 \td\pi_s(\widetilde{f})\phi_k^s &= \tfrac{1}{2}k(r\lambda+2m+k-1)\sqrt{-1}\phi_{k-1}^s.
\end{align*}
where for convenience we put $\phi_{-1}^s:=0$. Hence, the vector $\phi_0^s(t)=t^{2m}e^{-t}$ is a lowest weight vector, i.e. $\td\pi_s(\widetilde{f})\phi_0^s=0$. The lowest weight is given by $r\lambda+2m$.

For each $m \in {\mathbb{N}}$, 
 we define a linear map $\Phi_m$
 by 
\begin{align}
 \Phi_m:&
        L^2(\RR_+,t^{r\lambda-1}\td t)
        \otimes
        \calH^m(\SS)
        \to 
        L^2(\Xi,\td\mu),\,
\label{eqn:Phim}
\\
        &\Phi_m(f\otimes\phi)(x)=f(|x|)\phi(\textstyle\frac{x}{|x|}).  
\notag
\end{align}
This constructs
 each summand in \eqref{eqn:pideco}, 
 and $\Phi_m$ respects the actions
of $\widetilde{\SL(2,\RR)}\times Z_{G^\vee}(\fraks)$ 
as follows:

\begin{theorem}
\label{thm:sl2decoS}
\begin{enumerate}
\item[\textup{(1)}] The linear map 
$
     \Phi_m
$
 respects the action of ${\mathfrak{sl}}(2,\RR)$, 
 that is, 
\begin{align*}
 \td\pi(E)\circ\Phi_m &= \Phi_m\circ(\td\pi_{r\lambda+2m}(e)\otimes\id),\\
 \td\pi(H)\circ\Phi_m &= \Phi_m\circ(\td\pi_{r\lambda+2m}(h)\otimes\id),\\
 \td\pi(F)\circ\Phi_m &= \Phi_m\circ(\td\pi_{r\lambda+2m}(f)\otimes\id).
\end{align*}
\item[\textup{(2)}]
The $(\sl(2,\RR),\so(2))$-module
 $\Phi_m(\calV_{r\lambda+2m}\otimes\calH^m(\SS))$
 is contained in the space $L^2(\Xi)_{\frakk}$
 of $\frakk$-finite vectors of $L^2(\Xi,\td\mu)$. 
\item[\textup{(3)}]
The representation $(\pi,L^2(\Xi,\td\mu))$
 decomposes into a multiplicity-free
 sum
 of irreducible representations
 of $\widetilde{SL(2,\RR)} \times Z_{G^\vee}(\fraks)$
 as follows:
\[
  L^2(\Xi,\td\mu)
 \simeq 
 {\sum_{m=0}^\infty}\raisebox{0.15cm}{$^\oplus$}{\,\pi_{r \lambda + 2m} \boxtimes\calH^m(\SS)},
\]
and its underlying $(\frakg,\frakk)$-module
 decomposes
 under the action of the dual pair $(\fraks, \frakk^{\frakl})$ as 
\begin{align}
 L^2(\Xi,\td\mu)_\frakk 
 &\simeq
 \bigoplus_{m=0}^\infty
 {\calV_{r\lambda+2m}\boxtimes\calH^m(\SS)},
\label{eq:DualPairDecompositionSchrödinger}
\end{align}
where $\calV_{s}$ is the irreducible representation of $\fraks\simeq\sl(2,\RR)$ of lowest weight $s$ and $\calH^m(\SS)$ is an irreducible $\frakk^\frakl$-module for $r>2$ or $d>1$ or $m=0$ and decomposes into two irreducible non-isomorphic $\frakk^\frakl$-modules for $r=2$, $d=1$ and $m>0$.
\end{enumerate}
\end{theorem}

\begin{remark}
\begin{enumerate}
\item[\textup{(1)}]
Suppose $G'$ is a reductive subgroup of $G$.  
In general ${\mathfrak {k}}'$-finite vectors
 are not necessarily
 ${\mathfrak {k}}$-finite vectors
 in the irreducible unitary representation $\pi$
 of $G$.  
The first statement
 of Theorem \ref{thm:sl2decoS}
 constructs a discrete part of the branching law
 of the restriction $\pi|_{G'}$, 
 and the second statement implies
 that there is no continuous spectrum
 \cite{K98}.
 \item[\textup{(2)}]
For $\frakg=\so(2,k)$ and $\frakg={\mathfrak {sp}}(k,\RR)$
 this decomposition was given in \cite[Section 1.3]{KM07a}
 and extended to a Dunkl setting in \cite{BSO06}. 
See also \cite[Theorem 7.1]{KO03b}
 for the branching law 
 for the minimal representation of $O(p,q)$
to the subgroup $O(p,q_1) \times O(q_2)$
 when $p+q$ is even 
 and $q_1 + q_2 = q$.  
\end{enumerate}
\end{remark}

\begin{proof}[Proof of Theorem \ref{thm:sl2decoS}]
\begin{enumerate}
\item[\textup{(1)}]
Let $f\otimes\phi\in L^2(\RR_+,t^{\frac{rd}{2}-1}\td t)\otimes\calH^m(\SS)$. Then
\begin{align*}
 \td\pi(E)(\Phi_m(f\otimes\phi))(x) &= \sqrt{-1}\tr(x)f(|x|)\phi(\textstyle\frac{x}{|x|})\\
 &= \Phi_m(\td\pi_{r\lambda+2m}(e)f\otimes\phi)(x),\\
 \td\pi(H)(\Phi_m(f\otimes\phi))(x) &= (2D_x+r\lambda)\left[f(|x|)\phi(\textstyle\frac{x}{|x|})\right]\\
 &= \left(2|x|f'(|x|)+r\lambda f(|x|)\right)\phi(\textstyle\frac{x}{|x|})\\
 &= \Phi_m(\td\pi_{r\lambda+2m}(h)f\otimes\phi)(x).
\end{align*}
For the action of $F$ we use Lemma \ref{lem:BesselProdRule} to obtain
\begin{align*}
 & \td\pi(F)(\Phi_m(f\otimes\phi))(x)\\
 ={}& \sqrt{-1}\calBe\left[f(|x|)\phi(\textstyle\frac{x}{|x|})\right]\\
 ={}& \sqrt{-1}\left(\calBe f(|x|)\cdot\phi(\textstyle\frac{x}{|x|})+2\left(\left.P\left(\frac{\partial}{\partial x}f(|x|),\frac{\partial}{\partial x}\phi(\textstyle\frac{x}{|x|})\right)x\right|{\bf e}\right)\right.\\
 & \left.+f(|x|)\cdot\calBe\phi(\textstyle\frac{x}{|x|})\right).
\end{align*}
By using the formulae
\begin{align*}
 \calBe f(|x|) &= |x|f''(|x|)+r\lambda f'(|x|),\\
 \frac{\partial}{\partial x}f(|x|) &= f'(|x|){\bf e},\\
 \frac{\partial}{\partial x}\phi(\textstyle\frac{x}{|x|}) &= \frac{\partial}{\partial x}\left[\tr(x)^{-m}\phi(x)\right]\\
 &= -m|x|^{-1}\phi(\textstyle\frac{x}{|x|}){\bf e}+|x|^{-m}\frac{\partial\phi}{\partial x}(x), 
\end{align*}
we find that
\begin{align*}
 & \left(\left.P\left(\frac{\partial}{\partial x}f(|x|),\frac{\partial}{\partial x}\phi(\textstyle\frac{x}{|x|})\right)x\right|{\bf e}\right)\\
 ={}& f'(|x|)\left[-m|x|^{-1}\phi(\textstyle\frac{x}{|x|})(P({\bf e},{\bf e})x|{\bf e})+|x|^{-m}\left(\left.P\left({\bf e},\frac{\partial\phi}{\partial x}(x)\right)x\right|{\bf e}\right)\right]\\
 ={}& f'(|x|)\left[-m\phi(\textstyle\frac{x}{|x|})+|x|^{-m}\left(x\left|\frac{\partial\phi}{\partial x}(x)\right.\right)\right]\\
 ={}& f'(|x|)\left[-m\phi(\textstyle\frac{x}{|x|})+m|x|^{-m}\phi(x)\right] = 0.
\end{align*}
Further, for the last term we have with Lemma \ref{lem:BesselOnProductWithTrace}:
\begin{align*}
 \calBe\phi(\textstyle\frac{x}{|x|}) &= \calBe\left[\tr(x)^{-m}\phi(x)\right]\\
 &= -m(r\lambda+m-1)|x|^{-m-1}\phi(x)
 &= -m(r\lambda+m-1)|x|^{-1}\phi(\textstyle\frac{x}{|x|}).
\end{align*}
Putting things together gives
\begin{align*}
 & \td\pi(F)(\Phi(f\otimes\phi))(x)\\
 ={}& \sqrt{-1}\left(|x|f''(|x|)+r\lambda f'(|x|)-m(r\lambda+m-1)|x|^{-1}f(|x|)\right)\phi(\textstyle\frac{x}{|x|})\\
 ={}& \Phi_m(\td\pi_{r\lambda+2m}(f)f\otimes\phi)(x).  
\end{align*}
\item[\textup{(2)}]
We recall from \eqref{eqn:L2K}
that
\[
  L^2(\Xi)_{\frakk}=\calP(\Xi)e^{-|x|}.  
\]

Since $\phi_0^{r\lambda+2m}(t)=t^{2m}e^{-t}$ is a lowest weight vector
 in $\calV_{r\lambda+2m}$
 and 
\[
     \Phi_m(\phi_0^{r\lambda+2m}
     \otimes
     \calH^m(\SS))\in\calP(\Xi)e^{-|x|}=L^2(\Xi)_{\frakk},
\]
 the second assertion holds. 
\item[\textup{(3)}]
Note first that each summand $\calV_{r\lambda+2m}\boxtimes\calH^m(\SS)$ is an $\fraks$-isotypic component and hence not only $K^L$ but also the possibly disconnected group $Z_{G^\vee}(\fraks)$ acts on $\calH^m(\SS)$. Since $L^2(\Xi)_{\frakk}$ is discretely decomposable
 as an $\fraks\oplus \frakk^{\frakl}$-module,
 taking $\frakk$-finite vectors in \eqref{eqn:pideco}
 now yields the third statement.\qedhere
\end{enumerate}
\end{proof}

\begin{remark}
\label{rem:discdeco}
In the above proof,
 we have used a specific fact
 on the Schr{\"o}dinger model,
 namely,
 $L^2(\Xi)_{\frakk}$
 coincides with $\calP(\Xi)e^{-|x|}$.  
Alternatively,
 we can use a representation theoretic result,
 namely,
 Fact \ref{fact:deco} (1).  
In fact,
 any lowest weight $(\frakg,\frakk)$-module
 is ${\mathfrak{z}}(\frakk)$-admissible
 where ${\mathfrak{z}}(\frakk)$ denotes
 the center of $\frakk$.  
In our setting ${\mathfrak{z}}(\frakk)=\RR(-E+F)
 \subset {\mathfrak{sl}}(2,\RR)$, 
 and thus the assumption
 of Fact \ref{fact:deco}
 is fulfilled.   
\end{remark}

\subsection{Folding maps and the Schr\"odinger model}
\label{subsec:foldingmaps}
The classical Schr{\"o}dinger model 
 for the Weil representation
 of the metaplectic group $\Mp(n,\RR)$
 is realized on $L^2(\RR^n)$,
 whereas our Schr{\"o}dinger model 
 for $\frakg= {\mathfrak {sp}}(n,\RR)$
 is realized in a somewhat different space,
 namely,
 the Hilbert space $L^2(\Xi)$
 where $\Xi$ is the manifold
 consisting of symmetric matrices
 of rank one.  
In this subsection 
 we relate these two Hilbert spaces
 by the {\it{folding map}}
 (see \eqref{eqn:fold} below),
 and also explain the irreducible decompositions
 in Theorem \ref{thm:sl2decoS}
 not only for $\frakg= {\mathfrak {sp}}(k,\RR)$
 but also for the cases $\frakg=\su(k,k),\so^*(4k)$.  
Let $V=\Herm(k,\FF)$,
 $\FF=\RR,\CC$, or $\HH$.  
We let $d=\dim_\RR\FF=2\lambda$ is $1$, $2$ or $4$. 
The minimal $L$-orbit $\Xi$ is given as the image
 of the folding map 
\begin{equation}
\label{eqn:fold}
     p:\FF^k\setminus\{0\}\to\Xi,\,
       x\mapsto xx^*,
\end{equation}
where $x^*={}^t\!\overline{x}$ denotes the composition of transposition and conjugation. The folding map $p$ is a principal bundle with principal fiber $U(1;\FF)$ and hence
\begin{align*}
 p^*:L^2(\Xi,\td\mu)\stackrel{\sim}{\to}L^2(\FF^k)^{U(1;\FF)}.
\end{align*}
Correspondingly to 
$
    \operatorname{Herm}(k,\FF)
    \subset
    \operatorname{Herm}(kd,\RR)
    =\operatorname{Sym}(kd,\RR),
$
 there is a natural homomorphism 
 $\frakg \to {\mathfrak {sp}}(kd,\RR)$.  
The isomorphism $p^*$ intertwines
 the representation $\pi$
 of $G^{\vee}$
 on $L^2(\Xi,\td\mu)$
 with the restriction of the Weil representation
 of $\Mp(kd,\RR)$ on $L^2(\RR^{kd})\simeq L^2(\FF^k)$. 
We recall the dual pair correspondence
 with respect to 
\[
  \widetilde{SL(2,\RR)} \times O(n) 
 \to\Mp(n,\RR)
\]
amounts to a multiplicity-free decomposition
 of the Schr{\"o}dinger model 
 as a representation
 of $\widetilde{\SL(2,\RR)} \times O(n)$:
\begin{align*}
 L^2(\RR^n) 
 \simeq 
 {\sum_{j=0}^\infty}\raisebox{0.15cm}{$^\oplus$}
 {\pi_{j+\frac{n}{2}}\boxtimes\calH^j(\RR^n)}.
\end{align*}
Now we take $U(1;\FF)$-invariants on both sides and obtain:
\begin{enumerate}
\item[\textup{(1)}] $\frakg=\sp(k,\RR)$. We have $U(1;\FF)=O(1)=\{\pm\1\}$ and hence the $O(1;\FF)$-invariants are exactly the terms with sperical harmonics of even degree. This yields
\begin{align*}
 L^2(\Xi,\td\mu) &\simeq L^2(\RR^k)^{O(1)} 
                    \simeq {\sum_{m=0}^\infty}\raisebox{0.15cm}{$^\oplus$}
                    {\pi_{2m+\frac{k}{2}}\boxtimes\calH^{2m}(\RR^k)}
\end{align*}
which is \eqref{eq:DualPairDecompositionSchrödinger} since $r\lambda=\frac{k}{2}$.
\item[\textup{(2)}] $\frakg=\su(k,k)$. We have $U(1;\FF)=U(1)$ and in the decomposition
\begin{align*}
 \calH^j(\RR^{2k}) &= \bigoplus_{\alpha+\beta=j}{\calH^{\alpha,\beta}(\CC^k)}
\end{align*}
into $U(k)$-irreducibles the $U(1)$-invariants are exactly those $\calH^{\alpha,\beta}(\CC^k)$ with $\alpha=\beta$. 
This yields
\begin{align*}
 L^2(\Xi,\td\mu) &\simeq L^2(\RR^{2k})^{U(1)} 
 \simeq {\sum_{m=0}^\infty}\raisebox{0.15cm}{$^\oplus$}{\pi_{2m+k}\boxtimes\calH^{m,m}(\CC^k)}
\end{align*}
which is \eqref{eq:DualPairDecompositionSchrödinger} since $r\lambda=k$.
\item[\textup{(3)}] $\frakg=\so^*(4k)$. We have $U(1;\FF)=Sp(1)$ and in the decomposition
\begin{align*}
 \calH^j(\RR^{4k}) &= \bigoplus_{\substack{p+q=j\\p\geq q\geq0}}{\calH^{p,q}(\HH^k)\boxtimes\CC^{p-q+1}}
\end{align*}
into $Sp(k)$-irreducibles the $Sp(1)$-invariants are exactly those $\calH^{p,q}(\HH^k)$ with $p=q$. This yields
\begin{align*}
 L^2(\Xi,\td\mu) &\simeq L^2(\RR^{4k})^{Sp(1)} 
 \simeq {\sum_{m=0}^\infty}\raisebox{0.15cm}{$^\oplus$}{\pi_{2m+2k}\boxtimes\calH^{m,m}(\HH^k)}
\end{align*}
which is \eqref{eq:DualPairDecompositionSchrödinger} since $r\lambda=2k$.
\end{enumerate}
\section{A Fock space realization
 for \minholreps}
\label{sec:fockspace}

In this section we construct a Fock space $\calF(\XX)$
 on a complex submanifold $\XX$
 in the complex Jordan algebra $V_{\CC}$
 defined in Subsection \ref{sec:MinRepCplx},
 which is biholomorphic to the minimal nilpotent
 $K_{\CC}$-orbit $\OO_{\operatorname{min}}^{K_{\CC}}$
 in $\frakp_{\CC}$
 via the Cayley transform.  
For this we introduce a density
 on $\XX$ given explicitly
 by a K-Bessel function,
 and define an action of the conformal Lie algebra
 $\frakg$ on the space $\calP(\XX)$
 of regular functions.  
A remarkable feature
 is that the action is given not
 by pseudodifferential operators
 (cf. \cite{BK94})
 but by polynomial differential operators
 up to second order.  
We then find the reproducing kernel 
 of the Fock space $\calF(\XX)$, 
 and give a proof
 of the irreducibility and unitarizability
 of the $(\frakg,\frakk)$-module
 $\calP(\XX)$
 by using the {\it{Bessel operators}}.  
In the next section we see that the two representations 
 on $L^2(\Xi,\td \mu)$
 (Section \ref{sec:schmodel})
 and on $\calF(\XX)$
 (Section \ref{sec:fockspace})
 are isomorphic to each other.  

\subsection{Polynomials on $\XX$}

Recall that $\XX$ is the minimal non-zero $L_\CC$-orbit in $V_\CC$ which is a complexification of the real $L$-orbit $\Xi$
 through a primitive idempotent $c_1$
 in the Euclidean Jordan algebra $V$. 
Let $\Delta_j$ denote the principal minors of $V$ with respect to a fixed Jordan frame $c_1,\ldots,c_r$. For ${\bf m}\in\NN^r$ we define
\begin{align*}
 \Delta_{\bf m}(x) &:= \Delta_1(x)^{m_1-m_2}\cdots\Delta_{r-1}(x)^{m_{r-1}-m_r}\Delta_r(x)^{m_r}.
\end{align*}
All these polynomials are extended holomorphically to $V_\CC$. Denote by $\calP(V_\CC)$ the space of all holomorphic polynomials on $V_\CC$. We further let $\calP^{\bf m}(V_\CC)$ denote the subspace of $\calP(V_\CC)$ spanned by the polynomials $\Delta_{\bf m}(gx)$, $g\in L_\CC$.
The following decomposition
 is a Jordan theoretic reformulation
 of the Hua--Kostant--Schmid theorem
in the tube case:

\begin{theorem}
The space $\calP(V_\CC)$ of holomorphic polynomials on $V_\CC$ decomposes
 into a multiplicity-free sum of irreducible $L_\CC$-modules:
\begin{align*}
 \calP(V_\CC) &= \bigoplus_{\bf m}{\calP^{\bf m}(V_\CC)}.
\end{align*}
\end{theorem}

\begin{remark}\label{rem:KHSonXX}
Note that if $z\in\XX$, then $\Delta_2(z)=\ldots=\Delta_r(z)=0$ and hence $\Delta_{\bf m}\neq0$ on $\XX$ iff $m_2=\ldots=m_r=0$. 
Therefore the space $\calP(\XX)$ of restrictions
 of holomorphic polynomials on $V_\CC$ to $\XX$ decomposes under the $L_\CC$-action as
\begin{align*}
 \calP(\XX) &= \bigoplus_{m=0}^\infty{\calP^m(\XX)},\\
\intertext{where}
 \calP^m(\XX) &= \{p|_\XX:p\in\calP^{(m,0,\ldots,0)}(V_\CC)\}.
\end{align*}
In fact,
 it is clear that the restriction map $\calP^{(m,0,\ldots,0)}(V_\CC)\to\calP^m(\XX)$ is an isomorphism of $L_\CC$-modules
 since it is a surjective $L_\CC$-homomorphism
 and $\calP^{(m,0,\ldots,0)}(V_\CC)$ is irreducible.  
\end{remark}

\subsection{Construction of the Fock space}

We introduce a density $\omega$ on $\XX$ by
\begin{align*}
 \omega(z) &= \widetilde{K}_{\lambda-1}(|z|) & z\in\XX,
\end{align*}
where $|z|:=(z|\overline{z})^\frac{1}{2}$ 
and $\widetilde K_{\alpha}(z)
=(\frac z 2)^{-\alpha}K_{\alpha}(z)$
 is the renormalized K-Bessel function. 
In view of the integral formula \eqref{eq:IntFormulaX} and the asymptotic behaviour of the K-Bessel function (see Appendix \ref{app:BesselFcts}), the $L^2$-inner product
\begin{align*}
 \langle F,G\rangle &:= \int_\XX{F(z)\overline{G(z)}\omega(z)\td\nu(z)}
\end{align*}
is finite for any $F,G\in\calP(\XX)$
 and hence turns $\calP(\XX)$ into a pre-Hilbert space. We denote its completion by $\widetilde{\calF}(\XX)$. Let $\calO(\XX)$ be the space of holomorphic functions on the complex manifold $\XX$. In Theorem \ref{thm:NaturalFockSpace} we will prove that the space $\widetilde{\calF}(\XX)$ coincides with the Fock space (as defined in the introduction)
\begin{align}
 \calF(\XX) &= \left\{F\in\calO(\XX):\int_\XX{|F(z)|^2\omega(z)\td\nu(z)}<\infty\right\}.\label{eq:DefFockSpace}
\end{align}

\begin{proposition}
$\calF(\XX)$ is a closed subspace of $L^2(\XX,\omega\td\nu)$ and the point evaluation $\calF(\XX)\to\CC,\,F\mapsto F(z)$ is continuous for every $z\in\XX$. In particular, $\widetilde{\calF}(\XX)\subseteq\calF(\XX)$ and the point evaluation $\widetilde{\calF}(\XX)\to\CC,\,F\mapsto F(z)$ is continuous for every $z\in\XX$.
\end{proposition}

\begin{proof}
This is a local statement and hence, we may transfer it with a chart map to an open domain $\Omega\subseteq\CC^k$. Here the measure $\omega\td\nu$ is absolutely continuous with respect to the Lebesgue measure $\td z$ and hence it suffices to show that $\calO(\CC^d)\cap L^2(\CC^d,\td z)\subseteq L^2(\CC^d,\td z)$ is a closed subspace with continuous point evaluations. This is done e.g. in \cite[Proposition 3.1 and Corollary 3.2]{Hel62}.
\end{proof}

We recall
 that $\calB$ is a vector-valued holomorphic differential operator $\calB$
 introduced in \eqref{eq:DefHolBessel}.  
Then the density $\omega$
 satisfies the following.  
\begin{lemma}
\label{lem:Bw}
$
\calB \omega (z) = \frac{\overline z}{4}
\omega (z).  
$
\end{lemma}

\begin{proof}
In view that $\omega(z)=u(\frac{(z|\overline{z})}{4})$
 where $u(t)=\widetilde K_{\lambda-1} (2 \sqrt t)$, 
Lemma follows from Proposition \ref{prop:BesselHypergeomEq}.  
\end{proof}

\begin{proposition}
\label{prop:BesselAdjoint}
The adjoint $\calB^*$ of $\calB$ on $\calF(\XX)$
 is the multiplication operator by $\frac{z}{4}$.
\end{proposition}

\begin{proof}
Let $F, G \in \calF(\XX)$. 
Then by Proposition \ref{prop:PiCC} we know that
\begin{align*}
 \int_\XX{\calB F(z)\overline{G(z)}\omega(z)\td\nu(z)} &= \int_\XX{F(z)\calB(\overline{G(z)}\omega(z))\td\nu(z)}.
\end{align*}
The function $\overline{G(z)}$ is antiholomorphic and hence $\frac{\partial}{\partial z}\overline{G(z)}=0$. Using Lemma \ref{lem:BesselProdRule} we obtain
\begin{align*}
 \int_\XX{F(z)\calB(\overline{G(z)}\omega(z))\td\nu(z)} &= \int_\XX{F(z)\overline{G(z)}\calB\omega(z)\td\nu(z)}.
\end{align*}
 Now Proposition \ref{prop:BesselAdjoint} follows from
 Lemma \ref{lem:Bw}.  
\end{proof}

\subsection{The Bessel--Fischer inner product}

We introduce another inner product on the space $\calP(\XX)$ of polynomials, namely the \textit{Bessel--Fischer inner product}. For two polynomials $p$ and $q$ it is defined by
\begin{align*}
 [p,q] &:= \left.p(\calB)\overline{q}(4z)\right|_{z=0},
\end{align*}
where $\overline{q}(z)=\overline{q(\overline{z})}$ is obtained by conjugating the coefficients of the polynomial $q$. A priori it is not even clear that this sesquilinear form is positive definite.

\begin{proposition}\label{prop:FischerEqualsL2}
For $p,q\in\calP(\XX)$ we have
\begin{align}
 [p,q] &= \langle p,q\rangle.\label{eq:FischerEqualsL2}
\end{align}
\end{proposition}

The proof is similar to the proof of \cite[Proposition 3.8]{BSO06}

\begin{proof}
First note that for all $p,q\in\calP(\XX)$
\begin{align*}
 [(a|\textstyle\frac{z}{4})p,q] &= [p,(\overline{a}|\calB)q] & \mbox{for }a\in V_\CC,\\
 \langle(a|\textstyle\frac{z}{4})p,q\rangle &= \langle p,(\overline{a}|\calB)q\rangle & \mbox{for }a\in V_\CC.
\end{align*}
In fact, the second equation follows from Proposition \ref{prop:BesselAdjoint}. The first equation is immediate since the components $(a|\calB)$, $a\in V_\CC$, of the Bessel operator form a commuting family of differential operators on $\XX$. 
Therefore $(a|\calB)p(\calB)\overline{q}(4z)=4p(\calB)\overline{(\overline{a}|\calB)q}(4z)$ and the claim follows. To prove \eqref{eq:FischerEqualsL2} we proceed by induction on $\deg(q)$. First, if $p=q=\1$, the constant polynomial with value $1$, it is clear that $[p,q]=1$. With the integral formula \eqref{eq:IntFormulaX} we further find that
\begin{align*}
 c_{r,\lambda}\langle p,q\rangle &= c_{r,\lambda}\int_\XX{\omega(z)\td\nu(z)} = \int_0^\infty{\widetilde{K}_{\lambda-1}(t)t^{2r\lambda-1}\td t}\\
 &= 2^{2r\lambda-2}\Gamma\left(r\lambda\right)\Gamma\left((r-1)\lambda+1\right) = c_{r,\lambda}.
\end{align*}
where we have used the integral formula \eqref{eq:KBesselIntFormula} for the last equality. Thus, \eqref{eq:FischerEqualsL2} holds for $\deg(p)=\deg(q)=0$. If now $\deg(p)$ is arbitrary and $\deg(q)=0$ then $(\overline{a}|\calB)q=0$ and hence
\begin{align*}
 [(a|\textstyle\frac{z}{4})p,q] &= [p,(\overline{a}|\calB)q] = 0 & \mbox{and}\\
 \langle(a|\textstyle\frac{z}{4})p,q\rangle &= \langle p,(\overline{a}|\calB)q\rangle = 0.
\end{align*}
Therefore \eqref{eq:FischerEqualsL2} holds if $\deg(q)=0$. We note that \eqref{eq:FischerEqualsL2} also holds if $\deg(p)=0$ and $\deg(q)$ is arbitrary. In fact,
\begin{align*}
 [p,q] &= p(0)\overline{q(0)} = \overline{[q,p]}, & \mbox{and} && \langle p,q\rangle &= \overline{\langle q,p\rangle}
\end{align*}
and \eqref{eq:FischerEqualsL2} follows from the previous considerations. Now assume \eqref{eq:FischerEqualsL2} holds for $\deg(q)\leq k$. For $\deg(q)\leq k+1$ we then have $\deg((\overline{a}|\calB)q)\leq k$ and hence, by the assumption
\begin{align*}
 [(a|\textstyle\frac{z}{4})p,q] &= [p,(\overline{a}|\calB)q] = \langle p,(\overline{a}|\calB)q\rangle = \langle(a,\textstyle\frac{z}{4})p,q\rangle.
\end{align*}
This shows \eqref{eq:FischerEqualsL2} for $\deg(q)\leq k+1$ and $p(0)=0$, i.e. without constant term. But for constant $p$, i.e. $\deg(p)=0$ we have already seen that \eqref{eq:FischerEqualsL2} holds and therefore the proof is complete.
\end{proof}

The previous theorem provides us with a new expression for the inner product on $\widetilde{\calF}(\XX)$. The Bessel--Fischer inner product is more suitable for explicit computations. If we denote by $\calP^m(\XX)\subseteq\calP(\XX)$ the subspace of homogeneous polynomials of degree $m\in\NN$, then the following result is immediate with the Bessel--Fischer inner product.

\begin{corollary}\label{cor:PmOrthogonal}
The subspaces $\calP^m(\XX)$ are pairwise orthogonal.
\end{corollary}

\begin{proof}
Let $p\in\calP^m(\XX)$ and $q\in\calP^n(\XX)$ with $m\neq n$. We may assume without loss of generality that $m<n$. It is clear that for $a\in V_\CC$ we have $(a|\calB)\overline{q}\in\calP^{n-1}(\XX)$. 
Therefore $p(\calB)\overline{q}(4z)\in\calP^{n-m}(\XX)$. Since $m-n\neq0$, every polynomial in $\calP^{n-m}(\XX)$ vanishes at $z=0$ and hence, $[p,q]=0$.
\end{proof}

\begin{proposition}\label{prop:FockSpaceGlobalFcts}
\begin{align*}
 \left\{F|_\XX:F\in\calO(V_\CC),\int_\XX{|F(z)|^2\omega(z)\td\nu(z)}<\infty\right\} &\subseteq \widetilde{\calF}(\XX).
\end{align*}
\end{proposition}

\begin{proof}
Let $F\in\calO(V_\CC)$. Then $F$ has a Taylor expansion $F(z)=\sum_{m=0}^\infty{p_m(z)}$ into homogeneous polynomials $p_m$ of degree $m$ which converges uniformly on bounded subsets. We show that this series also converges in $\calF(\XX)$. Then, since point evaluation in $\calF(\XX)$ is continuous, it follows that $F$ as the limit of this series is also in $\calF(\XX)$.\\
For $R>0$ we put $\XX_R:=\{z\in\XX:|z|\leq R\}$. Since $\XX_R$ is bounded, the series $\sum_{m=0}^\infty{p_m}$ converges uniformly on $\XX_R$. Hence, we obtain
\begin{align*}
 \infty &> \int_\XX{|F(z)|^2\omega(z)\td\nu(z)}\\
 &= \lim_{R\to\infty}{\int_{\XX_R}{|F(z)|^2\omega(z)\td\nu(z)}}\\
 &= \lim_{R\to\infty}{\sum_{m,n=0}^\infty{\int_{\XX_R}{p_m(z)\overline{p_n(z)}\omega(z)\td\nu(z)}}}.
\end{align*}
We claim that for $m\neq n$ and $R>0$ we have
\begin{align*}
 \int_{\XX_R}{p_m(z)\overline{p_n(z)}\omega(z)\td\nu(z)} &= 0.
\end{align*}
In fact, we already know that $\langle p_m,p_n\rangle=0$. With the integral formula \eqref{eq:IntFormulaX} we find
\begin{align*}
 0 &= \int_{K^{L_\CC}}{\int_0^\infty{p_m(utc_1)\overline{p_n(utc_1)}\omega(utc_1)t^{2r\lambda-1}\td t}\td u}\\
 &= \int_{K^{L_\CC}}{p_m(uc_1)\overline{p_n(uc_1)}\td u}\cdot\int_0^\infty{\omega(tc_1)t^{m+n+2r\lambda-1}\td t}.
\end{align*}
Using the integral formula \eqref{eq:KBesselIntFormula} we find that the second factor is a product of Gamma functions and never vanishes. 
Therefore the first factor has to vanish. But then the same calculation yields
\begin{multline*}
 \int_{\XX_R}{p_m(z)\overline{p_n(z)}\omega(z)\td\nu(z)}\\
 = \int_{K{L_\CC}}{p_m(uc_1)\overline{p_n(uc_1)}\td u}\cdot\int_0^R{\omega(tc_1)t^{m+n+2r\lambda-1}\td t} = 0.
\end{multline*}
We then obtain that
\begin{align*}
 \lim_{R\to\infty}{\sum_{m=0}^\infty{\int_{\XX_R}{|p_m(z)|^2\omega(z)\td\nu(z)}}} < \infty.
\end{align*}
Now we can interchange the limits since the right hand side converges absolutely. This yields
\begin{align*}
 \sum_{m=0}^\infty{\|p_m\|^2} &< \infty
\end{align*}
which is nothing else but the convergence of the series $\sum_{m=0}^\infty{p_m}$ in $\calF(\XX)$.
\end{proof}

\subsection{The reproducing kernel}

We can now calculate the reproducing kernel of the Hilbert space $\widetilde{\calF}(\XX)$. For this we first calculate the reproducing kernels on the finite-dimensional subspaces $\calP^m(\XX)$.

\begin{proposition}\label{prop:RepKernelm}
The reproducing kernel $\KK^m(z,w)$ of the Hilbert space $\calP^m(\XX)$ is given by
\begin{align*}
 \KK^m(z,w) &= \frac{1}{4^mm!(\lambda)_m}(z|\overline{w})^m, & z,w\in\XX.
\end{align*}
\end{proposition}

\begin{proof}
Write $\KK^m_w(z)=\KK^m(z,w)$. We use the Fischer inner product to show that $p(w)=[p,\KK_w^m]=p(\calB)\KK_{\overline{w}}^m(4z)|_{z=0}$ for any polynomial $p\in\calP^m(\XX)$. For this we note that by \eqref{eq:BesselOnPowers} the action of the Bessel operator on the polynomials $(z|w)^k$ is given by
\begin{align*}
 \calB_z(z|w)^k &= k(\lambda+k-1)(z|w)^{k-1}w.
\end{align*}
Now suppose $p$ is a monomial, i.e. $p(z)=\prod_{j=1}^m{(a_j|z)}$ with $a_j\in V_\CC$. Iterating \eqref{eq:BesselOnPowers} we obtain
\begin{align*}
 [p,(-|\overline{w})^m] &= \left.p(\calB)(4z|w)^m\right|_{z=0} = 4^m\left.\left(\prod_{j=1}^m{(a_j|\calB)}\right)(z|w)^m\right|_{z=0}\\
 &= 4^m\Big(m(m-1)\cdots1\Big)\Big((\lambda+m-1)(\lambda+m-2)\cdots\lambda\Big)\prod_{j=1}^m{(a_j|w)}\\
 &= 4^mm!(\lambda)_mp(w).\qedhere
\end{align*}
\end{proof}

We give a closed formula
 of the reproducing kernel of $\widetilde{\calF}(\XX)$
 in terms of the renormalized I-Bessel function
 $\widetilde{I}_\alpha(z)=(\frac z 2)^{-\alpha}I_{\alpha}(z)$
  (see Appendix \ref{app:BesselFcts}).

\begin{theorem}
\label{thm:Frepro}
The reproducing kernel $\KK(z,w)$ of the Hilbert space $\widetilde{\calF}(\XX)$ is given by
\begin{align*}
 \KK(z,w) &= \Gamma(\lambda)\widetilde{I}_{\lambda-1}(\sqrt{(z|\overline{w})}), & z,w\in\XX.
\end{align*}
\end{theorem}

\begin{proof}
By the previous result
 $\KK^m(z,w)$ is the reproducing kernel
 of $\calP^m(\XX)$. 
Further we know by Corollary \ref{cor:PmOrthogonal} that the spaces $\calP^m(\XX)$ are pairwise orthogonal. 
Therefore by \cite[Proposition I.1.8]{Nee00}, the sum
\begin{align*}
 \sum_{m=0}^\infty{\KK^m(z,w)}
 = \sum_{m=0}^\infty{\frac{1}{4^mm!(\lambda)_m}(z|\overline{w})^m}
 = \Gamma(\lambda)\widetilde{I}_{\lambda-1}(\sqrt{(z|\overline{w})}).
\end{align*}
converges pointwise to the reproducing kernel $\KK(z,w)$ of the direct Hilbert sum $\calP(\XX)=\bigoplus_{m=0}^\infty{\calP^m(\XX)}$.
\end{proof}

The following consequence is a standard result for reproducing kernel spaces and can e.g. be found in \cite[page 9]{Nee00}.

\begin{corollary}
For every $F\in\widetilde{\calF}(\XX)$ and every $z\in\XX$ we have
\begin{align*}
 |F(z)| \leq \KK(z,z)^{\frac{1}{2}}\|F\|.
\end{align*}
\end{corollary}

\subsection{Unitary action on the Fock space}

In Subsection \ref{sec:MinRepCplx} we have already verified that the complexification $\td\pi_\CC$ of the action $\td\pi$ defines a Lie algebra representation on $C^\infty(\XX)$ by polynomial differential operators in $z$. 
Thus the action $\td \pi_{\CC}$
 preserves the subspace $\calP(\XX)$
 of holomorphic polynomials.
We shall define the action,
 to be denoted by $\td \rho$,
 as the conjugation
 of $\td\pi_\CC$
 by the Cayley transform $c \in \operatorname{Int}(\frakg_{\CC})$
 introduced in \eqref{eq:DefCayleyTransform}. 

\begin{definition}
On $\calP(\XX)$ we define a $\frakg$-action $\td\rho$ by
\begin{equation}
\label{eqn:drho}
  \td \rho := \td \pi_{\CC} \circ c.  
\end{equation}
By the formulae
 \eqref{eq:CayleyTransform1}, 
 \eqref{eq:CayleyTransform2}
 and \eqref{eq:CayleyTransform3}, 
 this definition amounts to 
\begin{align*}
 \td\rho(a,0,0) &= \textstyle\td\pi_\CC(\frac{a}{4},\sqrt{-1}L(a),a),\\
 \td\rho(0,L(a)+D,0) &= \textstyle\td\pi_\CC(\sqrt{-1}\frac{a}{4},D,-\sqrt{-1}a),\\
 \td\rho(0,0,a) &= \textstyle\td\pi_\CC(\frac{a}{4},-\sqrt{-1}L(a),a), 
\end{align*}
for $a \in V$
 and $D \in \frak{aut}(V)$
 in terms of the Jordan algebra.  
\end{definition}

\begin{remark}\label{rem:FockKAction}
By Lemma \ref{lem:Cayley}, 
the decomposition of $(\td\rho,\calP(\XX))$ into $\frakk$-types
 equals the decomposition of $(\td\pi_\CC,\calP(\XX))$ into $\frakl$-types. 
The action of $\frakl$ under $\td\pi$ is induced by the geometric action of $L$ on the orbit $\XX=L_\CC\cdot c_1$ up to multiplication by a character. 
In particular,
 $\calP(\XX)$ is $\frakk$-finite via $\td \rho$.  
In view of Remark \ref{rem:KHSonXX},
 $\calP(\XX)$ decomposes into $\frakk$-types as follows:
\begin{align*}
 \calP(\XX) &= \bigoplus_{m=0}^\infty{\calP^m(\XX)}.
\end{align*}
The unique (up to scalar) $\frakk^\frakl$-invariant vector
 in the $\frakk$-type $\calP^m(\XX)$ is the $m$-th power of the trace:
\begin{align*}
 \Psi_m(z) &:= \tr(z)^m.
\end{align*}
\end{remark}

The $\sl_2$-triple 
$
     (\widetilde{E},\widetilde{F},\widetilde{H})
    =(c^{-1}E, c^{-1}F, c^{-1}H)
$ acts on $\calP(\XX)$ by
\begin{align*}
 \td\rho(\widetilde{E}) &= \td\pi_\CC(E) = \sqrt{-1}\tr(z),\\
 \td\rho(\widetilde{H}) &= \td\pi_\CC(H) = 2\calE+r\lambda\\
 \td\rho(\widetilde{F}) &= \td\pi_\CC(F) = \sqrt{-1}\calBe,
\end{align*}
where $\calE=(x|\frac{\partial}{\partial x})$ is the Euler operator and $\calBe$ the identity component of the Bessel operator (see \eqref{eq:IdentityBessel}). Since $\calP^m(\XX)$ consists
 of homogeneous polynomials of degree $m$, 
 we have 
\begin{align}
\label{eqn:HPm}
 \td\rho(\widetilde{H})|_{\calP^m(\XX)} &= 2m+r\lambda.
\end{align}
Then it is easy
 to see the following mapping properties of the $\sl_2$ acting
 on the $\frakk$-types:
\begin{align}
 \td\rho(\widetilde{E}): \calP^m(\XX)\to\calP^{m+1}(\XX),
\label{eqn:EPm}
\\
 \td\rho(\widetilde{F}): \calP^m(\XX)\to\calP^{m-1}(\XX).  
\label{eqn:FPm}
\end{align}
Let us compute
 how they act on $\Psi_m$:
\begin{lemma}\label{lem:ActionOnPowersOfTrace}
For $z\in\XX$ and $m\in\NN$ we have
\begin{align*}
 \td\rho(\widetilde{E})\tr(z)^m &= \sqrt{-1}\tr(z)^{m+1},\\
 \td\rho(\widetilde{H})\tr(z)^m &= (r\lambda+2m)\tr^m(z),\\
 \td\rho(\widetilde{F})\tr(z)^m &= m(r\lambda+m-1)\sqrt{-1}\tr(z)^{m-1}.  \\
\end{align*}
\end{lemma}

\begin{proof}
Since $\tr(z)=(z|{\bf e})$ and $\frac{\partial}{\partial z}(z|{\bf e})={\bf e}$, we have
\begin{align*}
 \calB_{\bf e}\tr(z)^m &= m(m-1)\tr(z)^{m-2}(P({\bf e})z|{\bf e})+m\lambda\tr(z)^{m-1}({\bf e}|{\bf e})\\
 &= m(\textstyle\frac{rd}{2}+m-1)\tr(z)^{m-1}.\\
\end{align*}
The other statement for $\td \rho(\widetilde{E})$
 and $\td \rho(\widetilde H)$ are clear.  
\end{proof}
We are ready
 to give two basic properties
 of the $(\frakg,\frakk)$-module
 $\calP(\XX)$ via $\td \rho$,
 namely, 
 Proposition \ref{prop:irredPX} and Lemma \ref{lem:ActionOnPowersOfTrace}.  

\begin{proposition}
\label{prop:irredPX}
$\calP(\XX)$ is an irreducible $(\frakg,\frakk)$-module.
\end{proposition}

\begin{proof}
By Remark \ref{rem:FockKAction}, 
the $\frakk$-type decomposition of $\calP(\XX)$ is given by
\begin{align*}
 \calP(\XX) &= \bigoplus_{m=0}^\infty{\calP^m(\XX)}.
\end{align*}
Therefore it suffices to show
 that for each $m\in\NN$ there exists a vector $v\in\calP^m(\XX)$ and $X,Y\in\frakg_\CC$
 such that $0\neq\td\rho(X)v\in\calP^{m-1}(\XX)$
 ($m \ge 1$)
 and $0\neq\td\rho(Y)v\in\calP^{m+1}(\XX)$. 
But this follows immediately
{} from Lemma \ref{lem:ActionOnPowersOfTrace}.
\end{proof}

\begin{proposition}
\label{prop:PXinfuni}
The $(\frakg,\frakk)$-module $\calP(\XX)$ is infinitesimally unitary with respect to the $L^2$-inner product $\langle-,-\rangle$.
\end{proposition}

\begin{proof}
Using Proposition \ref{prop:BesselAdjoint} it is immediate that $\td\rho(a,0,a)$ and $\td\rho(0,L(a),0)$, $a\in V$, act on $\calP(\XX)$ by skew-symmetric operators. It remains to consider $\td\rho(a,D,-a)$ for $a\in V$ and $D\in\frakk^\frakl=\aut(V)$. Then $\td\rho(a,D,-a)=\td\pi_\CC(0,D+2\sqrt{-1}L(a),0)$.
\begin{enumerate}
\item[\textup{(a)}] We first treat the case of $\td\rho(0,D,0)=\td\pi_\CC(0,D,0)$, $D\in\frakk^\frakl$. Using Proposition \ref{prop:PiCC} we have
\begin{align*}
 \int_\XX{\td\rho(0,D,0)F(z)\cdot\overline{G(z)}\omega(z)\td\nu(z)} &= -\int_\XX{F(z)\cdot\overline{\td\rho(0,D,0)(G(z)\omega(z))}\td\nu(z)}.
\end{align*}
Since $\omega(z)$ is invariant under $K^L$ we have $\td\rho(0,D,0)\omega(z)=0$ and hence $\td\rho(0,D,0)$ is skew-symmetric.
\item[\textup{(b)}] Now consider $\td\rho(a,0,-a)=\td\pi_\CC(0,2\sqrt{-1}L(a),0)=2\sqrt{-1}\td\pi_\CC(0,L(a),0)$. It suffices to show that $\td\pi_\CC(0,L(a),0)=D_{az}+\frac{\lambda}{2}$ is symmetric with respect to the inner product in $L^2(\XX,\omega\td\nu)$. Using Proposition \ref{prop:PiCC} we find that
\begin{align*}
 & \int_\XX{\td\pi_\CC(0,L(a),0)F(z)\cdot\overline{G(z)}\omega(z)\td\nu(z)}\\
 ={}& -\int_\XX{F(z)\cdot\td\pi_\CC(0,L(a),0)(\overline{G(z)}\omega(z))\td\nu(z)}.\\
\intertext{Since $\overline{G(z)}$ is antiholomorphic, 
we have $D_{az}\overline{G(z)}=0$ and hence}
 ={}& -\int_\XX{F(z)\overline{G(z)}\cdot\td\pi_\CC(0,L(a),0)\omega(z)\td\nu(z)}.
\end{align*}
Now, $\td\pi_\CC(0,L(a),0)\omega(z)=\overline{\td\pi_\CC(0,L(a),0)\omega(z)}$. In fact, $\omega(z)=\phi(z|\overline{z})$ with $\phi\in C^\infty(\RR_+)$ real-valued and hence
\begin{align*}
 \overline{D_{az}\omega(z)} &= \overline{(az|\overline{z})\phi'(z|\overline{z})} = (a\overline{z}|z)\phi'(z|\overline{z})\\
 &= (\overline{z}|az)\phi'(z|\overline{z}) = D_{az}\omega(z).
\end{align*}
Therefore we obtain
\begin{align*}
 & \int_\XX{\td\pi_\CC(0,L(a),0)F(z)\cdot\overline{G(z)}\omega(z)\td\nu(z)}\\
 ={}& -\int_\XX{F(z)\overline{G(z)}\cdot\overline{\td\pi_\CC(0,L(a),0)\omega(z)}\td\nu(z)}\\
 ={}& -\overline{\int_\XX{\overline{F(z)}G(z)\cdot\td\pi_\CC(0,L(a),0)\omega(z)\td\nu(z)}}\\
\intertext{and the same argument as in the beginning, interchanging $F$ and $G$, shows that}
 ={}& \overline{\int_\XX{\overline{F(z)}\cdot\td\pi_\CC(0,L(a),0)G(z)\omega(z)\td\nu(z)}}\\
 ={}& \int_\XX{F(z)\cdot\overline{\td\pi_\CC(0,L(a),0)G(z)}\omega(z)\td\nu(z)}.
\end{align*}
Hence, also $\td\rho(a,0,-a)$ is skew-adjoint and the proof is complete.\qedhere
\end{enumerate}
\end{proof}

\begin{theorem}\label{thm:UnitaryRepOnFock}
The $(\frakg,\frakk)$-module $\calP(\XX)$ integrates to an irreducible unitary representation $\rho$ of the universal cover $\widetilde{G}$ of $G$ on $\widetilde{\calF}(\XX)$.
\begin{enumerate}
\item[\textup{(1)}] For $r>1$ this representation factors to a finite cover $G^\vee$ of $G$ given by $G^\vee:=\widetilde{G}/\Gamma$, where $\Gamma=\exp(k\pi\ZZ({\bf e},0,-{\bf e}))$ and $k\in\NN_+$ is an integer such that $k\frac{r\lambda}{2}=k\frac{rd}{4}\in\ZZ$.
\item[\textup{(2)}] For $r=1$ this representation factors to a finite cover of $G=\PP\SL(2,\RR)$ if and only if $\lambda\in\QQ$. In this case, a finite cover $G^\vee$ of $G$ to which $\calP(\XX)$ integrates is given by $G^\vee:=\widetilde{G}/\Gamma$, where $\Gamma=\exp(k\pi\ZZ({\bf e},0,-{\bf e}))$ and $k\in\NN_+$ is an integer such that $k\frac{r\lambda}{2}\in\ZZ$.
\end{enumerate}
\end{theorem}

\begin{proof}
By the previous results it only remains to check in which cases the minimal $\frakk$-type $\calP^0(\XX)=\CC\1$ integrates
 to a finite cover. Since the center of $\frakk$ is given by $Z(\frakk)=\RR({\bf e},0,-{\bf e})$ and $\frakk=Z(\frakk)\oplus[\frakk,\frakk]$ with $[\frakk,\frakk]$ semisimple, it suffices to check the action of $Z(\frakk)$. The $\frakk$-action on $\1$ is given by
\begin{align*}
 \td\rho(a,D,-a)\1 &= \td\pi_\CC(0,D+2\sqrt{-1}L(a),0)\1 = \frac{r\lambda}{2n}\Tr(2\sqrt{-1}L(a))\1\\
 &= \lambda\sqrt{-1}(a|{\bf e})\1.
\end{align*}
Therefore, the center $Z(\frakk)=\RR({\bf e},0,-{\bf e})$ acts by
\begin{align*}
 \td\rho({\bf e},0,-{\bf e})\1 &= r\lambda\sqrt{-1}\1.
\end{align*}
In $K$ we have $e^{\pi({\bf e},0,-{\bf e})}=\1$ and hence, the claim follows.
\end{proof}

\begin{remark}
The finite cover $G^\vee$ of $G$ constructed in Theorem \ref{thm:UnitaryRepOnFock} may not be minimal with the property that $\td\rho$ integrates to a representation of it. The minimal cover of $G$ to which $\td\rho$ integrates is determined in \cite[Theorem 2.30]{HKM11}.
\end{remark}

\begin{remark}
The reproducing kernel of the Fock space and the density $\omega(z)$ for the measure on $\XX$ only depend on $\lambda=\frac{d}{2}$ which is constant for the series $\frakg=\sp(k,\RR)$ ($d=1$), $\frakg=\su(k,k)$ ($d=2$) and $\frakg=\so^*(4k)$ ($d=4$) and therefore should give also a Fock model for the corresponding infinite-dimensional groups (see \cite{NO02} for the Schr\"odinger models).
\end{remark}

\subsection{Action of the $\sl_2$ and harmonic polynomials}

Let
\begin{align*}
 \calH^m(\XX) := \{p\in\calP^m(\XX):\calBe p=0\}
\end{align*}
be the space of \textit{harmonic polynomials} on $\XX$. 
Note that since $\Xi\subseteq\XX$ is totally real,
 the restriction to $\Xi$ defines an isomorphism $\calH^m(\XX)\to\calH^m(\Xi)$, where $\calH^m(\Xi)$ was defined
 in Subsection \ref{sec:SphericalHarmonics}. 
Therefore $\calH^m(\XX)$ is an irreducible representation of the compact group $Z_{G^\vee}(\fraks)$.

\begin{remark}
For $V=\RR$ the one-dimensional Jordan algebra
 we have $\XX=\CC^\times$
 and $\calBe=z\frac{\td^2}{\td z^2}+\lambda\frac{\td}{\td z}$. 
Hence, all solutions $u\in C^\infty(\XX)$ of $\calBe u=0$ are given by
\begin{align*}
 u(z) &= \begin{cases}C_1z^{1-\lambda}+C_2 & \mbox{ for $\lambda\neq1$,}\\C_1\ln(x)+C_2 & \mbox{ for $\lambda=1$,}\end{cases}
\end{align*}
$C_1,C_2\in\RR$. Since $\lambda>0$, only the constant functions are polynomial solutions
 and therefore $\calH^m(\XX)=0$ for $m>0$ and $\calH^0(\XX)=\CC\1$.
\end{remark}

\begin{proposition}\label{prop:PolyDecompIntoHarmonics}
Every polynomial $p\in\calP(\XX)$ decomposes uniquely into
\begin{align*}
 p &= \sum_{k=0}^m{\tr^k(z)h_{m-k}},
\end{align*}
where $h_{m-k}\in\calH^{m-k}(\XX)$ is a harmonic polynomial. The polynomials $h_{m-k}$ are explicitly given by
\begin{align*}
 h_{m-k} &= \sum_{j=0}^{m-k}{(-1)^j\frac{\Gamma(r\lambda+2m-2k)\Gamma(r\lambda+2m-2k-j-1)}{j!\,k!\,\Gamma(r\lambda+2m-k)\Gamma(r\lambda+2m-2k-1)}\tr^j(z)\calBe^{k+j}p}.
\end{align*}
In particular, we have
\begin{align*}
 \calP(\XX) &= \bigoplus_{k,m=0}^\infty{\tr^k(z)\calH^m(\XX)}.
\end{align*}
\end{proposition}

\begin{proof}
The proof works in the same way as in \cite[proof of Theorem 5.1]{BSO06}. We first show uniqueness by induction on $m$. For this assume
\begin{align}
 \sum_{k=0}^m{\tr^k(z)h_{m-k}} &= 0.\label{eq:UniquenessDecompEq}
\end{align}
for $h_{m-k}\in\calH^{m-k}(\XX)$. If $m=0$ then trivially $h_0=0$ and we are done. Now suppose $m>0$. Applying $\calBe^m$ to both sides yields, using Lemma \ref{lem:BesselOnProductWithTrace}:
\begin{align*}
 m!\,(r\lambda+m-1)\cdots(r\lambda)h_0 &= 0,
\end{align*}
whence $h_0=0$. Then \eqref{eq:UniquenessDecompEq} reads
\begin{align*}
 \sum_{k=0}^{m-1}{\tr^k(z)h_{m-k}} &= 0
\end{align*}
and the induction hypothesis applies.\\
To show the existence of the claimed decomposition as well as the explicit formula, we proceed in three steps:
\begin{enumerate}
\item[\textup{(1)}] Define
\begin{align*}
 Q_0p &:= \sum_{j=0}^m{q_{m,j}\tr^j(z)\calBe^jp},\\
\intertext{where}
 q_{m,j} &:= (-1)^j\frac{\Gamma(r\lambda+2m-j-1)}{j!\,\Gamma(r\lambda+2m-1)}.
\end{align*}
Then $Q_0p\in\calH^m(\XX)$. In fact, using Lemma \ref{lem:BesselOnProductWithTrace} we have
\begin{align*}
 \calBe Q_0p &= \sum_{j=0}^m{q_{m,j}\left(j(r\lambda+2(m-j)+j-1)\tr^{j-1}(z)\calBe^jp+\tr^j(z)\calBe^{j+1}p\right)}\\
 &= \sum_{j=1}^m{\left(j(r\lambda+2m-j-1)q_{m,j}+q_{m,j-1}\right)\tr^{j-1}(z)\calBe^jp} = 0.
\end{align*}
\item[\textup{(2)}] We now define operators $Q_k$ by applying $Q_0$ to $\calBe^kp$:
\begin{align*}
 Q_kp &:= Q_0(\calBe^kp) = \sum_{j=0}^{m-k}{q_{m-k,j}\tr(z)^j\calBe^{k+j}p}.
\end{align*}
Multiplication with $\tr(z)^k$ yields
\begin{align*}
 \tr(z)^kQ_kp &= \sum_{j=0}^{m-k}{q_{m-k,j}\tr(z)^{k+j}\calBe^{k+j}p}\\
 &= \sum_{j=k}^m{q_{m-k,j-k}\tr(z)^j\calBe^jp}.
\end{align*}
If we denote $a_{k,l}:=q_{m-k,j-k}$ for $k=0,\ldots,m$ and $j=k,\ldots,m$, then
\begin{align*}
 \left(\begin{array}{cccc}a_{0,0} & a_{0,1} & \cdots & a_{0,m}\\0 & a_{1,0} & \cdots & a_{1,m}\\\vdots & \ddots & \ddots & \vdots\\0 & \cdots & 0 & a_{m,m}\end{array}\right)\left(\begin{array}{c}p\\\tr(z)\calBe p\\\vdots\\\tr(z)^m\calBe^mp\end{array}\right) &= \left(\begin{array}{c}Q_0p\\\tr(z)Q_1p\\\vdots\\\tr(z)^mQ_mp\end{array}\right).
\end{align*}
Since $a_{k,k}=1$, the matrix on the left hand side is invertible and in particular there exist constants $b_0,\ldots,b_m$ which are independent of $p$ such that
\begin{align}
 p &= \sum_{k=0}^m{b_k\tr(z)^kQ_kp}.\label{eq:ExpansionIntoHarmonicsInTermsOfQ}
\end{align}
Since $Q_kp\in\calH^{m-k}(\XX)$ this shows the existence of the claimed decomposition with $h_{m-k}=b_kQ_kp$.
\item[\textup{(3)}] We now prove the explicit formula for $h_{m-k}$. For this we have to find the constants $b_0,\ldots,b_m$. If we substitute $\tr(z)^kQ_kp$ for $p$ in \eqref{eq:ExpansionIntoHarmonicsInTermsOfQ} then we obtain
\begin{align*}
 \tr(z)^kQ_kp &= \sum_{j=0}^m{b_j\tr(z)^jQ_j(\tr(z)^kQ_kp)}.
\end{align*}
Since we have already proved uniqueness, we find
\begin{align*}
 \tr(z)^kQ_kp &= b_k\tr(z)^kQ_k(\tr(z)^kQ_kp)\\
\intertext{and hence}
 Q_kp &= b_kQ_k(\tr(z)^kQ_kp)\\
 &= b_k\sum_{j=0}^{m-k}{q_{m-k,j}\tr(z)^j\calBe^{k+j}(\tr(z)^kQ_kp)}.\\
\intertext{Applying Lemma \ref{lem:BesselOnProductWithTrace} again gives}
 &= b_kq_{m-k,0}\cdot\frac{k!\,\Gamma(r\lambda+2m-k)}{\Gamma(r\lambda+2m-2k)}Q_kp.
\end{align*}
There are clearly polynomials $p\in\calP^m(\XX)$ with $Q_kp\neq0$, e.g. for $\tr(z)^kq$, $q\in\calH^m(\XX)$ it follows from the uniqueness that $q=b_kQ_k(\tr(z)^kq)$. Hence,
\begin{align*}
 b_k &= \frac{\Gamma(r\lambda+2m-2k)}{k!\,\Gamma(r\lambda+2m-k)}
\end{align*}
and the claimed formula for $h_{m-k}=b_kQ_kp$ follows.\qedhere
\end{enumerate}
\end{proof}

Recall from Remark \ref{rem:KHSonXX} that $\calP^m(\XX)\simeq\calP^{(m,0,\ldots,0)}(V_\CC)$. Denote by $d_{\bf m}$ the dimension of $\calP^{\bf m}(V_\CC)$. Then, using the results of \cite[Section XIV.5]{FK94}, we obtain
\begin{align*}
 d_m &:= d_{(m,0,\ldots,0)} = \frac{(\frac{n}{r})_m(r\lambda)_m}{m!\,(\lambda)_m}.
\end{align*}
For convenience we also put $d_{-1}:=0$. The following dimension formula for $\calH^m(\XX)$ is now immediate with Proposition \ref{prop:PolyDecompIntoHarmonics}:

\begin{corollary}
$$ \dim\,\calH^m(\XX) = d_m-d_{m-1}. $$
\end{corollary}

\begin{example}
For $V=\Sym(k,\RR)$ we have $r=k$, $d=1$ and $n=\frac{k}{2}(k+1)$. Hence, $\lambda=\frac{1}{2}$ and
\begin{align*}
 d_m &= \frac{(\frac{k}{2})_m(\frac{k+1}{2})_m}{m!\,(\frac{1}{2})_m} = \frac{(n)_m}{(2m)!} = {n+2m-1\choose2m} = {n+2m-1\choose n-1}
\end{align*}
which is exactly the dimension of the space of homogeneous polynomials of degree $2m$ in $n$ variables. For the dimension of $\calH^m(\XX)$ we obtain
\begin{align*}
 \dim\,\calH^m(\XX) &=  {n+2m-1\choose n-1}-{n+2m-3\choose n-1}
\end{align*}
which is the well-known formula for the dimension of $\calH^{2m}(\RR^n)$.
\end{example}

For fixed $p\in\calH^m(\XX)$ the span of the polynomials $\tr(z)^kp$, $k\in\NN$, is invariant under the action of $\sl(2,\RR)$
 by Lemma \ref{lem:ActionOnPowersOfTrace} and Lemma \ref{lem:BesselOnProductWithTrace} and defines an irreducible representation of $\sl(2,\RR)$ of lowest weight $r\lambda+2m$. We denote the corresponding representation on $\calW_{r\lambda+2m}:=\linspan\{\tr(z)^k:k\in\NN\}$ by $\rho_{r\lambda+2m}$. 
With Lemma \ref{lem:ActionOnPowersOfTrace} and Lemma \ref{lem:BesselOnProductWithTrace} we find
\begin{align*}
 \td\rho_{r\lambda+2m}(\widetilde{E})\tr^k(z) &= \sqrt{-1}\tr^{k+1}(z),\\
 \td\rho_{r\lambda+2m}(\widetilde{H})\tr^k(z) &= (r\lambda+2m+2k)\tr^k(z),\\
 \td\rho_{r\lambda+2m}(\widetilde{F})\tr^k(z) &= k(r\lambda+2m+k-1)\sqrt{-1}\tr^{k-1}(z).
\end{align*}
The lowest weight vector is given by $\tr^0(z)=\1$. Putting things together gives:

\begin{theorem}\label{thm:DualPairDecompositionFock}
Under the action of $(\fraks,\frakk^\frakl)$ the representation $(\td\rho,\calP(\XX))$ decomposes as
\begin{align*}
 \calP(\XX) &\simeq \bigoplus_{m=0}^\infty{\calW_{r\lambda+2m}\boxtimes\calH^m(\XX)},
\end{align*}
where $\calW_{r\lambda+2m}$ denotes the irreducible representation of $\fraks\simeq\sl(2,\RR)$ of lowest weight $r\lambda+2m$ and $\calH^m(\XX)$ is an irreducible representation of $Z_{G^\vee}(\fraks)$.
\end{theorem}

\begin{remark}
The corresponding decomposition
 in the Schr\"odinger model
 was given in Theorem \ref{thm:sl2decoS}.
\end{remark}

Using the theory of spherical harmonics we can now prove that the completion $\widetilde{\calF}(\XX)$ of the space of polynomials and the intrinsically defined Fock space $\calF(\XX)$ (see \eqref{eq:DefFockSpace}) agree.

\begin{theorem}\label{thm:NaturalFockSpace}
We have $\widetilde{\calF}(\XX)=\calF(\XX)$.
\end{theorem}

\begin{proof}
We first treat the case $V=\RR$. Here we have $\XX=\CC^\times$ with norm given by
\begin{align*}
 \|F\|^2 &= \const\cdot\int_\CC{|F(z)|^2\omega(z)|z|^{2(\lambda-1)}\td z}.
\end{align*}
Using the asymptotic behavior of the K-Bessel function near $z=0$ (see Appendix \ref{app:BesselFcts}) we find that in the Laurent expansion near $z=0$ of function $F\in\calF(\XX)$ all negative terms have to vanish and hence $F\in\calO(\CC)$. By Proposition \ref{prop:FockSpaceGlobalFcts} such a function already belongs to $\widetilde{\calF}(\XX)$ and the proof is complete.\\
Now suppose that $r>1$. Put $K'=Z(K^\vee)\times K^L$, where $K^\vee\subseteq G^\vee$ denotes the maximal compact subgroup of $G^\vee$ corresponding to $\frakk$. We note that $Z(K^\vee)\simeq\SO(2)$ is given by $\{\exp(t(E-F)):t\in\RR\}$ and that $Z(K^\vee)$ is the intersection of $K^\vee$ and the analytic subgroup of $G^\vee$ with Lie algebra $\fraks\simeq\sl(2,\RR)$. Then by Theorem \ref{thm:DualPairDecompositionFock} we have
\begin{align*}
 \widetilde{\calF}(\XX)_{K'} &= \widetilde{\calF}(\XX)_{K^\vee} = \calP(\XX).
\end{align*}
Note that by Proposition \ref{prop:PolyDecompIntoHarmonics} we further know that
\begin{align*}
 \calP(\XX) &= \bigoplus_{k,m=0}^\infty{\tr^k(z)\calH^m(\XX)}.
\end{align*}
Now, the representation $\rho|_{K'}$ on $\widetilde{\calF}(\XX)$ (see Theorem \ref{thm:UnitaryRepOnFock}) extends to a unitary $K'$-representation on $\calF(\XX)$ given by
\begin{align*}
 (\exp(t(E-F)),k)\cdot F(z) &= e^{r\lambda t\sqrt{-1}}F(e^{2t\sqrt{-1}}k^{-1}z), & z\in\XX.
\end{align*}
We now calculate the $K'$-finite vectors in $\calF(\XX)$ by finding the $K'$-finite vectors in $\calO(\XX)$ and then intersecting with $L^2(\XX,\omega\td\nu)$. Note that the polar coordinates map
\begin{align*}
 q:\RR_+\times\SS\stackrel{\sim}{\to}\Xi\subseteq V,\,(r,x)\mapsto rx,
\end{align*}
extends to a holomorphic embedding
\begin{align*}
 q_\CC:\CC^\times\times\SS_\CC\to\XX\subseteq V_\CC,\,(s,z)\mapsto sz,
\end{align*}
onto the open subset $\{z\in\XX:\tr(z)\neq0\}\subseteq\XX$, where $\SS_\CC=K^L_\CC/(K^L_\CC)_{c_1}\subseteq V_\CC$, $K^L_\CC\subseteq G_\CC$ denoting the complexification of $K^L$ in $G_\CC$. Then $q_\CC$ induces a restriction map
\begin{align*}
 q_\CC^*:\calO(\XX)\to\calO(\CC^\times\times\SS_\CC).
\end{align*}
Under this restriction the Lie algebra $\frakk^\frakl_\CC$ acts on $\calO(\SS_\CC)$ via vector fields on $\SS_\CC$ and $\so(2,\CC)$ acts by the Euler operator on $\calO(\CC^\times)$. 
Therefore we obtain the decomposition
\begin{align*}
 \calO(\CC^\times\times\SS_\CC)_{K'} &= \bigoplus_{k\in\ZZ}{\bigoplus_{m=0}^\infty{\tr(z)^k\calH^m(\XX)}}.
\end{align*}
We claim that $\tr(z)^{k}$ only extends to a holomorphic function on $\XX$ if $k\geq0$. For this we show that there exists $z\in\XX$ with $\tr(z)=0$. (Note that the following argument only works for $r>1$.) Let $x_0\in V_{12}$ be any element with $|x_0|^2=2$. Then $x_0^2=c_1+c_2$. We claim that $z=c_1+\sqrt{-1}x_0-c_2\in\XX$, but clearly $\tr(z)=0$. To show that $z\in\XX$ we prove that $z=\exp(2\sqrt{-1}c_2\Box x_0)c_1\in L_\CC\cdot c_1=\XX$. In fact, we find
\begin{align*}
 (c_2\Box x_0)c_1 &= \frac{1}{2}x_0,\\
 (c_2\Box x_0)x_0 &= c_2,\\
 (c_2\Box x_0)c_2 &= 0
\end{align*}
and the claim follows. Hence, we obtain
\begin{align*}
 \calF(\XX)_{K'} &= \bigoplus_{k,m=0}^\infty{\tr(z)^k\calH^m(\XX)} = \widetilde{\calF}(\XX)_{K'}
\end{align*}
and therefore its completions $\calF(\XX)$ and $\widetilde{\calF}(\XX)$ have to agree.
\end{proof}

From now on we only use the notation $\calF(\XX)$
 for the Fock space.
In Section \ref{sec:sbtransform},
 we shall give another equivalent definition
 (an \lq{extrinsic definition}\rq\
 with respect to the embedding 
 $\XX \hookrightarrow V_{\CC}$). 
\section{The Segal--Bargmann transform}
\label{sec:sbtransform}

Generalizing the classical Segal--Bargmann transform,
 we explicitly construct an intertwining operator
 $\BB_{\Xi}$ 
between the Schr\"odinger model (Section \ref{sec:schmodel})
 and the Fock model (Section \ref{sec:fockspace})
 of the \minholrep\
 of $G^\vee$ in terms of its integral kernel. 
As an application,
 we establish the equivalence
 of three different definitions
 for the Fock space $\calF(\XX)$,
 including an \lq{intrinsic one}\rq\
 and an \lq{extrinsic one}\rq\
 (Corollary \ref{cor:Fock3}).  
Further
 the Segal--Bargmann transform $\BB_{\Xi}$
 brings us naturally
 to a generalization of the classical Hermite polynomials
 as the preimages of the monomials
 in the Fock model $\calF(\XX)$.

\subsection{Definition and properties}

Let $\widetilde{I}_{\alpha}(t):=(\frac {t}{2})^{-\alpha}I_{\alpha}(t)$
 be the renormalized I-Bessel function.  
Then $\widetilde{I}_{\alpha}(t)$ is an entire function
 on $\CC$
 (see Appendix \ref{app:BesselFcts}).  
We set an entire function $B$ on $\CC$ by
\begin{equation}
\label{eqn:Bker}
 B(t) := \Gamma(\lambda)\widetilde{I}_{\lambda-1}(2\sqrt{t}).  
\end{equation}
Clearly
 we have 
\[
     B(0)=1. 
\]
For $x,z \in V_{\CC}$, 
 we write simply $B(x|z)=B((x|z))$.
We are ready to define an integral transform
 $\BB_\Xi:C_c(\Xi)\to \calO(V_\CC)$ 
 by 
\begin{align}
\label{eqn:BO}
 \BB_\Xi\psi(z) &:= e^{-\frac{1}{2}\tr(z)}\int_\Xi{B(x|z)e^{-\tr(x)}\psi(x)\td\mu(x)} & z\in\XX,
\end{align}
for a compactly supported continuous function $\psi \in C_c(\Xi)$.  
We call $\BB_\Xi$ is the Segal--Bargmann transform.  

\begin{lemma}\label{lem:SBok}
The Segal--Bargmann transform $\BB_\Xi$ 
 is well-defined as a linear map
$
 L^2(\Xi)\to\calO(V_\CC).
$
\end{lemma}

\begin{proof}
Since the kernel function $e^{-\frac{1}{2}\tr(z)}B(x|z)e^{-\tr(x)}$
 is obviously analytic in $z$, 
it suffices to show that its $L^2$-norm in $x$ has a uniform bound
 on $\|z\|\leq R$
 for an arbitrary fixed $R>0$. 
Using the asymptotic behaviour
 of the I-Bessel function $\widetilde{I}_\alpha(t)$ as $t\to\infty$
 (see Appendix \ref{app:BesselFcts}) we obtain
\begin{align*}
 |B(t)| &\lesssim |t|^{\frac{1-2\lambda}{4}}e^{2|t|^{\frac{1}{2}}} \lesssim |t|^{\max(0,\frac{1-2\lambda}{4})}e^{2|t|^{\frac{1}{2}}}.
\end{align*}
Since $B(t)$ is analytic up to $t=0$,
 we have a uniform bound on $\CC$:
\begin{align*}
 |B(t)| &\lesssim (1+|t|^{\max(0,\frac{1-2\lambda}{4})})e^{2|t|^{\frac{1}{2}}}.
\end{align*}
Then for $x\in\Xi$, $z\in V_\CC$ with $\|z\|\leq R$, we find
\begin{align*}
 |e^{-\frac{1}{2}\tr(z)}B(x|z)e^{-\tr(x)}| &\lesssim (1+|(x|z)|^{\max(0,\frac{1-2\lambda}{4})})e^{2|z|^{\frac{1}{2}}|x|^{\frac{1}{2}}-|x|}\\
 &\leq (1+|Rx|^{\max(0,\frac{1-2\lambda}{4})})e^{2R^{\frac{1}{2}}|x|^{\frac{1}{2}}-|x|}
\end{align*}
which is $L^2$ in $x$ with norm independent
 of $z$ and the claim follows.
\end{proof}

Next,
 we show
 that $\BB_\Xi$ intertwines the action $\td\pi$ of $\frakg$
 on the space $L^2(\Xi)^\infty$ 
 of smooth vectors with the action $\td\rho$ on $\calF(\XX)^\infty$.

\begin{theorem}
\label{thm:3.2}
For any $X \in \frakg$,
\[
\BB_{\Xi} \circ \td \pi(X) = \td \rho(X) \circ \BB_\Xi
\quad
\text{ on }L^2(\Xi)^{\infty}.  
\]
\end{theorem}

\begin{proof}
Since $\td \rho = \td \pi_{\CC} \circ c$
 by definition, 
Theorem \ref{thm:3.2} can be restated
 in terms of the Jordan algebra
 as follows:
\begin{align*}
 \BB_\Xi\circ\td\pi(a,0,0) &= \textstyle\td\pi_\CC(\frac{a}{4},\sqrt{-1}L(a),a)\circ\BB_\Xi,\\
 \BB_\Xi\circ\td\pi(0,L(a)+D,0) &= \textstyle\td\pi_\CC(\sqrt{-1}\frac{a}{4},D,-\sqrt{-1}a)\circ\BB_\Xi,\\
 \BB_\Xi\circ\td\pi(0,0,a) &= \textstyle\td\pi_\CC(\frac{a}{4},-\sqrt{-1}L(a),a)\circ\BB_\Xi, 
\end{align*}
 for any $a \in V$ and $D \in {\frak{aut}}(V)$.  
The verification
 of these formulae is an easy,
 though lengthy calculation. 
We outline the method
 and the crucial steps. 
First note that the operators $\td\pi(X)$, $X\in\frakg$, are skew-symmetric on $L^2(\Xi,\td\mu)$. Therefore we have
\begin{align*}
 (\BB_\Xi\circ\td\pi(X)f)(z) &= e^{-\frac{1}{2}\tr(z)}\int_\Xi{B(x|z)e^{-\tr(x)}(\td\pi(X)f)(x)\td\mu(x)}\\
 &= -e^{-\frac{1}{2}\tr(z)}\int_\Xi{\td\pi(X)\left[B(x|z)e^{-\tr(x)}\right]f(x)\td\mu(x)}.
\end{align*}
Using the explicit formulas for the action $\td\pi(X)$ and especially for the action of $X\in\overline{\frakn}$ the product rule for the Bessel operator (see Lemma \ref{lem:BesselProdRule}), the action of $\td\pi(X)$ on $B(x|z)e^{-\tr(x)}$ can be computed. Further, we know that $\calB_xB(x|z)=zB(x|z)$ and $\calB_zB(x|z)=xB(x|z)$. The rest is standard.
\end{proof}

Next we examine how $\BB_{\Xi}$ acts 
 on the lowest weight $\frakk$-type.  
By a little abuse
 of notation
 that will be justified soon,
 we set 
\begin{equation}
\label{eqn:Psi0}
     \Psi_0:=\BB_\Xi\psi_0, 
\end{equation}
where 
we recall from \eqref{eqn:psi0}
 that $\psi_0(x)=e^{-\operatorname{tr}(x)}$.  

\begin{lemma}
\label{lem:3.3}
$\Psi_0(0) = 1.$
\end{lemma}

\begin{proof}
{}From \eqref{eqn:normpsi} we have
\begin{align*}
 \Psi_0(0) &= B(0)\int_\Xi{e^{-2\,\tr(x)}\td\mu(x)}.
\qedhere
\end{align*}
\end{proof}

\begin{proposition}
\label{prop:3.4}
\begin{enumerate}
\item[\textup{(1)}]
$\BB_\Xi\psi_0 = \1.$
\item[\textup{(2)}]
$\BB_{\Xi}$ induces a $(\frakg,\frakk)$-isomorphism
 between $L^2(\Xi)_{\frakk}$ and $\calP(\XX)$.
\end{enumerate}
\end{proposition}

\begin{proof}
\begin{enumerate}
\item[\textup{(1)}]
First we regard $\Psi_0$ as a holomorphic function 
 on $V_\CC$.  
Since the action on the Fock model
 is given by $\td \rho=\td \pi_\CC \circ c$
 and the Cayley transform $c$
 sends $\frakk$ to a totally real subspace
 of $\frakl_\CC$ by
\begin{align*}
 \frakk\to\frakl_\CC,(a,D,-a)\mapsto D+2\sqrt{-1}L(a).  
\end{align*}
Thorem \ref{thm:3.2} implies
 that $\Psi_0$ is $\frakl_\CC$-invariant. 
Hence,
 it has to be constant on every $L_\CC$-orbit. 
Since $\Psi_0$ is holomorphic on $V_\CC$
 and $V_\CC$ decomposes into finitely many $L_\CC$-orbits,
 it follows that $\Psi_0$ is constant on $V_\CC$. Hence the first statement follows from Lemma \ref{lem:3.3}.
\item[\textup{(2)}]
We know 
 that the underlying $(\frakg,\frakk)$-module
 $L^2(\Xi)_{\frakk}$ of the Schr{\"o}dinger model
 $(\pi,L^2(\Xi))$
 is irreducible.  
Further,
 $\calP(\XX)$ is an irreducible
 $(\frakg,\frakk)$-module
 by Proposition \ref{prop:irredPX}. 
Since $\BB_\Xi$ is non-zero
 and $\BB_\Xi$ intertwines the actions $\td\pi$ and $\td\rho$,
 $\BB_{\Xi}$ gives an isomorphism
 of $(\frakg,\frakk)$-modules.  
\end{enumerate}
\end{proof}

\begin{theorem}
\label{thm:SBunitary}
The Segal--Bargmann transform 
$\BB_\Xi$ is a unitary isomorphism $L^2(\Xi)\to\calF(\XX)$.
\end{theorem}

\begin{proof}
It only remains to show
 that $\BB_\Xi$ is isometric
 between $L^2(\Xi)_{\frakk}$ and $\calP(\XX)$,
 because $L^2(\Xi)_{\frakk} \subseteq L^2(\Xi)$
 and $\calP(\XX)\subseteq\calF(\XX)$ are dense.\\

Since both $L^2(\Xi)_{\frakk}$ and $\calP(\XX)$
 are irreducible,
 infinitesimally unitary $(\frakg,\frakk)$-modules,
 $\BB_{\Xi}$ is a scalar multiple
 of a unitary operator.  
Since
\[
\langle {\bf {1}}, {\bf {1}} \rangle_{L^2(\XX,\omega \td \nu)}
=1
=\langle \psi_0,\psi_0\rangle_{L^2(\Xi,\td \mu)}, 
\]
$\BB_{\Xi}$ must be
 a unitary operator by Proposition \ref{prop:3.4}.  
\end{proof}

\begin{corollary}
The inverse Segal--Bargmann transform is given by
\begin{align*}
 \BB_\Xi^{-1}F(x) &= e^{-\tr(x)}\int_\XX{B(x|\overline{z})e^{-\frac{1}{2}\tr(\overline{z})}F(z)\omega(z)\td\nu(z)}
\end{align*}
\end{corollary}

\begin{proof}
\begin{align*}
 \langle\BB_\Xi^{-1}F,\psi\rangle &= \langle F,\BB_\Xi\psi\rangle\\
 &= \int_\XX{F(z)\overline{\BB_\Xi\psi(z)}\omega(z)\td\nu(z)}\\
 &= \int_\XX{\int_\Xi{F(z)e^{-\frac{1}{2}\overline{\tr(z)}}\overline{B(x|z)}e^{-\tr(x)}\overline{\psi(x)}\td\mu(x)}\omega(z)\td\nu(z)}\\
 &= \int_\Xi{\int_\XX{e^{-\frac{1}{2}\tr(\overline{z})}B(x|\overline{z})e^{-\tr(x)}F(z)\omega(z)\td\nu(z)}\overline{\psi(x)}\td\mu(x)}\qedhere
\end{align*}
\end{proof}

We can now use the Segal--Bargmann transform to obtain a different description of the Fock space.

\begin{theorem}\label{thm:FockAsRestriction}
$$ \calF(\XX) = \left\{F|_\XX:F\in\calO(V_\CC),\int_\XX{|F(z)|^2\omega(z)\td\nu(z)}<\infty\right\}. $$
\end{theorem}

\begin{proof}
The inclusion $\supseteq$ holds by Proposition \ref{prop:FockSpaceGlobalFcts}. The other inclusion now follows
 with Lemma \ref{lem:SBok}
 since $\BB_\Xi:L^2(\Xi)\to\calF(\XX)$ is an isomorphism.
\end{proof}

\begin{remark}
We note that the restriction map
 $\calO(V_{\CC}) \to \calO(\XX)$
 is not surjective, 
 and therefore the above equivalence is non-trivial.  
\end{remark}
Combining Theorem \ref{thm:NaturalFockSpace}
 with Theorem \ref{thm:FockAsRestriction},
 we have obtained three equivalent definitions
 of the Fock space $\calF(\XX)$ as follows:
\begin{corollary}
\label{cor:Fock3}
The following three subspaces
 are the same.  
\begin{itemize}
\item[$\bullet$]
$\calO(\XX) \cap L^2(\XX, \omega \td \nu)$.  
\item[$\bullet$]
The completion of $\calP(\XX)$ in $L^2(\XX, \omega \td \nu)$.  
\item[$\bullet$]
$\{F|_{\XX}:F \in \calO(V_{\CC})\} \cap L^2(\XX, \omega \td \nu)$.  
\end{itemize}
\end{corollary}

\subsection
{Relations with the classical Segal--Bargmann
 transform}
\label{subsec:relclassic}

In the cases $V=\Herm(k,\FF)$
 with $\FF=\RR,\CC,\HH$,
 we can relate our Segal--Bargmann transform $\BB_\Xi$
 directly with the classical Segal--Bargmann transform
by using the folding map. 

As in Subsection \ref{subsec:foldingmaps}
 we let $d=\dim_\RR\FF=2\lambda$
 and define the complexification of the folding map
 \eqref{eqn:fold} by  
\[
     p_\CC:\FF_\CC^k:=\FF^k\otimes_\RR\CC\to\XX,\,z\mapsto zz^*.  
\]
The $p_\CC$ is a principal bundle
 with structure group $\ZZ_2$, 
 $\CC^{\times}$, 
 and $SL(2,\CC)$, 
 as $p:\FF^k\setminus\{0\} \to \Xi$
 is the one with $U(1;\FF)$, 
 $\FF=\RR$, $\CC$, and $\HH$, respectively.  
Note that the conjugation for 
 $z^*=\overline{{}^{t\!}z}$
 is taken as the conjugation in $\FF$, 
 and not the one corresponding to the complexification $\FF_\CC=\FF\otimes_\RR\CC$. 

Let $\calF(\CC^n)$ denote the classical Fock space
 on $\CC^n$ 
 with respect to the Gaussian measure
 $e^{-|z|^2}\td z$, 
and $\BB:L^2(\RR^n)\to\calF(\CC^n)$ the classical Segal--Bargmann transform given by
\begin{align*}
 \BB u(z) &= e^{-\frac{1}{2}z^2}\int_{\RR^n}{e^{2z\cdot x}e^{-x^2}u(x)\td x}.
\end{align*}
We consider the following diagram:
\[\xymatrix{
 L^2(\Xi,\td\mu) \ar^/-5em/{p^*}[rrrr] \ar_{\BB_\Xi}[d] &&&& **[l] L^2(\FF^k)^{U(1;\FF)} \subseteq L^2(\FF^k) \simeq L^2(\RR^{dk})\!\!\!\!\!\! \ar^{\BB}[d]\\
 \calF(\XX) \ar^/-5.5em/{p_\CC^*}[rrrr] &&&& **[l] \calF(\FF_\CC^k)^{U(1;\FF)} \subseteq \calF(\FF_\CC^k) \simeq \calF(\CC^{dk}),\!\!\!\!\!\!
}\]

\begin{theorem}
\label{thm:foldB}
$p_{\CC}^{\ast} \circ \BB_{\Xi}$
 is a scalar multiple of $\BB \circ p^{\ast}$.  
\end{theorem}

We shall give a proof 
 of this theorem 
 by comparing the integral kernels
 of $\BB_{\Xi}$ and $\BB$.  
Conversely,
 we may use the above diagram
 for the definition of the Segal--Bargmann transform
 for the minimal representations
 arising from the Euclidean Jordan algebra
 $V=\operatorname{Herm}(k,\FF)$,
 $\FF=\RR$, 
 $\CC$ or $\HH$.  
 
\begin{proof}
The integral kernel for $\BB \circ p^{\ast}$ is 
 obtained by integrating the kernel
 of the classical Segal--Bargmann transform over $U(1;\FF)$,
 i.e. over its orbit $S^{d-1}$, a $(d-1)$-dimensional sphere, using the integral formula
\begin{equation*}
 \int_{S^{d-1}}{e^{r\omega\cdot\eta}\td\omega} = 2\pi^{\frac{d}{2}}\widetilde{I}_{\frac{d}{2}-1}(r),
\end{equation*}
it amounts to a scalar multiple
 of the integral kernel of $\BB_{\Xi}$,
 which is by \eqref{eqn:BO}
\begin{align*}
 \Gamma(\lambda)e^{-\frac{1}{2}\tr(z)}\widetilde{I}_{\lambda-1}(2\sqrt{(z|x)})e^{-\tr(x)}.  
\end{align*}
Thus Theorem \ref{thm:foldB} is proved.  
\end{proof}

\begin{example}
\label{ex:foldB}
In the case $\FF =\RR$,
 we have $\frakg=\sp(k,\RR)$.  
Then $p^{\ast}$ induces an isomorphism 
 $L^2(\Xi, \td \mu)
\overset \sim \to L_{\operatorname{even}}^2(\RR^k)$, 
 the even part of the metaplectic representation
 on even $L^2$-functions
 on $\RR^k$, 
 and $p_{\CC}^{\ast}$ induces 
 an isomorphism $\calF(\XX)
 \overset \sim \to \calF_{\operatorname{even}}(\CC^k)$
 the even part of the Fock space
 on $\CC^k$.  
The kernel function
 for the classical Segal--Bargmann transform
 $\BB:L^2(\RR^k)\to\calF(\CC^k)$
 and that for our Segal--Bargmann transform
 is related by the integration
 of $S^0$
 (two points),
 namely,
\begin{equation*}
 e^r+e^{-r} = (2\pi^{\frac{1}{2}})(\tfrac{1}{\sqrt{\pi}}\cosh r).  
\end{equation*}
\end{example}

\subsection{Generalized Hermite functions}

Now let $B:=(e_j)_j\subseteq V$ be any basis of $V$. 
For a multiindex $\alpha\in\NN^B$ we use the notation
\begin{align*}
 z^\alpha &:= \prod_j{(e_j|z)^{\alpha_j}},\\
 \calB^\alpha &:= \prod_j{(e_j|\calB)^{\alpha_j}}.
\end{align*}

\begin{remark}
The monomials $z^\alpha$ do not form an orthogonal system
 in $\calP(\XX)$.
In fact,
 the monomials $z^\alpha$ are not linearly independent
 in $\calP(\XX)$
 since there are polynomials
 that vanish on $\XX$. 
The space $\calP(\XX)$ of regular functions
 on $\XX$ is defined to be the quotient of the space $\CC[z]$
 of (holomorphic) polynomials on $V_\CC$
 by the ideal generated by all polynomials vanishing on $\XX$. 
\end{remark}

The \textit{generalized Hermite functions} on $\Xi$ are defined by
\begin{align*}
 h_\alpha(x) &:= e^{\tr(x)}\calB^\alpha e^{-2\,\tr(x)}, & x\in\Xi,
\end{align*}
for $\alpha\in\NN^B$. In particular, for $\alpha=0$ we have $h_0=\psi_0$. Note that
\begin{align*}
 h_\alpha(x) &= H_\alpha(x)e^{-\tr(x)},
\end{align*}
where $H_\alpha(x)$ is a polynomial of degree $|\alpha|$. 
We call $H_\alpha(x)$ the \textit{generalized Hermite polynomial}.

Now we can show that the Segal--Bargmann transform $\BB_\Xi$ maps the Hermite functions $h_\alpha$ onto the monomials $z^\alpha$.

\begin{proposition}
$\BB_\Xi h_\alpha = z^\alpha.$
\end{proposition}

\begin{proof}
Since the Bessel operator $\calB$ is symmetric with respect to the inner product on $L^2(\Xi)$, we obtain
\begin{align*}
 \BB_\Xi h_\alpha(z) &= e^{-\frac{1}{2}\tr(z)}\int_\Xi{B(x|z)e^{-\tr(x)}h_\alpha(x)\td\mu(x)}\\
 &= e^{-\frac{1}{2}\tr(z)}\int_\Xi{B(x|z)\cdot\calB^\alpha e^{-2\,\tr(x)}\td\mu(x)}\\
 &= e^{-\frac{1}{2}\tr(z)}\int_\Xi{\calB_x^\alpha B(x|z)\cdot e^{-2\,\tr(x)}\td\mu(x)}\\
 &= e^{-\frac{1}{2}\tr(z)}\int_\Xi{z^\alpha B(x|z)e^{-2\,\tr(x)}\td\mu(x)}\\
 &= z^\alpha\Psi_0(z).\qedhere
\end{align*}
\end{proof}

\begin{corollary}
Each $\frakk$-type $W_m\subseteq L^2(\Xi,\td\mu)$ is spanned by the functions $h_\alpha$ for $|\alpha|=m$.
\end{corollary}

\begin{proof}
Since in the Fock model each $\frakk$-type $\calP^m(\XX)$ is spanned by the monomials $z^\alpha$ for $|\alpha|=m$ this is clear by the previous theorem.
\end{proof}

\begin{remark}
Suppose, the basis $B=(e_j)_j$ is chosen such that $e_1={\bf e}$ is the identity of the Jordan algebra. Then for $\alpha=(m,0\ldots,0)$ we obtain
\begin{align*}
 h_\alpha(x) &= e^{\tr(x)}\calBe^m e^{-2\,\tr(x)}.
\end{align*}
Further,
\begin{align*}
 z^\alpha &= \tr(z)^m
\end{align*}
and hence, $z^\alpha$ is the unique
 (up to scalar) $\frakk^\frakl$-invariant vector
 in the $\frakk$-type $\calP^m(\XX)$. 
Since $\BB_\Xi$ is an intertwining operator and $\BB_\Xi h_\alpha=z^\alpha$, we obtain that $h_\alpha$ is the unique (up to scalar) $\frakk^\frakl$-invariant vector in the $\frakk$-type $W_m$. By \cite{HKM11} we know that also the Laguerre function
\begin{align*}
 \ell_m^\lambda(x) &= e^{-\tr(x)}L_m^\lambda(2\,\tr(x))
\end{align*}
is a $\frakk^\frakl$-invariant vector in the $\frakk$-type $W_m$
 and hence, $h_\alpha$ and $\ell_m^\lambda$ have to be proportional to each other, which means that
\begin{align*}
 e^{\tr(x)}\calBe^m e^{-2\,\tr(x)} = \const\cdot e^{-\tr(x)}L_m^{\lambda}(2\,\tr(x)).
\end{align*}
\end{remark}
\section{The unitary inversion operator}
\label{sec:UnitaryInversionOperator}

The Schr{\"o}dinger model
 of the minimal representation $\pi$
 on $L^2(\Xi)$
 has an advantage
 that the representation space is simple, 
 namely, 
 it is the Hilbert space
 consisting of arbitrary $L^2$-functions on $\Xi$.  
Another advantage is 
 that the group action 
 of a maximal parabolic subgroup is also simple.  
Thus the {\it{unitary inversion operator}}
 $\calF_\Xi$ 
 (see \eqref{eqn:FO} for the definition below)
 plays a key role
 in the global action of $G^{\vee}$ on $L^2(\Xi)$.  
See \cite[Chapter 1]{KMa11}
 for a comparison of different models
 of minimal representations.  
The operator $\calF_\Xi$ is 
 essentially the Euclidean Fourier transform 
 on the metaplectic representation $L^2({\mathbb{R}}^n)$
 for ${\mathfrak {g}}={\mathfrak {sp}}(n,{\mathbb{R}})$.  
The program to find the integral kernel
 of the unitary inversion operator $\calF_\Xi$
 has been carried out
 in \cite{KM07a} for ${\mathfrak {g}}={\mathfrak {so}}(2,n)$
 and in \cite{KMa11} also
 for a non-Hermitian Lie algebra
 ${\mathfrak {g}}={\mathfrak {so}}(p,q)$
 ($p+q$ even)
 in terms of the Bessel function
 (or the \lq{Bessel distribution}\rq).  
In this section 
 we take another approach
 to find an explicit integral kernel
 of $\calF_\Xi$ 
 as an application 
 of the results in Section \ref{sec:sbtransform}
 on the Segal--Bargmann transform
 under the assumption
 that $G$ is a simple Hermitian Lie group of tube type. 

In the framework of Jordan algebras,
 the unitary inversion operator is the action
 $\pi(\widetilde{j})$
 of the inversion element
 $\widetilde{j}=\exp_{\widetilde{G}}(\frac{\pi}{2}({\bf e},0,-{\bf e}))\in\widetilde{G}$ 
 up to a phase factor (see \cite[Section 3.3]{HKM11}). 
More precisely, we set
\begin{align}
\label{eqn:FO}
 \calF_\Xi &:= e^{-\pi\sqrt{-1}\frac{r\lambda}{2}}\pi(\widetilde{j}).
\end{align}
The operator $\calF_\Xi$ is unitary on $L^2(\Xi,\td\mu)$ of order $2$, i.e. $\calF_\Xi^2=\id$.

Let $\widetilde {J}_{\alpha}(z):=(\frac z 2)^{-\alpha}J_{\alpha}(z)$
 be the renormalized J-Bessel function,
 which is an entire function on $\CC$
 (see Appendix \ref{app:BesselFcts}).  
We define an entire function $F$ on $\CC$ by
\begin{align}
\label{eqn:Fkernel}
 F(z) &= 2^{-r\lambda}B(-z) = 2^{-r\lambda}\Gamma(\lambda)\widetilde{J}_{\lambda-1}(2\sqrt{z})
\end{align}
and write $F(x|y)=F((x|y))$, $x,y\in\Xi$,
 for short.  

Denote by $L^2(\Xi)_{\frakk}$
 the space of $\frakk$-finite vectors
 of $L^2(\Xi)$. 
We know that $L^2(\Xi)_{\frakk}=\calP(\Xi)e^{-\tr(-)}$,
 where $\calP(\Xi)$ denotes the space of restrictions
 of polynomials on $V$ to $\Xi$. 
\begin{proposition}
The formula
\begin{align*}
 \calT\psi(x) &:= \int_\Xi{F(x|y)\psi(y)\td\mu(y)}
\end{align*}
defines an operator $L^2(\Xi)_{\frakk}\to C^\infty(\Xi)$.
\end{proposition}

\begin{proof}
Use the integral formula \eqref{eq:IntFormulaO} and the asymptotic behaviour of the J-Bessel function $\widetilde{J}_\alpha(z)$ (see Appendix \ref{app:BesselFcts}) to show that the integral converges uniformly for $x$ in a bounded subset and $\psi(x)=p(x)e^{-\tr(x)}$, $p\in\calP(\Xi)$.
\end{proof}

\begin{proposition}
The operator $\calT$ extends
 to a unitary operator $\calT:L^2(\Xi)\to L^2(\Xi)$
 with $\calT\psi_0=\psi_0$.  
Further,
 $\calT$ leaves $L^2(\Xi)_{\frakk}$ invariant
 and intertwines the $\frakg$-action
 with the $\frakg$-action composed with
 $\Ad(\widetilde{j}):\frakg\to\frakg$.
\end{proposition}

\begin{proof}
By Proposition \ref{prop:BesselHypergeomEq} the operator $\calT$ intertwines the Bessel operator $\calB$ with the coordinate multiplication $-x$. 
Since both actions together generate the $\frakg$-action the intertwining property follows. Further, for $x\in\Xi\subseteq\XX$ we find
\begin{align*}
 \calT\psi_0(x) &= \int_\Xi{F(x|y)e^{-\tr(y)}\td\mu(y)}\\
 &= \int_\Xi{B(-(2x|y))e^{-2\tr(y)}\td\mu(y)}\\
 &= e^{-\tr(x)}\BB_\Xi\psi_0(-2x)\\
 &= e^{-\tr(x)}.
\end{align*}
Since $L^2(\Xi)_{\frakk}=\td\pi(\calU(\frakg))\psi_0$,
 it follows
 that $\calT$ maps $L^2(\Xi)_{\frakk}$
 into $L^2(\Xi)_{\frakk}$. 
Now, since invariant Hermitian forms on $L^2(\Xi)_{\frakk}$
 are unique up to a scalar,
 we find that $\calT$ is a unitary isomorphism.
\end{proof}

\begin{theorem}\label{thm:UnitaryInversionKernel}
$\calF_\Xi=\calT$.
\end{theorem}

\begin{proof}
By the previous proposition $\calF_\Xi\circ\calT^{-1}$ extends to a unitary isomorphism $L^2(\Xi)\to L^2(\Xi)$ which intertwines the $\frakg$-action. 
Therefore by Schur's Lemma, $\calF_\Xi$ is a scalar multiple of $\calT$. Since $\calF_\Xi\psi_0=\psi_0=\calT\psi_0$ this gives the claim.
\end{proof}

\begin{remark}
Since the group $G$ is generated by $j$ and a maximal parabolic subgroup whose action in the $L^2$-model is simple, the action of $j$ in the $L^2$-model is of special interest. For $V=\Sym(k,\RR)$, i.e. $\frakg=\sp(k,\RR)$, the operator $\calF_\Xi$ is basically the Euclidean Fourier transform, whereas for $V=\RR^{1,k-1}$, i.e. $\frakg=\so(2,k)$, the integral kernel of $\calF_\Xi$ was first calculated
 in Kobayashi--Mano \cite{KM05}, 
 and generalized to the non-Euclidean case
 $\RR^{p-1,q-1}$
 in \cite{KMa11}. 
The integral kernel of $\calF_{\Xi}$
 may involve distributions
 with singular support
 for non-Hermitian group $O(p,q)$
 ($p,q \ge 2$).  
For the case $V=\RR$, i.e. $\frakg=\sl(2,\RR)$, the operator $\calF_\Xi$ depends on the parameter $\lambda\in(0,\infty)$ and is the Hankel transform studied in B.~Kostant \cite{Kos00}.
\end{remark}

\begin{remark}
Since the functions $x\mapsto F(x|y)$, $y\in\Xi$, are eigenfunctions of the Bessel operator, the unitary inversion operator gives an expansion of any function $\psi\in L^2(\Xi)$ into eigenfunctions of the Bessel operator.
\end{remark}

Define $(-1)^*$ on $\calF(\XX)$ by $(-1)^*F(z)=F(-z)$.

\begin{proposition}\label{prop:InversionOnFockSpace}
$$ \BB_\Xi\circ\calF_\Xi = (-1)^*\circ\BB_\Xi. $$
\end{proposition}

\begin{proof}
We have $\td\rho(t({\bf e},0,-{\bf e}))=\td\pi_\CC(2t\sqrt{-1}(0,\1,0))=2t\sqrt{-1}(D_z+\frac{r\lambda}{2})$. 
Therefore we obtain
\begin{align*}
 \rho(e^{t({\bf e},0,-{\bf e})})F(z) &= e^{r\lambda\sqrt{-1}t}F(e^{2t\sqrt{-1}}z).
\end{align*}
For $t=\frac{\pi}{2}$ we obtain the action of $\widetilde{j}$ which is given by $e^{\pi\sqrt{-1}\frac{r\lambda}{2}}(-1)^*$.
\end{proof}

\begin{proposition}
$$ \calF_\Xi h_\alpha = (-1)^{|\alpha|}h_\alpha. $$
\end{proposition}

\begin{proof}
Since $\BB_\Xi h_\alpha=z^\alpha$ and $(-z)^\alpha=(-1)^{|\alpha|}z^\alpha$ the claim follows.
\end{proof}

\begin{theorem}[Bochner type identity]
For any $p\in\calH^m(\SS)$ we have
\begin{align*}
 \calF_\Xi(pe^{-\tr(x)}) &= e^{m\pi\sqrt{-1}}pe^{-\tr(x)}.
\end{align*}
\end{theorem}

\begin{proof}
Write $\calF_\Xi=e^{-r\lambda\frac{\pi}{2}\sqrt{-1}}e^{\frac{\pi}{2}\sqrt{-1}(\tr(x)-\calBe)}$. For $p\in\calH^m(\SS)$ we calculate with Lemma \ref{lem:BesselProdRule}:
\begin{align*}
 \calBe(pe^{-\tr(x)}) &= \calBe pe^{-\tr(x)}+2\left(\left.P\left(\frac{\partial p}{\partial x},\frac{\partial e^{-\tr(x)}}{\partial x}\right)\right|{\bf e}\right)+p\calBe e^{-\tr(x)}\\
\intertext{and since $\calBe p=0$, $\frac{\partial e^{-\tr(x)}}{\partial x}=-e^{-\tr(x)}{\bf e}$ and $\calBe e^{-\tr(x)}=(\tr(x)-r\lambda)e^{-\tr(x)}$, we obtain}
 &= -2\left(x\left|\frac{\partial p}{\partial x}\right.\right)+(\tr(x)-r\lambda)pe^{-\tr(x)}\\
 &= (\tr(x)-r\lambda-2m)pe^{-\tr(x)},
\end{align*}
because $(x|\frac{\partial}{\partial x})=\calE$ is the Euler operator which acts on $\calP_m(\SS)$ by the scalar $m$. Now it follows that
\begin{align*}
 (\tr(x)-\calBe)(pe^{-\tr(x)}) &= (r\lambda+2m)pe^{-\tr(x)}.
\end{align*}
Exponentiating this gives the claim.
\end{proof}
\section{Heat kernel and Segal--Bargmann transform}
\label{sec:heat}

We recall from Lemma \ref{lem:BesseleElliptic}
 that the second order differential operator
 $\calBe=({\bf e}|\calB)$ is an elliptic,
 self-adjoint operator on $L^2(\Xi, \td \mu)$.  
In this section,
 we consider the corresponding heat equation
\begin{align}
 (\calBe-\partial_t)u &= 0.\label{eq:HeatEquation}
\end{align}
We find the heat kernel \eqref{eqn:heatker}
 to the equation \eqref{eq:HeatEquation},
 and see that the Segal--Bargmann transform 
 can be obtained also
 by using the \lq{restriction principle}\rq.

\subsection{The heat equation and the heat kernel}

Recall the function $B(z|w)$, $z,w\in\XX$, from \eqref{eqn:Bker} which occurs in the kernel of the Segal--Bargmann transform and the unitary inversion operator. The following reproducing property will be needed later:

\begin{lemma}
\label{lem:5.3}
For $x\in\Xi$ and $z\in\XX$ the following identity holds
\begin{align*}
 2^{-r\lambda}\int_\Xi{e^{-\tr(\xi)}B(z|\xi)B(-(x|\xi))\td\mu(\xi)} &= e^{\tr(z)-\tr(x)}B(-(x|z)).
\end{align*}
\end{lemma}

\begin{proof}
For any $\psi\in C_c^\infty(\Xi)$ and $z\in\XX$ we have
\begin{align*}
 \BB_\Xi\calF_\Xi\psi(z) &= e^{-\frac{1}{2}\tr(z)}\int_\Xi{B(x|z)e^{-\tr(x)}\calF_\Xi\psi(x)\td\mu(x)}\\
 &= 2^{-r\lambda}e^{-\frac{1}{2}\tr(z)}\int_\Xi{\int_\Xi{B(x|z)B(-(x|y))e^{-\tr(x)}\psi(y)\td\mu(y)}\td\mu(x)}\\
 &= 2^{-r\lambda}e^{-\frac{1}{2}\tr(z)}\int_\Xi{\int_\Xi{e^{-\tr(x)}B(x|z)B(-(x|y))\td\mu(x)}\psi(y)\td\mu(y)}.
\end{align*}
On the other hand, by Proposition \ref{prop:InversionOnFockSpace} we obtain
\begin{align*}
 \BB_\Xi\calF_\Xi\psi(z) &= \BB_\Xi\psi(-z)\\
 &= e^{\frac{1}{2}\tr(z)}\int_\Xi{B(-(y|z))e^{-\tr(y)}\psi(y)\td\mu(y)}.
\end{align*}
Therefore the integral kernels have to coincide, which gives
\begin{align*}
 2^{-r\lambda}e^{-\frac{1}{2}\tr(z)}\int_\Xi{e^{-\tr(x)}B(x|z)B(-(x|y))\td\mu(x)} &= e^{\frac{1}{2}\tr(z)}B(-(y|z))e^{-\tr(y)}.
\end{align*}
This is the claimed formula.
\end{proof}

We define 
\begin{align}
\label{eqn:heatker}
 \Gamma(t,x,y) &:= (2t)^{-r\lambda}e^{-\frac{1}{t}(\tr(x)+\tr(y))}B\left(\left.\frac{x}{t}\right|\frac{y}{t}\right), & t>0,x,y\in\Xi.
\end{align}
Note that $\Gamma(t,x,y)>0$ for $t>0$ and $x,y\in\Xi$.
We now show that $\Gamma(t,x,y)$ is the heat kernel to the heat equation \eqref{eq:HeatEquation}.

\begin{theorem}\label{thm:HeatKernel}
The kernel $\Gamma(t,x,y)$ has the following properties:
\begin{enumerate}
\item[\textup{(1)}] $ \Gamma(t,x,y) = 2^{-2r\lambda}\int_\Xi{e^{-t\cdot\tr(\xi)}B(-(x|\xi))B(-(y|\xi))\td\mu(\xi)}. $
\item[\textup{(2)}] $ \int_\Xi{\Gamma(t,x,y)\td\mu(y)} = 1. $
\item[\textup{(3)}] $ \int_\Xi{\Gamma(s,x,z)\Gamma(t,y,z)\td\mu(z)} = \Gamma(s+t,x,y). $
\item[\textup{(4)}] For every $y\in\Xi$ the function $\Gamma(t,x,y)$ solves the heat equation \eqref{eq:HeatEquation}
\end{enumerate}
\end{theorem}

The proof is standard
(e.g. \cite{Roe03} for the Dunkl--Laplacian).
For the sake of completeness,
 we give a proof.  

\begin{proof}
\begin{enumerate}
\item[\textup{(1)}] 
This is immediate from Lemma \ref{lem:5.3}.
\item[\textup{(2)}] We have
\begin{align*}
 \int_\Xi{\Gamma(t,x,y)\td\mu(y)} &= (2t)^{-r\lambda}e^{-\frac{1}{t}\tr(x)}\int_\Xi{B\left(\left.\frac{x}{t}\right|\frac{y}{t}\right)e^{-\frac{1}{t}\tr(y)}\td\mu(y)}\\
\intertext{and substituting $2z=\frac{y}{t}$ we obtain}
 &= e^{-\frac{1}{t}\tr(x)}\int_\Xi{B\left(\left.\frac{2x}{t}\right|z\right)e^{-\tr(z)}\td\mu(z)}\\
 &= \BB_\Xi\psi_0\left(\frac{2x}{t}\right) = 1.
\end{align*}
\item[\textup{(3)}] We substitute (1) for the first factor in the integrand. This yields
\begin{align*}
 & \int_\Xi{\Gamma(s,x,z)\Gamma(t,y,z)\td\mu(z)}\\
 ={}& (8t)^{-r\lambda}\int_\Xi{\int_\Xi{e^{-s\tr(\xi)}B(-(x|\xi))B(-(z|\xi))e^{-\frac{1}{t}(\tr(y)+\tr(z))}B\left(\left.\frac{y}{t}\right|\frac{z}{t}\right)\td\mu(\xi)}\td\mu(z)}\\
\intertext{and substituting $z=t\eta$ gives}
 ={}& 8^{-r\lambda}\int_\Xi{\left(\int_\Xi{e^{-\tr(\eta)}B(-(\eta|t\xi))B\left(\left.\frac{y}{t}\right|\eta\right)\td\mu(\eta)}\right)e^{-\frac{1}{t}\tr(y)}e^{-s\tr(\xi)}B(-(x|\xi))\td\mu(\xi)}.\\
\intertext{Now Lemma \ref{lem:5.3} gives}
 ={}& 2^{-2r\lambda}\int_\Xi{e^{-(s+t)\tr(\xi)}B(-(x|\xi))B(-(y|\xi))\td\mu(\xi)}\\
 ={}& \Gamma(s+t,x,y)
\end{align*}
by (1) again.
\item[\textup{(4)}] This follows from (1) by differentiating under the integral.\qedhere
\end{enumerate}
\end{proof}

The kernel $\Gamma(t,x,y)$ can be used to construct solutions to the heat equation \eqref{eq:HeatEquation}. In fact, the heat semigroup $e^{t\calBe}$, $t\geq0$, is explicitly given in terms of the integral kernel $\Gamma(t,x,y)$ as follows:
\begin{align*}
 e^{t\calBe}f(x) &= \int_\Xi{\Gamma(t,x,y)f(y)\td\mu(y)}.
\end{align*}
Using this observation we now interpret the Segal--Bargmann transform purely in terms of the $\sl_2$-triple $(E,F,H)$ which was introduced in \eqref{eq:sl2triple}.

\begin{theorem}
$$ \BB_\Xi = 2^{r\lambda}e^{\frac{1}{2}\tr}e^{\calBe} = 2^{\frac{rd}{2}}e^{-\frac{1}{2}\sqrt{-1}\td\pi(E)}e^{-\sqrt{-1}\td\pi(F)}. $$
\end{theorem}

\begin{proof}
It is clear by definition that
\begin{align*}
 \BB_\Xi f(z) &= 2^{r\lambda}e^{\frac{1}{2}\tr(z)}\int_\Xi{\Gamma(1,z,x)f(x)\td\mu(x)}.
\end{align*}
Since $\td\pi(E)=\sqrt{-1}\tr(x)$ and $\td\pi(F)=\sqrt{-1}\calBe$ the claimed formula holds.
\end{proof}

\subsection{The Segal--Bargmann transform with the heat kernel}
\label{subsec:5.2}

The formula $\calR_\Xi F(x)=e^{-\frac{1}{2}\tr(x)}F(x)$ defines an operator $\calP(\XX)\to L^2(\Xi)$ and hence we obtain a densely defined unbounded operator $\calR_\Xi:\calF(\XX)\to L^2(\Xi)$. 
Therefore it makes sense to consider its adjoint $\calR_\Xi^*:L^2(\Xi)\to\calF(\XX)$ as a densely defined unbounded operator.

\begin{proposition}
\label{prop:RRconti}
For $f\in L^2(\Xi)$ we have
\begin{align*}
 \calR_\Xi\calR_\Xi^*f(x) &= 2^{2r\lambda}\int_\Xi{\Gamma(2,x,y)f(y)\td\mu(y)}
\end{align*}
and $\calR_\Xi\calR_\Xi^*f\in L^2(\Xi)$. This defines a continuous linear operator with operator norm $\|\calR_\Xi\calR_\Xi^*\|\leq2^{2r\lambda}$.
\end{proposition}

\begin{proof}
We have
\begin{align*}
 \calR_\Xi^*f(z) &= \langle\calR_\Xi^*f|\KK_z\rangle\\
 &= \langle f|\calR_\Xi\KK_z\rangle\\
 &= \int_\Xi{f(y)\overline{\calR_\Xi\KK_z(y)}\td\mu(y)}\\ 
 &= \int_\Xi{B\left(\left.\frac{y}{2}\right|\frac{z}{2}\right)e^{-\frac{1}{2}\tr(y)}f(y)\td\mu(y)}\\
 &= 2^{2r\lambda}e^{\frac{1}{2}\tr(z)}\int_\Xi{\Gamma(2,z,y)f(y)\td(y)}
\end{align*}
and the formula follows. Now, H\"older's inequality gives
\begin{align*}
 \int_\Xi{|\Gamma(2,x,y)f(y)|\td\mu(y)} &\leq \left(\int_\Xi{\Gamma(2,x,y)\td\mu(y)}\right)^{\frac{1}{2}}\left(\int_\Xi{\Gamma(2,x,y)|f(y)|^2\td\mu(y)}\right)^{\frac{1}{2}}\\
 &= \left(\int_\Xi{\Gamma(2,x,y)|f(y)|^2\td\mu(y)}\right)^{\frac{1}{2}}.
\end{align*}
where we have used Theorem \ref{thm:HeatKernel}\,(2). Then we find, using Fubini's theorem:
\begin{align*}
 \|\calR_\Xi\calR_\Xi^*f\|^2 &= 2^{4r\lambda}\int_\Xi{\left|\int_\Xi{|\Gamma(2,x,y)f(y)\td\mu(y)}\right|^2\td\mu(x)}\\
 &\leq 2^{4r\lambda}\int_\Xi{\int_\Xi{\Gamma(2,x,y)|f(y)|^2\td\mu(y)}\td\mu(x)}\\
 &= 2^{4r\lambda}\int_\Xi{|f(y)|^2\td\mu(y)} = 2^{4r\lambda}\|f\|^2
\end{align*}
and the proof is complete.
\end{proof}

By Proposition \ref{prop:RRconti}, 
the operator $\calR_\Xi\calR_\Xi^*$
 is a continuous, 
 positive operator.  
Hence the operator $|\calR_\Xi|:=\sqrt{\calR_\Xi\calR_\Xi^*}$
 is well-defined. 
We now show that the Segal--Bargmann transform can be constructed only from the restriction map $\calR_\Xi$.

\begin{proposition}\label{prop:RestrictionPrinciple}
$\calR_\Xi^*=\BB_\Xi\circ\sqrt{\calR_\Xi\calR_\Xi^*}$.
\end{proposition}

\begin{proof}
The previous proposition and the properties of the heat kernel yield
\begin{align*}
 |\calR_\Xi|f(x) &= \sqrt{\calR_\Xi\calR_\Xi^*}f(x)\\
 &= 2^{r\lambda}\int_\Xi{\Gamma(1,x,y)f(y)\td\mu(y)}.
\end{align*}
Hence, we obtain
\begin{align*}
 \BB_\Xi(|\calR_\Xi|f)(z) &= \int_\Xi{e^{-\frac{1}{2}\tr(z)}B(z,y)e^{-\tr(y)}|\calR_\Xi|f(y)\td\mu(y)}\\
 &= 2^{2r\lambda}e^{\frac{1}{2}\tr(z)}\int_\Xi{\int_\Xi{\Gamma(1,z,y)\Gamma(1,y,x)f(x)\td\mu(x)}\td\mu(y)}\\
 &= 2^{2r\lambda}e^{\frac{1}{2}\tr(z)}\int_\Xi{\Gamma(2,z,x)f(x)\td(x)}\\
 &= \calR_\Xi^*f(z).\qedhere
\end{align*}
\end{proof}
\section{Example: $\frakg=\so(2,n)$}
\label{sec:exso2n}

We study the example $\frakg=\so(2,n)$ in more detail
 and discuss also the relation
 with the results in \cite{KM05, KM07a}.

\subsection{The Schr{\"o}dinger model}

Let $V=\RR^{1,n-1}$ be the Euclidean Jordan algebra with multiplication
\begin{align*}
 x \cdot y &= (x_1y_1+x'\cdot y',x_1y'+y_1x'), 
\end{align*}
 for $x=(x_1,x')$, $y=(y_1,y') \in \RR^{1,n-1}=\RR \oplus \RR^{n-1}$.  
The unit element is given by ${\bf e}=(1,0,\ldots,0)$. 
Then $L=\Str(V)_0 \simeq \RR_+ \times \SO(1,n-1)_0$
 and $G=\Co(V)_0$ is the adjoint group
 of ${\mathfrak {so}}(2,n)$. 

Take the primitive idempotent
 $c_1= \frac 1 2 (1,0,\cdots, 0,1)$.  
The $L$-orbit $\Xi$ through the primitive element $c_1$
 in $V$ is given by the future light cone
 in the Minkowski space $\RR^{1,n-1}$:
\begin{align*}
 \Xi &= \{x\in\RR^{1,n-1}:x_1=\sqrt{x_2^2+\cdots+x_n^2}>0\}.
\end{align*}
The trace form of $\RR^{1,n-1}$ as a Jordan algebra
 takes the form
\[
  (x|y)=2(x_1 y_1 + (x',y'))
       =4 (x',y')
\]
on $\Xi$, 
where $(\,,\,)$ denotes the standard inner product
 on $\RR^{n-1}$.  
Since the volume
 of the Euclidean sphere in $\RR^{n-1}$
 of radius $\frac 1 2$
 is given by 
\[
  \frac{1}{2^{n-2}} \frac{2 \pi^{\frac {n-1}{2}}}{\Gamma(\frac {n-1}{2})},
\]
our normalization 
of the measure $\td \mu$
 (see \eqref{eq:IntFormulaO})
 on the orbit $\Xi$ is given by
\begin{align}
\td \mu 
=& 
\frac{2^{n-2}}{\Gamma(n-2)} 
\frac{2^{n-3} \Gamma(\frac{n-1}{2})}{\pi^{\frac {n-1}{2}}}
r^{n-3} \td r \td \omega
\label{eqn:someasure}
\\
=&
\frac{2^{n-2}}{\Gamma(\frac{n-2}2) \pi^{\frac {n-2}{2}}}
r^{n-3} \td r \td \omega
\notag
\end{align}
in polar coordinates 
$
     \RR_+\times S^{n-2}\to\Xi,\,
     (r,\omega)\mapsto(r,r\omega).  
$ 
We set 
\begin{align*}
 \varepsilon_i &= \begin{cases}+1 & \mbox{for $i=1$,}\\-1 & \mbox{for $2\leq i\leq n$,}\end{cases} & \Box &= \sum_{i=1}^n{\varepsilon_i\frac{\partial^2}{\partial x_i^2}}, & E &= \sum_{i=1}^n{x_i\frac{\partial}{\partial x_i}}.
\end{align*}
Then the Bessel operator is of the form (see \cite[Proposition 2.36]{HKM11})
\begin{align*}
 \calB &= -\frac{1}{4}\sum_{i=1}^n{\calB_ie_i}
\end{align*}
with
\begin{align*}
 \calB_i &= \varepsilon_ix_i\Box-(2E+n-2)\frac{\partial}{\partial x_i},
\end{align*}
which are exactly the {\it{fundamental differential operators}}
 $P_i$ 
 on the isotypic cone
 introduced in \cite[(1.1.3)]{KMa11}
 with the signature $(n_1,n_2)=(1,n-1)$
 in our setting here.  

The unitary operator 
 in the Schr{\"o}dinger model
 on $L^2(\RR^{n-1}, \frac {\td x}{|x|})$
 corresponding to the inversion element
 $w_0=\exp(\frac \pi 2\begin{pmatrix} 0 & -1 \\ 1 & 0\end{pmatrix})$
 was obtained previously by Kobayashi--Mano
 as the integral transform
 against the following kernel function:
\begin{equation}
\label{eqn:KMinv}
 \frac{\widetilde{J}_{\frac{n-4}{2}}(4 \sqrt{(x',y')})}
      {e^{\frac{n-2}{2}\pi i}\pi^{\frac{n-2}{2}}}
  \frac{d y'}{|y'|}.  
\end{equation}
See \cite[Theorem D]{KM05}
 as a special value of the holomorphic semigroup,
 or alternatively
 as a special case
 of the indefinite orthogonal group
 $O(n_1+1,n_2+1)$
with $(n_1, n_2)=(1,n-1)$
 in \cite[Theorem 1.3.1]{KMa11}.  

In view of $\widetilde {j}=w_0^{-1}$,
 our Fourier inversion operator $\calF_{\Xi}$
 takes the form 
\[
\calF_{\Xi}
=e^{-\frac{\pi i (n-2)}{2}}\pi(\widetilde j)
=e^{\frac{\pi i (n-2)}{2}}\pi(w_0).  
\]
Hence the kernel function \eqref{eqn:KMinv}
 gives (of course)
 the same formula
 of the Fourier inversion operator
 ${\mathcal{F}}_{\Xi}$
 in Theorem \ref{thm:D}
 because
\begin{align*}
&2^{-r \lambda}\Gamma(\lambda)
 \widetilde J_{\lambda -1}(2\sqrt{(x|y)})\td \mu(y)
\\
=&2^{-(n-2)}\Gamma(\frac{n-2}{2})
  \widetilde J_{\frac{n-4}{2}}(2\sqrt{2(x',y')})\td \mu(y)
\\
=&\frac{1}{\pi^{\frac{n-2}{2}}} 
  \widetilde J_{\frac{n-4}{2}}(2\sqrt{2(x',y')})\frac{\td y}{|y'|}
\end{align*}
by \eqref{eqn:someasure}.  

\subsection{The Fock model}

The complex orbit $\XX$ through the primitive idempotent $c$
 is given by
\begin{align*}
 \XX &= \{z\in\CC^n:z_1^2=z_2^2+\cdots+z_n^2\}\setminus\{0\}, 
\end{align*}
which contains $\Xi$
 as a totally real submanifold. Put
\begin{align*}
 \SS_\XX := \{z\in\XX:|z|=\sqrt{2(|z_1|^2+\cdots+|z_n|^2)}=1\}
\end{align*}
and note that this is a compact symmetric space for the group $K^{L_\CC}\simeq\SO(n)\times U(1)$.
In view of the polar decomposition
\begin{align*}
 \RR_+\times\SS_\XX\stackrel{\sim}{\to}\Xi,\,(r,\eta)\mapsto r\eta,
\end{align*}
the measure $\td\nu$ on $\Xi$ takes the form (cf. \eqref{eq:IntFormulaX})
\begin{align*}
 \td\nu = \frac{1}{2^{2n-6}\Gamma(n-2)\Gamma(\frac{n}{2})}r^{2n-5}\td r\td\eta
\end{align*}
where $\td\eta$ denotes the unique $K^{L_\CC}$-invariant measure on $\SS_\XX$ with total volume $1$. Therefore the inner product on the Fock space $\calF(\XX)$ is given by
\begin{align*}
 \langle F,G\rangle &= \frac{1}{2^{2n-6}\Gamma(n-2)\Gamma(\frac{n}{2})}\int_{\SS_\XX}{\int_0^\infty{F(r\eta)\overline{G(r\eta)}r^{2n-5}\td r}\td\eta}, & F,G\in\calF(\XX).
\end{align*}
The Segal--Bargmann transform in Theorem 
 \ref{thm:c}
 amounts to 
\[
({\mathbb{B}}_{\Xi}f)(z)
=\frac{2^{n-2}}{\pi^{\frac{n-2}{2}}}
\exp (-(z',z')^{\frac 1 2})
\int_{0}^{\infty} \int_{S^{n-2}} 
\widetilde I_{\frac{n-4}{2}}(4 r^{\frac 1 2}(z',\omega)^{\frac 12})
e^{-2r}
f(r \omega)
r^{n-3}
\td r \td \omega, 
\]
for $z =(z_1,z') \in \XX$, 
 and $f \in L^2(\Xi)$.  

The ring of regular functions
 on $\XX$
 is given by 
\begin{align*}
 \calP(\XX) &= \CC[z_1, \cdots, x_n]/\langle z_1^2-z_2^2-\cdots-z_n^2\rangle,
\end{align*}
where $\langle z_1^2-z_2^2-\cdots-z_n^2\rangle$ denotes the ideal generated by $z_1^2-z_2^2-\cdots-z_n^2$. 
Hence, 
the $\frakk$-types $\calP^m(\XX)$ are given by
\begin{align*}
 \calP^m(\XX) &= \CC_m[\CC^{n-1}]\oplus z_1\CC_{m-1}[\CC^{n-1}]
\end{align*}
and hence
\begin{align*}
 \dim\calP^m(\XX) &= \dim\CC_m[\CC^{n-1}]+\dim\CC_{m-1}[\CC^{n-1}]\\
 &= {n+m-2\choose m}+{n+m-3\choose m-1}\\
 &= {n+m-1\choose n}-{n+m-3\choose n-1}\\
 &= \dim\calH^m(\RR^n).
\end{align*}
In fact, an $\so(n)$-equivariant isomorphism $\Phi:\calP^m(\XX)\stackrel{\sim}{\to}\calH^m(\RR^n)$ can be constructed as follows: Let $p\in\calP^m(V_\CC)$, then the polynomial $p(\sqrt{-1}z_1,z_2,\ldots,z_n)$ has an expansion into classical spherical harmonics
\begin{align*}
 p(\sqrt{-1}z_1,z_2,\ldots,z_n) &= \sum_{k=0}^{\lfloor\frac{m}{2}\rfloor}{h_{m-2k}(z)(z_1^2+\cdots+z_n^2)^k},
\end{align*}
where $h_{m-2k}\in\calH^{m-2k}(\RR^n)$. Then put $\Phi(p|_\XX):=h_m$. This map is well-defined since $-z_1^2+z_2^2+\cdots+z_n^2$ vanishes on $\XX$. It shows in particular that every polynomial $p\in\calP(\XX)$ has a unique extension to a polynomial on $V_\CC$ in the kernel of the differential operator $-\frac{\partial^2}{\partial z_1^2}+\frac{\partial^2}{\partial z_2^2}+\cdots+\frac{\partial^2}{\partial z_n^2}$.

\appendix
\renewcommand{\theequation}{A.\arabic{equation}} 

\section{Appendix: Special Functions}

\subsection{Renormalized Bessel functions}\label{app:BesselFcts}

Following \cite{KM07a, KMa11},
 we renormalize the Bessel functions
 $J_\alpha(z)$, $I_\alpha(z)$ and $K_\alpha(z)$
 by:
\begin{align*}
 \widetilde{J}_\alpha(z) &= \left(\frac{z}{2}\right)^{-\alpha}J_\alpha(z),\\
 \widetilde{I}_\alpha(z) &= \left(\frac{z}{2}\right)^{-\alpha}I_\alpha(z),\\
 \widetilde{K}_\alpha(z) &= \left(\frac{z}{2}\right)^{-\alpha}K_\alpha(z).
\end{align*}
In the analysis
 of minimal representations,
 these functions
 appear naturally rather
 than the usual Bessel functions.  
We refer the reader to 
\cite[\S7.2]{KMa11}
 for a concise summary
 of the renormalized Bessel functions.   
Among others, 
 $\widetilde{J}_\alpha(z)$
 and $\widetilde{I}_\alpha(z)$ are entire functions, 
$\widetilde{J}_\alpha(\sqrt{-1}z)=\widetilde{I}_\alpha(z)$
 and $\widetilde{J}_\alpha(-z)=\widetilde{J}_\alpha(z)$, $\widetilde{I}_\alpha(-z)=\widetilde{I}_\alpha(z)$. 
In particular, 
 $\widetilde{J}_\alpha(\sqrt{z})$ and $\widetilde{I}_\alpha(\sqrt{z})$ are entire functions. 
Their Taylor expansions are given by 
\begin{align*}
 \widetilde{J}_\alpha(2\sqrt{z}) &= \sum_{n=0}^\infty{\frac{(-1)^n}{\Gamma(n+\alpha+1)n!}z^n},\\
 \widetilde{I}_\alpha(2\sqrt{z}) &= \sum_{n=0}^\infty{\frac{1}{\Gamma(n+\alpha+1)n!}z^n}.
\end{align*}
The function $\widetilde{J}_\alpha(z)$ solves the differential equation
\begin{align*}
 zu''+(2\alpha+1)u'+zu &= 0,
\end{align*}
whereas the functions $\widetilde{I}_\alpha(z)$ and $\widetilde{K}_\alpha(z)$ are linear independent solutions to the differential equation
\begin{align*}
 zu''+(2\alpha+1)u'-zu &= 0.
\end{align*}

The renormalized I-Bessel function (and also the corresponding J-Bessel function) has the following asymptotic behaviour:
\begin{align*}
 \widetilde{I}_\alpha(0) &= \frac{1}{\Gamma(\alpha+1)},\\
 |\widetilde{I}_\alpha(z)| &\lesssim |z|^{-\alpha-\frac{1}{2}}e^{|z|} & \mbox{as $|z|\to\infty$.}
\end{align*}

The asymptotic behaviour of the K-Bessel function is given by
\begin{align*}
 \widetilde{K}_\alpha(x) &= \left\{\begin{array}{ll}\frac{\Gamma(\alpha)}{2}\left(\frac{x}{2}\right)^{-2\alpha}+o(x^{-2\alpha}) &\mbox{if $\alpha>0$}\\-\log(\frac{x}{2})+o(\log(\frac{x}{2})) &\mbox{if $\alpha=0$}\\\frac{\Gamma(-\alpha)}{2}+o(1) &\mbox{if $\alpha<0$}\end{array}\right. & \mbox{as }x &\rightarrow 0,\\
 \widetilde{K}_\alpha(x) &= \frac{\sqrt{\pi}}{2}\left(\frac{x}{2}\right)^{-\alpha-\frac{1}{2}}e^{-x}\left(1+\mathcal{O}\left(\frac{1}{x}\right)\right) & \mbox{as }x &\rightarrow\infty.
\end{align*}

We further have the following integral formula for $\Re(\beta+1),\Re(\beta-2\alpha+1)>0$ and $\Re(a)>0$ (see \cite[formula 6.561\,(16)]{GR65}):
\begin{align}
 \int_0^\infty{\widetilde{K}_\alpha(ax)x^\beta\td x} &= 2^{\beta-1}a^{-\beta-1}\Gamma\left(\frac{\beta+1}{2}\right)\Gamma\left(\frac{\beta-2\alpha+1}{2}\right).\label{eq:KBesselIntFormula}
\end{align}

\subsection{The Gau\ss\ hypergeometric function ${_2F_1}(a,b;c;z)$}

The Gau\ss\ hypergeometric function ${_2F_1}(a,b;c;z)$ is for $|z|<1$ defined by
\begin{align*}
 {_2F_1}(a,b;c;z) &= \sum_{n=0}^\infty{\frac{(a)_n(b)_n}{n!\,(c)_n}z^n}.
\end{align*}
If $a=-n\in-\NN$, then the series is finite and ${_2F_1}(-n,b;c;z)$ is a polynomial and hence an entire function in $z\in\CC$.

The Gau\ss\ hypergeometric function ${_2F_1}(a,b;c;z)$ solves the following second order ordinary differential equation:
\begin{align*}
 (1-z)zu''(z)+(c-(a+b+1)z)u'(z)-abu(z) &= 0.
\end{align*}

For $c=\frac{1}{2}$ and $a=-n\in-\NN$ we have by \cite[equations (6.3.5), (6.4.23), (3.1.1), (6.4.9)]{AAR99}
\begin{align}
 {_2F_1}(-n,b;\textstyle{\frac{1}{2}};z^2) &= \frac{n!}{(\frac{1}{2})_n}P_n^{(-\frac{1}{2},b-n-\frac{1}{2})}(1-2z^2)\notag\\
 &= (-1)^n\frac{n!}{(\frac{1}{2})_n}P_n^{(b-n-\frac{1}{2},-\frac{1}{2})}(2z^2-1)\notag\\
 &= (-1)^n\frac{(2n)!}{(\frac{1}{2})_n(b+\frac{1}{2})_n}P_{2n}^{(b-n-\frac{1}{2},b-n-\frac{1}{2})}(z)\notag\\
 &= (-1)^n\frac{(2n)!\,(b-n+\frac{1}{2})_n}{(\frac{1}{2})_n(2b-2n)_n(b+\frac{1}{2})_n}C_{2n}^{b-n}(z),\label{eq:2F1asCnlambda}
\end{align}
where $P_{n}^{(\alpha,\beta)}$ is the Jacobi polynomial
 and $C_n^\lambda(z)$ denotes the Gegenbauer polynomial.
\begin{sidewaystable}[H]
\begin{center}
\begin{tabular}{|c|c|c|c|c|c|}
  \cline{1-6}
  $V$ & $\frakg=\co(V)$ & $\frakl=\str(V)$ & $\frakk^\frakl=\aut(V)$ & $\SS=K^L/K^L_{c_1}$ & $\lambda$\\
  \hline\hline
  $\RR$ & $\sl(2,\RR)$ & $\RR$ & $0$ & $\{\1\}$ & $\in(0,\infty)$\\
  $\Sym(k,\RR)$ ($k\geq2$) & $\sp(k,\RR)$ & $\sl(k,\RR)\oplus\RR$ & $\so(k)$ & $\PP^{k-1}(\RR)=\SO(k)/S(O(1)\times O(k-1))$ & $1/2$\\
  $\Herm(k,\CC)$ ($k\geq2$) & $\su(k,k)$ & $\sl(k,\CC)\oplus\RR$ & $\su(k)$ & $\PP^{k-1}(\CC)=\SU(k)/S(U(1)\times U(k-1))$ & $1$\\
  $\Herm(k,\HH)$ ($k\geq2$) & $\so^*(4k)$ & $\sl(k,\HH)\oplus\RR$ & $\sp(k)$ & $\PP^{k-1}(\HH)=\Sp(k)/(\Sp(1)\times\Sp(k-1))$ & $2$\\
  $\RR^{1,k-1}$ ($k\geq3$) & $\so(2,k)$ & $\so(1,k-1)\oplus\RR$ & $\so(k-1)$ & $S^{k-2}=\SO(k-1)/\SO(k-2)$ & $(k-2)/2$\\
  $\Herm(3,\OO)$ & $\mathfrak{e}_{7(-25)}$ & $\mathfrak{e}_{6(-26)}\oplus\RR$ & $\mathfrak{f}_4$ & $\PP^2(\OO)=F_4/\Spin(9)$ & $4$\\
  \hline
\end{tabular}
\caption{Simple Euclidean Jordan algebras and their corresponding Lie algebras \label{tb:Groups}}
\end{center}
\end{sidewaystable}

\bibliographystyle{amsplain}
\bibliography{bibdb}

\vspace{30pt}

\textsc{Joachim Hilgert\\Institut f\"ur Mathematik, Universit\"at Paderborn, Warburger Str. 100, 33098 Paderborn, Germany.}\\
\textit{E-mail address:} \texttt{hilgert@math.uni-paderborn.de}\\

\textsc{Toshiyuki Kobayashi\\Kavli IPMU and Graduate School of Mathematical Sciences, The University of Tokyo, 3-8-1 Komaba, Meguro, Tokyo, 153-8914, Japan.}\\
\textit{E-mail address:} \texttt{toshi@ms.u-tokyo.ac.jp}\\

\textsc{Jan M\"ollers\\Institut for Matematiske Fag, Aarhus Universitet, Ny Munkegade 118, Bygning 1530, Lokale 423, 8000 Aarhus C, Denmark.}\\
\textit{E-mail address:} \texttt{moellers@imf.au.dk}\\

\textsc{Bent {\O}rsted\\Institut for Matematiske Fag, Aarhus Universitet, Ny Munkegade 118, Bygning 1530, Lokale 431, 8000 Aarhus C, Denmark.}\\
\textit{E-mail address:} \texttt{orsted@imf.au.dk}

\end{document}